%% file: Splicing-V4.tex
\documentclass[10pt]{amsart}
\usepackage{amsmath,amsthm,amsfonts,amscd,flafter,epsf,amssymb,graphicx,color}
\usepackage{epsfig,manfnt}
\usepackage{amsfonts}
\usepackage{amssymb}
\usepackage{amsmath}
\usepackage{amsthm}
\usepackage{latexsym}
\usepackage{amscd}
\usepackage{flafter}
\usepackage{epsf}
\usepackage{tikz,tikz-cd}
\input{diagram}

\input{command1}

\input{standardcommand}

\begin{document}

\title[Floer homology and splicing knot complements]
{Floer homology and splicing knot complements}%
\author{Eaman Eftekhary}%
\address{School of Mathematics, Institute for Research in Fundamental Sciences (IPM),
P. O. Box 19395-5746, Tehran, Iran}%
\email{eaman@ipm.ir}



\begin{abstract}
We  obtain a formula for  the Heegaard Floer homology (hat theory) of
the three-manifold $Y(K_1,K_2)$ obtained by splicing the complements of the 
knots $K_i\subset Y_i$, $i=1,2$, in terms of the knot Floer
homology of $K_1$ and $K_2$. We also present a few applications. If $h_n^i$ denotes 
the rank of the Heegaard Floer group $\ov\HFKT$ for the knot obtained by $n$-surgery 
over $K_i$ we show that the rank of $\ov{\HFT}(Y(K_1,K_2))$ is bounded below by 
$$\big|(h_\infty^1-h_1^1)(h_\infty^2-h_1^2)-
(h_0^1-h_1^1)(h_0^2-h_1^2)\big|.$$
We also show that if splicing the complement of a knot $K\subset Y$ with the trefoil 
complements gives a homology sphere $L$-space then $K$ is trivial and $Y$ is a homology 
sphere $L$-space.
\end{abstract}
\maketitle
\tableofcontents
\section{Introduction}\label{sec:introduction}
Heegaard Floer homology, introduced by Ozsv\'ath and Szab\'o in \cite{OS-3mfld},
has been the source of powerful techniques for the study of objects in
low dimensional topology.  It is  interesting to investigate if  Heegaard Floer homology
is capable to distinguish the standard sphere from other homology spheres.
Since the Heegaard Floer groups of the connected sum of two homology spheres 
are obtained as the tensor product of the Heegaard Floer groups associated with 
the two pieces, the question is reduced to determining prime homology spheres 
with trivial Heegaard Floer groups. The Poincar\'e homology sphere 
$\Sig(2,3,5)$ with either orientation 
is the unique known example of a non-trivial prime 
homology sphere $Y$ with  $\ov{\HFT}(Y;\Z)=\Z$.
A conjecture of Ozsv\'ath and Szab\'o predicts that this is in fact the only possible 
example.
\\

In this paper we study the  Heegaard Floer groups of a homology sphere $Y$ 
which contains an incompressible torus. We may
use the incompressible torus to decompose $Y$, fill out 
the torus boundary of each one of the two pieces by gluing a solid torus 
and obtain two new homology  spheres $Y_1$ and $Y_2$. By requiring 
$Y_1$ and $Y_2$ to be homology spheres the gluing of the solid tori is
determined; the decomposition  determines a knot $K_i$  in $Y_i$, $i=1,2$ and
 $Y=Y(K_1,K_2)$ is obtained by splicing the complements of $K_1$ and 
$K_2$ in $Y_1$ and $Y_2$, respectively. A formula is obtained for $\ov{\HFT}(Y;\Fbb)$
 in terms of the knot Floer {{objects}} associated with $K_1\subset Y_1$ and 
 $K_2\subset Y_2$, where $\Fbb$ denotes the field $\Z/2\Z$ with two elements. \\

The more precise statement of the splicing formula obtained in this paper is as follows.
Let $K\subset Y$ denote a null-homologous knot inside a three-manifold $Y$.
For every $n\in\Z\cup\{\infty\}$ let $Y_n=Y_n(K)$ denote the three-manifold
 obtained by performing 
 $n$-surgery on $K$, and let $K_n\subset Y_n$
denote the knot in $Y_n$ which is the core of the neighbourhood replaced for
 nd$(K)\subset Y$ in constructing $Y_n$.
Denote the homology group $\widehat{\text{HFK}}(Y_n,K_n;\Fbb)$
by $\Hbb_n(K)$, and its dimension as a vector space over $\Fbb$ by $h_n(K)$. 
In particular, 
$\Hbb_\infty(K)=\widehat{\text{HFK}}(Y,K;\Fbb)$
and $\Hbb_0(K)=\widehat{\text{HFL}}(Y,K;\Fbb)$ are the knot Floer homology and 
the longitude Floer homology of $K$, respectively  (see \cite{OS-knot, Ef-Whitehead}).\\

\begin{figure}
\def\svgwidth{12cm}
\begin{center}
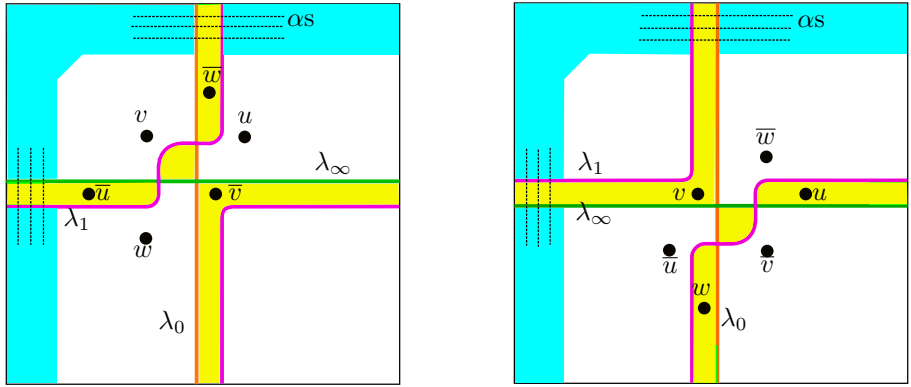
 \caption{\label{fig:punctures}\label{fig:general-triple}\
{The curves $\lambda_0$ (orange), $\lambda_1$ (pink) and $\lambda_\infty$ (green) 
and the punctures  are chosen following one of the above two patterns. 
 The  punctures $u,v$ and $w$ are used to define 
 $\pphi_0,\pphi_1$ and $\pphi_\infty$, while $\ubar,\vbar$ and $\wbar$ are used to 
 define $\pphibar_0,\pphibar_1$  and $\pphibar_\infty$.}}
\end{center}
\end{figure}

Choose a Heegaard diagram 
$$H=(\Sig,\alphas=\{\alpha_1,...,\alpha_g\},\ov\betas=\{\beta_1,...,\beta_{g-1}\})$$
for the knot complement  $Y\setminus K$, and let $\lambda_\bullet$ denote
an oriented  longitude which has framing coefficient $\bullet\in \Z\cup\{\infty\}$.
One can choose the curves $\lambda_\bullet$ (which are disjoint from the curves in 
$\ov\betas$) so that the pairs $(\lambda_0,\lambda_\infty)$, 
$(\lambda_1,\lambda_\infty)$ and
$(\lambda_0,\lambda_1)$ have  single intersection points in the Heegaard diagram.
For $\bullet\in\{0,1,\infty\}$ set
$$\betas_\bullet=\{\beta_1^\bullet,...,\beta_{g-1}^\bullet,\lambda_\bullet\},$$
where $\beta_i^\bullet$ is an isotopic copy of the curve $\beta_i$. The pictures on the 
left-hand-side and the right-hand-side of Figure~\ref{fig:punctures} illustrates 
two possible general arrangements for the curves $\lambda_0,\lambda_1$ and 
$\lambda_\infty$. 
In Figure~\ref{fig:punctures} and other figures in this paper the surface orientation  
is chosen opposite from the standard orientation of the page in order to stay compatible 
with the orientation convention of \cite{OS-3mfld}.\\
  
The two Heegaard quadruples
\begin{displaymath}
(\Sig,\alphas,\betas_0,\betas_1,\betas_\infty;u,v,w)\ \ \text{and}\ \ 
(\Sig,\alphas,\betas_0,\betas_1,\betas_\infty;\ubar,\vbar,\wbar)
\end{displaymath}
obtained in this way then correspond to the exact triangles
\begin{equation}\label{eq:exact-triangles}
\begin{diagram}
\Hbb_0(K)&&\lTo{\pphi_1(K)} & & \Hbb_\infty(K) &
&\Hbb_0(K)&&\lTo{\pphibar_1(K)} && \Hbb_\infty(K)\\
&\rdTo{\pphi_\infty(K)} &&\ruTo{\pphi_0(K)}&& 
\text{and} &&\rdTo{\pphibar_\infty(K)} &&\ruTo{\pphibar_0(K)}&\\
&&\Hbb_1(K)&&&&&&\Hbb_1(K)&&
\end{diagram}, 
\end{equation}
respectively. The ranks of both $\pphi_\bullet(K)$ and $\pphibar_\bullet(K)$ are equal to 
$$a_\bullet(K)=\frac{h_0(K)+h_1(K)+h_\infty(K)-2h_\bullet(K)}{2},\ \ \ \ \ 
\bullet\in\{0,1,\infty\}.$$ 
The exactness of the above two triangles imply that the induced maps
\begin{displaymath}
{\Coker(\pphi_0(K))}\lra {\Ker(\pphi_\infty(K))}\ \ \ \text{and}\ \ \ 
{\Coker(\pphibar_0(K))}\lra {\Ker(\pphibar_\infty(K))} 
\end{displaymath}
by $\pphi_1(K)$ and $\pphibar_1(K)$ are isomorphisms. 
Note that both the domain and the target of the aforementioned isomorphisms are 
of dimension $a_1(K)$. Take 
$\theta(K)\co\Hbb_0(K)\ra \Hbb_\infty(K)$ 
(respectively, $\thetabar(K)\co\Hbb_0(K)\ra \Hbb_\infty(K)$) to be
an arbitrary extension of the inverse of the isomorphism induced by 
$\pphi_1(K)$ (respectively, $\pphibar_1(K)$), so that the ranks of both $\theta(K)$ and 
$\thetabar(K)$ are equal to $a_1(K)$. \\

\begin{figure}
\def\svgwidth{13cm}
\begin{center}
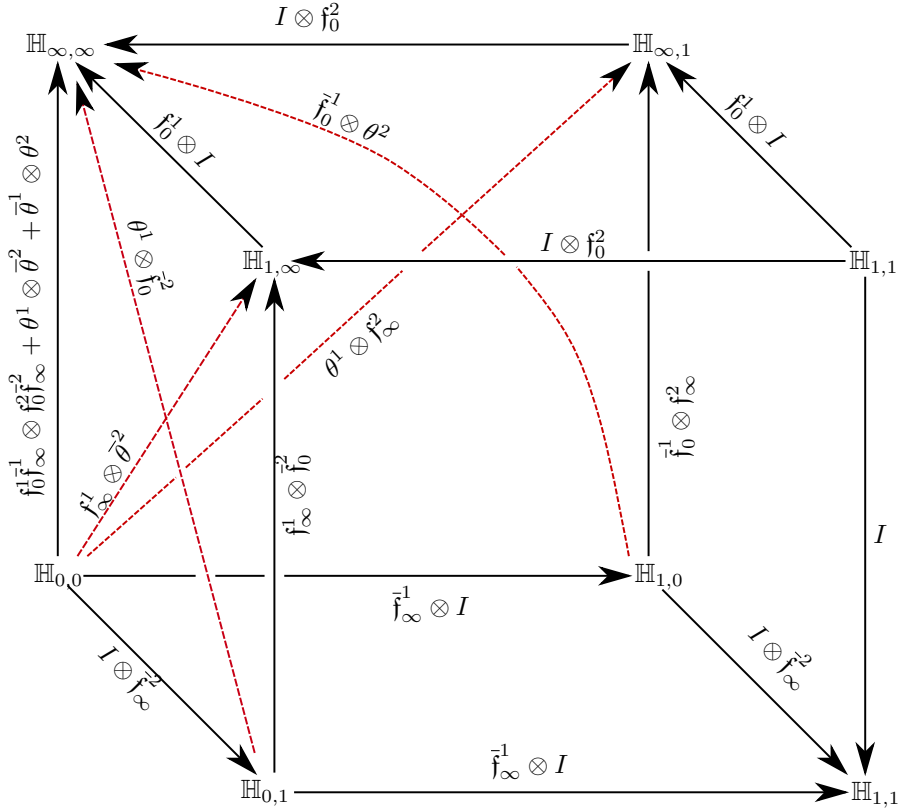
 \caption{\label{fig:cube}
{The above cube determines the chain complex $\Hbbc(K_1,K_2)$ and its differential 
$d_\Hbbc$.}}
\end{center}
\end{figure}

Suppose that a pair of knots $K_1$ and $K_2$ are given. For every $\star,\bullet\in\{0,1,\infty\}$
and $i=1,2$ set $\Hbb_\bullet^i=\Hbb_\bullet(K_i)$,
$\Hbb_{\star,\bullet}=\Hbb_\star^1\otimes \Hbb_\bullet^2$,
$\pphi_\bullet^i=\pphi_\bullet(K_i)$, 
 $\pphibar_\bullet^i=\pphibar_\bullet(K_i)$,
$\theta^i=\theta(K_i)$ and 
 $\thetabar^i=\thetabar(K_i)$. 
 Consider the chain complex
 $(\Hbbc(K_1,K_2),d_\Hbbc)$ constructed as follows.
 The $\Fbb$-module 
 $\Hbbc(K_1,K_2)$ is the direct sum of the modules which appear on the vertices
 of the cube illustrated in the diagram of Figure~\ref{fig:cube}.\\
 
Each directed edge (including the dashed edges) 
in the aforementioned diagram determines a homomorphism from $\Hbbc(K_1,K_2)$
to itself, which is trivial on all summands except for the one which corresponds
to its start point. The map takes the summand corresponding to its start point 
to the summand corresponding to its end point by the homomorphism 
which labels the directed edge.
The differential $d_{\Hbbc}$ of the complex $\Hbbc(K_1,K_2)$ is defined to be the sum of the
homomorphisms which correspond to the directed 
edges of the cube in Figure~\ref{fig:cube}. One should of course make sure that 
$d_\Hbbc\circ d_\Hbbc=0$. However, this follows quickly from the exactness 
of the triangle in (\ref{eq:exact-triangles}).

\begin{thm}\label{thm:cube}
With the above notation fixed, the Heegaard Floer homology of the three-manifold $Y(K_1,K_2)$
obtained by splicing the knot complements $Y_1\setminus K_1$ and $Y_2\setminus
K_2$ is given by 
\begin{displaymath}
\begin{split}
\ov\HFT(Y(K_1,K_2);\Fbb)\simeq H_*\left(\Hbbc(K_1,K_2),d_\Hbbc \right).
\end{split}
\end{displaymath}
\end{thm}

We use the combinatorial description of Heegaard Floer homology by 
Sarkar and Wang (\cite{SW}), which is also adapted for knots in $S^3$ in 
\cite{MOS}, and in \cite{MOST}. These combinatorial descriptions  help
us avoid several technical issues that arise when one  glues
holomorphic curves.\\

For the knots $K_1\subset Y_1$ and $K_2\subset Y_2$ as above define
\begin{displaymath}
\begin{split}
&\chi(K_1,K_2):=\big(h_\infty^1-h_1^1\big)\big(h_\infty^2-h_1^2\big)-
\big(h_0^1-h_1^1\big)\big(h_0^2-h_1^2\big).\\
\end{split}
\end{displaymath}
As a corollary of Theorem~\ref{thm:cube} we prove the following.
 \begin{cor}\label{cor:1}
For $Y=Y(K_1,K_2)$ as above we have 
 \begin{displaymath}
 \begin{split}
&\rank\left(\ov{\HFT}(Y;\Fbb)\right)\geq \max\left\{|\chi(K_1,K_2)|,
|\chi(\ovl{K_1},\ovl{K_2})|\right\},
 \end{split}
 \end{displaymath}
 where $\ovl{K_i}\subset \ovl{Y_i}=-Y_i$ denotes the mirror of $K_i$ in the 
 three-manifold $\ovl{Y_i}=-Y_i$ for $i=1,2$.
 \end{cor}

When one of the two knots is the trefoil, the formula is simplified significantly.
In particular, we prove the following corollary in Section~\ref{sec:example}.\\

\begin{cor}\label{cor:2}
Let $R$ denote the right handed trefoil.
With the above notation fixed 
\begin{displaymath}
\begin{split}
&\rank(\ov{\HFT}(Y(R,K)))\leq h_0(K)+h_1(K)\ \ \text{and}\\
&\rank(\ov{\HFT}(Y(R,K)))\geq 4 \max\{h_0(K),h_1(K),h_\infty(K)\}-
\left(h_0(K)+h_1(K)+2h_\infty(K)\right).
\end{split}
\end{displaymath}
Moreover, if $K$ is non-trivial $Y(R,K)$ is not an $L$-space.
\end{cor}

It is shown by Hedden and Levine that splicing non-trivial knots inside homology 
sphere $L$-spaces never produces an $L$-space \cite{Matt-Adam}. Meanwhile, 
the knot $K$ in Corollary~\ref{cor:2} lives in an arbitrary homology sphere. In this regard,
Corollary~\ref{cor:2} goes beyond the result of Hedden and Levine.\\   

\begin{remark}
The splicing formula of Theorem~\ref{thm:cube} is  different from the splicing formula 
from the original 
$\mathrm{arXiv}$ version of the paper. The results of a few other papers of the author are 
based on the splicing formula of this paper. The results of \cite{Ef-Csurgery} 
 remain unchanged, since the  formula presented in Subsection 5.1 
 (more precisely, Equation~(\ref{eq:complex-cube})) 
 which is used in \cite{Ef-Csurgery} remains unchanged. 
 The proof of the main theorem of \cite{Ef-essential} no longer goes through. 
 Fixing the argument requires developing some technology, including a description 
 of the bordered Floer homology for a knot complement only in terms of the knot 
 chain complex associated with the knot. The modifications 
 will appear in an upcoming revision of \cite{Ef-essential}. 
\end{remark}


The author would like to thank Zolt\'an Szab\'o for his suggestions. In particular, it
was his idea that the difficulties with gluing could disappear if the combinatorial
approach of \cite{SW} is used. He would also like to thank  Matt Hedden and Adam Levine,
as well as the anonymous referee 
for pointing out  (several) mistakes in the earlier versions of this paper.\\

\newpage
\section{Graphs of chain complexes}\label{sec:graphs}
\subsection{Oriented graphs and chain complexes}\label{subsec:graphs}
Let $G$ denote an oriented graph without oriented loops, 
which consists of a set $V(G)$ of vertices and 
a set $$E(G)\subset V(G)\times V(G)$$ of directed edges. For every $e=(v_1,v_2)\in E(G)$
we let $v_s(e)=v_2$ and $v_t(e)=v_1$. The edge $e$ is thus oriented from its 
starting vertex $v_s(e)$ towards its terminal vertex $v_t(e)$. The condition that $G$ does not 
contain any oriented loops implies that there is no sequence $e_1,\dots ,e_k\in E(G)$ with 
the property 
$$v_t(e_i)=v_s(e_{i+1}),\ i=1,\dots,k-1\ \ \ \text{and}\ \ \ v_t(e_k)=v_s(e_1).$$
\begin{defn}
Let $G$ denote an oriented graph without any oriented loops as above. 
A collection $\{(C_v,d_v)\}_{v\in V(G)}$ 
of chain complexes, together with the chain maps 
$$\big\{f_e\co C_{v_s(e)}\lra C_{v_t(e)}\ \big|\ e\in E(G)\big\}$$
is called a {\emph{graph of complexes}} if for every $v_1,v_2\in V(G)$ 
\begin{equation}\label{eq:G-complex}
\sum_{\substack{e_1,e_2\in E(G)\\
v_s(e_1)=v_1, \ v_t(e_2)=v_2\\
v_t(e_1)=v_s(e_2)}} f_{e_2}\circ f_{e_1}=0.
\end{equation}
\end{defn}
Associated with a graph of complexes as above let $C_G=\bigoplus_{v\in V(G)} C_v$ and define 
the differential $d_G\co C_G\lra C_G$ as follows. For $c\in C_v\subset C_G$ let 
$$d_G(c)=\sum_{w\in V(G)}d_{G,w}(c),$$ where $d_{G,w}(c)\in C_w$ is defined by 
\begin{displaymath}
d_{G,w}(c)=\begin{cases}
d_v(c)\ \ \ &\text{if }w=v\\
f_e(c)\ \ \ &\text{if there exists }e\in E(G)\ \text{with }v_s(e)=v\ \text{and }v_t(e)=w\\
0\ \ \ &\text{otherwise.}
\end{cases}
\end{displaymath} 
\begin{defn}
The chain complex $(C_G,d_G)$ is called the chain complex {\emph{associated with the graph
$G$ of chain complexes}}.
\end{defn}

The condition (~\ref{eq:G-complex}) implies that $d_G\circ d_G=0$, i.e. that 
$(C_G,d_G)$ is a chain complex, since each $f_e$ is a chain map. The chain complex 
$(C_G,d_G)$ is usually represented by drawing the oriented graph $G$, labelling 
each vertex $v\in V(G)$ by the chain complex $(C_v,d_v)$ (or simply by $C_v$ if there is 
no confusion) and labelling each oriented edge $e$ by the chain map $f_e$. \\

Let $G$ denote an oriented graph without any loops. It is then possible to label the vertices 
of $G$ by $1,2,\dots, n$ so that for each $e\in E(G)$ we have $v_s(e)<v_t(e)$ (as numbers
in $\{1,...,n\}$). Correspondingly, the chain complexes associated with the vertices of 
$G$ may be labelled $(C_1,d_1),\dots, (C_n,d_n)$.
Let $H$ denote the graph with vertices $1,\dots, n$ and edges 
$$E(H)=\left\{(i,j)\ |\ i,j\in\{1,\dots,n\}\ \text{and}\ i>j\right\}.$$
For $e\in E(H)$ let $g_e=f_e$ if $e\in E(G)$ and $g_e=0$ otherwise. Associated with 
$\{C_i\}_{i\in V(H)}$ and $\{g_e\}_{e\in E(H)}$ we thus find the complex $(C_H,d_H)$,
which is identified with $(C_G,d_G)$. In other words, we may always assume that the 
underlying graph in a graph of complexes is the complete oriented graph $H$. The condition 
(\ref{eq:G-complex}) in this case is equivalent to 
\begin{displaymath}
\sum_{i>k>j} g_{(i,k)}\circ g_{(k,j)}=0\ \ \ \ \ \ \ \ \ \ \ \ \forall\ i,j\in\{1,\dots,n\}.
\end{displaymath}

\subsection{Replacing chain complexes with their homology}
When the ring of coefficients is $\Fbb=\Z/2\Z$ we would like to replace each complex 
$(C_i,d_i)$ in $(C_H,d_H)$ with $(H_*(C_i,d_i),0)$
at the expense of modifying the chain maps $\{g_e\}_{e\in E(H)}$ so that the homology 
of the chain complex associated with the graph of chain complexes remains intact.
Let us begin with a lemma.

\begin{lem}\label{lem:simplification}
Suppose that a chain complex $(C,d_C)$ is decomposed, as a vector space over
$\Fbb$, as $C\simeq A\oplus A\oplus B$ for some vector spaces $A$ and $B$. Suppose 
 that the differential $d_C$ of $C$ has the following block form in this decomposition 
$$d_C=\left(\begin{array}{ccc}0&I_A&f_1\\ 0&0&f_2\\ g_1&g_2&h
\end{array}\right).$$
Then $d_B=h+g_2f_1 \co B\ra B$ is a differential and $H_*(C,d_C)=H_*(B,d_B)$.
\end{lem}
\begin{proof}
Since $d_C$ is a differential, $f_1g_2=0$ and the matrix 
$P=\left(\begin{array}{ccc}I&0&0\\ 0&I&f_1\\ g_2&0&I
\end{array}\right)$ is thus its own inverse.
Since $d_C^2=0$ we get
\begin{displaymath}
P\left(\begin{array}{ccc}0&I_A&f_1\\ 0&0&f_2\\ g_1&g_2&h
\end{array}\right)P=
\left(\begin{array}{ccc}0&I_A&0\\ 0&0&0\\ 0&0&h+g_2f_1
\end{array}\right).
\end{displaymath}
This completes the proof of the lemma.
\end{proof}
We  refer to the procedure which changes the chain complex $(C,d_C)$ to 
the chain complex $(B,d_B)$ as the {\emph{cancellation}} of the two subspaces 
$A\oplus 0\oplus 0\simeq A$ and $0\oplus A\oplus 0\simeq A$ of 
$C$ against each-other.\\

The differential $d_i$ of  $C_i$ may be used to decompose 
$C_i$ as $A_i^1\oplus H_i \oplus A_i^2$, where
$A_i^1$ and $A_i^2$ are two copies of the same $\Fbb$-module $A_i$,
so that $d_i$ takes the form
$$d_i=\left(\begin{array}{ccc}
0&0&I_{A_i}\\ 0&0&0\\ 0&0&0
\end{array}\right)\ \ \ \ \ \ \ \ i=1,\dots, n.$$
Note that $H_i=H_*(C_i,d_i)$ is in fact the homology of the complex $C_i$. 
In particular, $H_i\subset \Ker(d_i \co C_i\ra C_i)$.
Since $d_{v_t(e)}\circ g_e=g_e\circ d_{v_s(e)}$, in this basis 
the matrix block presentation of $g_e$ is of the form
$$g_e=\left(\begin{array}{ccc}
M_e&P_e&N_e\\ 0& G_e& Q_e\\ 0&0& M_e
\end{array} \right)\ \ \ \ \ \ \ \ \forall\ e\in E(H).$$

Initially, the block presentation for $d_H$ is of the form
\begin{displaymath}
d_H=\left(\begin{array}{ccccc}
d_1&0&0&\hdots &0\\
g_{(2,1)}&d_2&0&\hdots &0\\
g_{(3,1)}&g_{(3,2)}&d_3& \hdots &0\\
\vdots &\vdots & \vdots & \ddots & \vdots\\
g_{(n,1)}&g_{(n,2)}&g_{(n,3)}&\hdots &d_n
\end{array}\right).
\end{displaymath}
Replacing the above $3\times 3$ block presentations for $g_{(i,j)}$ and $d_i$ 
the homomorphism $d_H$ takes
 a $3n\times 3n$ block presentation, where $n$ of the block entries are the identity matrices 
corresponding to $d_1,\dots, d_n$. Lemma~\ref{lem:simplification} may be used inductively to 
cancel  $A_i^1$ against $A_i^2$ for $i=1,\dots, n$, and modify the remaining 
blocks correspondingly.
Straight-forward linear algebra implies the following lemma.

\begin{lem}\label{lem:G-complex}
Fix the above notation and for $i,j\in\{1,\dots,n\}$ let 
\begin{displaymath}
h_{(i,j)}=G_{(i,j)}+\sum_{\ell\geq 1}\ \sum_{i>k_1>k_2>\dots >k_\ell>j}
Q_{(i,k_1)}N_{(k_1,k_2)}\dots N_{(k_{\ell-1},k_\ell)} P_{(k_\ell,j)}.
\end{displaymath}
Then the homology of the chain complex associated with $H$, the complexes
$\{(C_i,d_i)\}_{i=1}^n$ and the chain maps $\{g_e\}_{e\in E(H)}$ is isomorphic 
 to the homology of the chain complex associated with $H$, the complexes
$\{(H_i,0)\}_{i=1}^n$ and the homomorphisms $\{h_e\}_{e\in E(H)}$.
\end{lem}  

For $\ell=1$ set $k=k_1$. 
For $h_j\in \Ker(G_{(k,j)} \co H_j\ra H_k)$, $P_{(k,j)}(h_j)=d_k(a_k)$ for 
some $a_k\in A_k$. The element $a_k$ may of course be modified by 
adding to $a_k$ an element $h_k\in H_k$. From here,
$Q_{(i,k)}P_{(k,j)}(h_j)$ is equal to 
$g_{(i,k)}(a_k)$ upto the addition of an element in $g_{(i,k)}(H_k)$. In particular, we find 
a natural well-defined map 
$$\theta_{(i>k>j)} \co \Ker(G_{(k,j)})\lra \Coker(G_{(i,k)}),$$
and $Q_{(i,k)}P_{(k,j)}$ is an extension of $\theta_{(i>k>j)}$ to a homomorphism
from $H_j$ to $H_k$. It is however important to note that simultaneous replacement 
of the maps $Q_{(i,k)}P_{(k,j)}$ with arbitrary extensions of $\theta(i>k>j)$ in 
Lemma~\ref{lem:simplification} is not  {\emph{a priori}} possible.  
\\  

In this paper, we will face situations where each complex $C_i$ is of the form 
$C_i^1\otimes C_i^2$ and each chain map $g_{(i,j)} \co C_j\lra C_i$ is of the form 
$g^1_{(i,j)}\otimes g^2_{(i,j)}$, where $g_{(i,j)}^1 \co C_j^1\lra C_i^1$
and $g_{(i,j)}^2 \co C_j^2\lra C_i^2$ are chain maps. 
In this situation, we may choose the decompositions 
$C^r_i=A^r_i\oplus H^r_i \oplus A^r_i$, for $r=1,2$. Subsequently, note that 
$$H_i=H_i^1\otimes H_i^2\ \ \text{and}\ \ 
A_i=\left(A_i^1\otimes A_i^2\right)\oplus \left(A_i^1\otimes H_i^2\right)
\oplus \left(H_i^1\otimes A_i^2\right)\oplus \left(A_i^1\otimes A_i^2\right).$$
Moreover, corresponding to $g_{(i,j)}^r$ we obtain the blocks 
$M_{(i,j)}^r,N_{(i,j)}^r,P_{(i,j)}^r,Q_{(i,j)}^r$ and $G_{(i,j)}^r$, for $r=1,2$.\\

We close this section by a pair of simple lemmas
addressing this situation.\\
\begin{lem}\label{lem:theta-1}
In the situation above ,
\begin{displaymath}
\begin{split}
Q_{(i,k)}P_{(k,j)}=&Q_{(i,k)}^1P_{(k,j)}^1\otimes Q_{(i,k)}^2P_{(k,j)}^2\\
&+Q_{(i,k)}^1P_{(k,j)}^1\otimes G_{(i,k)}^2G_{(k,j)}^2\\
&+G_{(i,k)}^1G_{(k,j)}^1\otimes Q_{(i,k)}^2P_{(k,j)}^2.
\end{split}
\end{displaymath}
\end{lem}
\begin{proof}
Choose $h_j\in H_j=H_j^1\otimes H_j^2$. The image of $h_j$ under $g_{(k,j)}$ is in 
$$(A_k^1\oplus H_k^1\oplus 0)\otimes (A_k^2\oplus H_k^2\oplus 0)\subset C_k.$$
In particular, we find
$$P_{(k,j)}(h_j)\in \left(A_k^1\otimes A_k^2\right)\oplus \left(A_k^1\otimes H_k^2\right)
\oplus \left(H_k^1\otimes A_k^2\right)\oplus 0\subset A_k.$$
For $h_j=h_j^1\otimes h_j^2$ the corresponding decomposition of $P_{(k,j)}(h_j)$
is of the form
\begin{displaymath}
\left(P_{(k,j)}^1(h_j^1)\otimes P_{(k,j)}^2(h_j^2),
P_{(k,j)}^1(h_j^1)\otimes G_{(k,j)}^2(h_j^2),
G_{(k,j)}^1(h_j^1)\otimes P_{(k,j)}^2(h_j^2),0\right).
\end{displaymath}\label{eq:B-component}
Note that $Q_{(i,k)}P_{(k,j)}(h_j)$ is the 
$H_j=H_j^1\otimes H_j^2$-component of $g_{(i,k)}P_{(k,j)}(h_j)$.
The components in the above presentation  are thus mapped 
by $Q_{(i,k)}P_{(k,j)}$ to 
\begin{displaymath}
\begin{split}
&Q_{(i,k)}^1P_{(k,j)}^1(h_j^1)\otimes Q_{(i,k)}^2P_{(k,j)}^2(h_j^2),\\
&Q_{(i,k)}^1P_{(k,j)}^1(h_j^1)\otimes G_{(i,k)}^2G_{(k,j)}^2(h_j^2),\\
&G_{(i,k)}^1G_{(k,j)}^1(h_j^1)\otimes Q_{(i,k)}^2P_{(k,j)}^2(h_j^2)\ \ \text{and}\ \ 0,
\end{split}
\end{displaymath}
respectively. This completes the proof of the lemma.
\end{proof}
An interesting particular case of the above lemma is when one of the chain maps 
$g_{(i,k)}^1$ or $g_{(k,j)}^1$ is the identity, where we find 
\begin{displaymath}
\begin{split}
g_{(i,k)}^1=Id\ \Rightarrow\ &Q_{(i,k)}P_{(k,j)}
=G_{(k,j)}^1\otimes Q_{(i,k)}^2P_{(k,j)}^2 \ \ \text{and}\\
g_{(k,j)}^1=Id\ \Rightarrow\ &Q_{(i,k)}P_{(k,j)}
=G_{(i,k)}^1\otimes Q_{(i,k)}^2P_{(k,j)}^2,
\end{split}
\end{displaymath}
respectively.

\begin{lem}\label{lem:theta-2}
With the above notation fixed, if $i>k>l>j$ and $g_{(i,k)}^1$ and $g_{(l,j)}^2$ are both 
the identity map we find 
\begin{displaymath}
Q_{(i,k)}N_{(k,l)}P_{(l,i)}=\left(Q^1_{(k,l)}P^1_{(l,j)}\right)\otimes 
\left(Q^2_{(i,k)}P^2_{(k,l)}\right).
\end{displaymath}
\end{lem}
\begin{proof}
Following the proof of Lemma~\ref{lem:theta-1}, 
for $h_j=h_j^1\otimes h_j^2\in H_j^1\otimes H_j^2$ one finds 
$$P_{(l,j)}(h_j)=\left(0,P_{(l,j)}(h_j^1)\otimes h_j^2,0\right)\in A_l.$$
The image of this element of $A_l$ under $g_{(k,l)}^1\otimes g_{(k,l)}^2$ is 
precisely 
$$g_{(k,l)}^1P_{(l,j)}^1(h_j^1)\otimes g_{(k,l)}^2(h_j^2)\subset 
C_j^1\otimes (H_j^2\oplus A_j^2).$$

On the other hand, the domain of $C_{(i,k)}$ (where it is non-zero) is the subset 
$$0\oplus 0\oplus (H_k^1\otimes A_k^2)\oplus 0\subset A_k.$$
In other words, only the component of $g_{(k,l)}^1P_{(l,j)}^1(h_j^1)\otimes g_{(k,l)}^2(h_j^2)$
which lands in $H_k^2\otimes A_k^2$ survives under the map $Q_{(i,k)}$. The aforementioned 
component is precisely $Q_{(k,l)}^1P_{(l,j)}^1(h_j^1)\otimes P_{(k,l)}^2(h_j^2)$, and the image 
of this element under $Q_{(i,k)}$ is precisely 
$Q_{(k,l)}^1P_{(l,j)}^1(h_j^1)\otimes Q_{(i,k)}^2P_{(k,l)}^2(h_j^2)$. This completes the proof 
of the lemma.
\end{proof}

\newpage
\section{A pair of exact triangles}
\subsection{The chain maps}
Let $K\subset Y$ denote a null-homologous knot and fix a Heegaard diagram
$$\ov{H}=(\Sig,\alphas=\{\alpha_1,...,\alpha_g\},\ov\betas=\{\beta_1,...,
\beta_{g-1}\})$$
 for the knot complement $Y\setminus K$. Set
$\betas_\bullet=\{\beta_1^\bullet,...,\beta_{g-1}^\bullet,\lambda_\bullet\}$, where 
$\beta_i^\bullet$ is an isotopic copy of the curve $\beta_i$, and $\lambda_\bullet$ is 
chosen so that the Heegaard triple $(\Sig,\alphas,\betas_\bullet)$ corresponds to the 
three-manifold obtained from $Y$ by $\bullet$-surgery on the knot $K$. Choose the 
curves $\lambda_0,\lambda_1$, and $\lambda_\infty$ so that each two of them 
have a unique transverse intersection point. The orientation on $K$ induces an 
orientation on the three curves $\lambda_0,\lambda_1$, and $\lambda_\infty$. \\

We assume that 
the intersection pattern of $\lambda_0,\lambda_1$ and $\lambda_\infty$ is  one of 
the two patterns illustrated in Figure~\ref{fig:punctures}. This gives the Heegaard 
quadruples
\begin{displaymath}
H=(\Sig,\alphas,\betas_0,\betas_1,\betas_\infty; u,v,w)\ \ \ \text{and}\ \ \ 
\ovl{H}=(\Sig,\alphas,\betas_0,\betas_1,\betas_\infty;\ovl{u},\ovl{v},\ovl{w}).
\end{displaymath}
Note that there is an identification 
$\ov\CFT(\Sig,\alphas,\betas_\bullet;u,v,w)=
\ov\CFT(\Sig,\alphas,\betas_\bullet;\ubar,\vbar,\wbar)$ for 
 $\bullet\in\{0,1,\infty\}$. 
Moreover, for $\bullet,\star\in\{0,1,\infty\}$ the complexes
$\ov\CFT(\Sig,\betas_\bullet,\betas_\star;u,v,w)$,
and $\ov\CFT(\Sig,\betas_\bullet,\betas_\star;\ovl{u},\ovl{v},\ovl{w})$ are 
 identical and the corresponding homology group is 
$\widehat{\mathrm{HF}}(\#^{g-1}(S^1\times S^2))$.
The top generator $\Theta=\Theta_{\bullet,\star}$
 in this Heegaard Floer homology group may be
used to define two holomorphic triangle maps  
(see \cite{OS-3mfld} for more details on the definition of holomorphic 
triangle maps).\\
\begin{defn}\label{Def:phi}
Associated with the Heegaard triples
$$H_{\bullet}=H\setminus \betas_\bullet\ \ \text{and}\ \
 \ovl{H}_{\bullet}=\ovl{H}\setminus \betas_\bullet$$
 define the maps
 \begin{displaymath}
 \begin{split}
 &\phi(H_{0}),\phi(\overline{H}_{0}) \co \widehat{\mathrm{CF}}
(\Sig,\alphas,\betas_1;u,v,w)
\lra \widehat{\mathrm{CF}}(\Sig,\alphas,\betas_\infty;u,v,w),\\
 &\phi(H_{1}),\phi(\overline{H}_{1}) \co \widehat{\mathrm{CF}}
(\Sig,\alphas,\betas_\infty;u,v,w)
\lra \widehat{\mathrm{CF}}(\Sig,\alphas,\betas_0;u,v,w)\ \ \text{and}\\
 &\phi(H_{\infty}),\phi(\overline{H}_{\infty}) \co \widehat{\mathrm{CF}}
(\Sig,\alphas,\betas_0;u,v,w)
\lra \widehat{\mathrm{CF}}(\Sig,\alphas,\betas_1;u,v,w)
 \end{split}
 \end{displaymath}
to be the holomorphic triangle maps corresponding to the  triply punctured 
Heegaard triples $H_0,\ovl{H}_0, H_1,\ovl{H}_1,H_{\infty}$ and $\ovl{H}_{\infty}$, 
respectively, defined using the top generators $\Theta_{\bullet,\star}$.
Denote the induced maps in homology by 
 $\phi_*(H_{\bullet})$ and $\phi_*(\overline{H}_{\bullet})$ and set
\begin{displaymath}
\pphi_\bullet(K):=\phi_*(H_{\bullet})\ \ \text{and}\ \  
 \pphibar_\bullet(K):=\phi_*(\ovl{H}_{\bullet})\ \ \ \ \ \text{for}\ 
 \bullet\in\{0,1,\infty\}. 
\end{displaymath} 
 
\end{defn}

\subsection{Behaviour under Heegaard moves}
The group $\widehat{\mathrm{HF}}(\Sig,\alphas,\betas_\bullet;u,v,w)$, denoted by 
$\Hbb_\bullet(K)$, is independent of the particular Heegaard 
diagram used for the  definition. We  have thus defined the maps 
$$\pphi_0(K),\pphibar_0(K) \co \Hbb_1(K)\ra \Hbb_\infty(K)\ \ \text{and}\ \ 
\pphi_\infty(K),\pphibar_\infty(K) \co \Hbb_0(K)\ra \Hbb_1(K).$$
The definition of the map $\pphi_0(K)$ depends on a Heegaard triple 
$(\Sig,\alphas,\betas_1,\betas_\infty; u,v,w)$ associated with the knot $K$.
Changing $H$ to another Heegaard triple changes $\Hbb_1(K)$ and $\Hbb_\infty(K)$ 
by an isomorphism which is determined by the corresponding Heegaard moves which
change one Heegaard diagram to the other. We would now like to show that the 
the corresponding change in the triangle maps $\pphi_0(H)$ and $\pphibar_0(H)$
respects the above isomorphisms. This justifies using the names $\pphi_0(K)$ and 
$\pphibar_0(K)$ for the above two homomorphisms.
The same statement would be true for 
$\pphi_\bullet(K)$ and $\pphibar_\bullet(K)$\\ 

 Let 
$\{*\}=\{0,1,\infty\}\setminus\{\bullet,\star\}$.
Suppose that two  marked Heegaard triples
$$H_*=(\Sig,\alphas,\betas_\bullet,\betas_\star,u,v,w)\ \ \text{and}\ \ 
H'_*=(\Sig',\alphas',\betas_\bullet',\betas_\star',u',v',w')$$ 
correspond to the same knot $K\subset Y$ for a pair
 $(\bullet,\star)\in\{(\infty,1),(1,0)\}$.  
Similarly, one may consider the Heegaard diagrams $\ovl{H}_*$ and $\ovl{H}'_*$.
 Suppose 
  furthermore that the maps
\begin{displaymath}
\begin{split}
&\imath_\bullet \co \widehat{\mathrm{HF}}(\Sig,\alphas,\betas_\bullet;u,v,w)
\ra \widehat{\mathrm{HF}}(\Sig',\alphas',\betas'_\bullet;u',v',w')\ \ \text{and}\\
&\imath_\star \co \widehat{\mathrm{HF}}(\Sig,\alphas,\betas_\star;u,v,w)\ra
\widehat{\mathrm{HF}}(\Sig',\alphas',\betas'_\star;u',v',w')
\end{split}
\end{displaymath}
are the isomorphisms of the corresponding Heegaard Floer homology groups
associated with the Heegaard moves (and the change of almost complex structure)
changing one Heegaard diagram to the other.
\begin{thm}\label{thm:1}
With the above notation fixed, 
$$\pphi_*(H_*)\circ \imath_{\bullet}
=\imath_{\star}\circ \pphi_*(H'_*),\ \ \text{and}\ \
\pphibar_*(\overline{H}_*)\circ \imath_{\bullet}
=\imath_{\star}\circ \pphibar_*(\overline{H}'_*).$$
\end{thm}
\begin{proof}
The proof consists of some standard steps in Heegaard Floer theory, which are 
sketched below for the Heegaard moves. \\

Note that the first Heegaard triple may be changed to the second Heegaard triple
by a sequence of Heegaard moves, supported in the complement of the marked 
points, of the following types.
\begin{itemize}
\item Changing the almost complex structure on the surface $\Sig$.
\item Isotopies of the curves in $\alphas$ which  are supported away from a 
neighbourhood $U$ of $\lambda_\bullet\cap\lambda_\star$ containing the 
marked points $u,v$ and  $w$, so that the curves in each collection remain disjoint.
 \item Handle slides among the curves in $\alphas$ supported away from $U$.
 \item Simultaneous handle slides among $\betas_\bullet\setminus \{\lambda_\bullet\}$ and 
 $\betas_\star\setminus \{\lambda_\star\}$ supported away from $U$.
\item Stabilization and destabilization of the Heegaard triple away from $U$.
\end{itemize}
  
 The independence of the induced map in homology from the choice of the 
 path of almost complex structures follows the corresponding 
 argument of \cite{Ozs-Stip}. Corresponding to each one of the above Heegaard moves, 
 we obtain a holomorphic 
 square map in the level of chain complexes, comprising to a chain homotopy map 
 between the compositions of the chain maps we are interested in.
 More precisely, performing isotopy or handle slide in 
 $\alphas$ would result in a new set of simple closed curves, which may be denoted 
 by $\alphas'$, by slight abuse of notation. The punctured Heegaard $4$-tuple
 $$(\Sig,\alphas,\alphas',\betas_\bullet,\betas_\star;u,v,w)$$
 determines a homomorphism
  \begin{displaymath}
 \begin{split}
 \ov\Phi \co \ov\CFT(\Sig,\alphas,\alphas';u,v,w)\otimes
 \ov\CFT(\Sig,\alphas',\betas_\bullet;u,v,w)\otimes &
 \ov\CFT(\Sig,\betas_\bullet,\betas_\star;u,v,w)\\
 &\lra
 \ov\CFT(\Sig,\alphas,\betas_\star;u,v,w),
 \end{split}
\end{displaymath}
 which is defined by counting holomorphic squares with Maslov index $-1$.
 Using the top closed elements in the complexes 
 $\ov\CFT(\Sig,\alphas,\alphas';u,v,w)$ and 
 $\ov\CFT(\Sig,\betas_\bullet,\betas_\star;u,v,w)$ we obtain a 
 corresponding map 
 $$\Phi \co \ov\CFT(\Sig,\alphas',\betas_\bullet;u,v,w)
 \lra \ov\CFT(\Sig,\alphas,\betas_\star;u,v,w).$$
 
 Let us denote the differentials of the chain complexes 
$$\ov\CFT(\Sig,\alphas',\betas_\bullet;u,v,w) \ \ \text{and}\ \ 
\ov\CFT(\Sig,\alphas,\betas_\star;u,v,w)$$ by 
$d_{\alpha',\beta_\bullet}$ and $d_{\alpha,\beta_\star}$ respectively.
The Heegaard triples $(\Sig,\alphas,\alphas',\betas_\bullet)$ and 
$(\Sig,\alphas,\alphas',\betas_\star)$ determine chain equivalences
\begin{displaymath}
\begin{split}
&\imath(\alphas,\alphas',\betas_\bullet) \co 
\ov\CFT(\Sig,\alphas',\betas_\bullet;u,v,w)
\lra \ov\CFT(\Sig,\alphas,\betas_\bullet;u,v,w)\ \ \text{and}\\
&\imath(\alphas,\alphas',\betas_\star) \co 
\ov\CFT(\Sig,\alphas',\betas_\star;u,v,w)
\lra \ov\CFT(\Sig,\alphas,\betas_\star;u,v,w).
\end{split}
\end{displaymath}
Moreover, we obtain holomorphic triangle maps associated with the 
Heegaard triples $(\Sig,\alphas,\betas_\bullet,\betas_\star)$ and 
$(\Sig,\alphas',\betas_\bullet,\betas_\star)$
which are denoted by 
\begin{displaymath}
\begin{split}
&\phi(\alphas,\betas_\bullet,\betas_\star) \co 
\ov\CFT(\Sig,\alphas,\betas_\bullet;u,v,w)
\lra \ov\CFT(\Sig,\alphas,\betas_\star;u,v,w)\ \ \text{and}\\
&\phi(\alphas',\betas_\bullet,\betas_\star) \co 
\ov\CFT(\Sig,\alphas',\betas_\bullet;u,v,w)
\lra \ov\CFT(\Sig,\alphas',\betas_\star;u,v,w).
\end{split}
\end{displaymath}

 Considering different types of degenerations for a square of Maslov index 
 $0$, we obtain the relation
 \begin{displaymath}
 \begin{split}
d_{\alpha\beta_\star}&\circ \Phi+\Phi\circ d_{\alpha'\beta_\bullet}= 
\imath(\alphas,\alphas',\betas_\star)\circ 
 \phi(\alphas,\betas_\bullet,\betas_\star)+
\phi(\alphas',\betas_\bullet,\betas_\star)
\circ \imath(\alphas,\alphas',\betas_\bullet).
\end{split}
\end{displaymath}
The induced relation in homology gives the claim for the invariance of 
$\phi_*(H_*)$ under handle slides in $\alpha$. The corresponding argument 
for  $\phi_*(\overline{H}_*)$ is done by changing the marked 
points.\\

The invariance under handle slides among the $\beta$ curves is proved similarly,
and we only highlight the important modifications. Let 
$\betas_\bullet'$ and $\betas_\star'$ be obtained from $\betas_\bullet$ and 
$\betas_\star$ by handle-slides which correspond to a handle slide in 
$\widehat{\betas}$. We thus have the following square of chain maps
\begin{displaymath}
\begin{diagram}
\ov{\CFT}(\Sig,\alphas,\betas_\bullet;u,v,w)&\rTo{\phi(\alphas,\betas_\bullet,\betas_\star)}
&\ov{\CFT}(\Sig,\alphas,\betas_\star;u,v,w)\\
\dTo{\imath_\bullet(\alphas,\betas_\bullet,\betas_\bullet')}&
\rdTo{\phi(\alpha,\betas_\bullet,\betas_\star')}&
\dTo_{\imath_\bullet(\alphas,\betas_\star,\betas_\star')}\\
\ov{\CFT}(\Sig,\alphas,\betas_\bullet';u,v,w)&\rTo{\phi(\alphas,\betas_\bullet',\betas_\star')}
&\ov{\CFT}(\Sig,\alphas,\betas_\star';u,v,w),
\end{diagram}
\end{displaymath}
while the quadruples $(\Sig,\alphas,\betas_\bullet,\betas_\bullet',\betas_\star';u,v,w)$ and 
$(\Sig,\alphas,\betas_\bullet,\betas_\star,\betas_\star';u,v,w)$ determine a pair of 
holomorphic square maps 
$$\Phi_1,\Phi_2 \co \ov{\CFT}(\alpha,\betas_\bullet;u,v,w)\lra 
\ov\CFT(\alphas,\betas_\star';u,v,w).$$
Considering different possible degenerations of holomorphic squares of Maslov index $0$
gives the relations
\begin{displaymath}
\begin{split}
d_{\alpha \beta_\star'}\circ \Phi_1 +\Phi_1\circ d_{\alpha \beta_\bullet}&=
\phi(\alphas,\betas_\bullet',\betas_\star')\circ \imath(\alphas,\betas_\bullet,\betas_\bullet')+
\phi(\alphas,\betas_\bullet,\betas_\star')\\
d_{\alpha \beta_\star'}\circ \Phi_2 +\Phi_2\circ d_{\alpha \beta_\bullet}&=
\imath(\alphas,\betas_\star,\betas_\star')\circ \phi(\alphas,\betas_\bullet,\betas_\star)+
\phi(\alphas,\betas_\bullet,\betas_\star')\\
\end{split}
\end{displaymath}
If we set $\Phi=\Phi_1+\Phi_2$ we thus find
$$
d_{\alpha \beta_\star'}\circ \Phi+\Phi\circ d_{\alpha \beta_\bullet}=
\phi(\alphas,\betas_\bullet',\betas_\star')\circ \imath(\alphas,\betas_\bullet,\betas_\bullet')+
\imath(\alphas,\betas_\star,\betas_\star')\circ \phi(\alphas,\betas_\bullet,\betas_\star)$$
which completes the proof of the invariance under handle slides of the $\beta$ curves
for $\phi_*(H_*)$. The argument for $\phi(\ovl{H}_*)$ is completely similar.\\

The proof of the invariance under stabilization and destabilization follows 
the general argument of \cite{Ozs-Stip} as well. 
\end{proof}
\begin{remark}
 This theorem should be compared with the naturality theorem of Ozsv\'ath and  
 Stipsicz in \cite{Ozs-Stip}.
\end{remark}

\begin{lem}\label{lem:exact-triangle}
With the above notation fixed, the triangles
\begin{equation}\label{eq:lem-exact-triangles}
\begin{diagram}
\Hbb_0(K)&&\lTo{\pphi_1(K)}&&\Hbb_\infty(K)\\
&\rdTo{\pphi_0(K)}&&\ruTo{\pphi_\infty(K)}&\\
&&\Hbb_1(K)&&
\end{diagram}\ \ \text{and}\ \ 
\begin{diagram}
\Hbb_0(K)&&\lTo{\pphibar_1(K)}&&\Hbb_\infty(K)\\
&\rdTo{\pphibar_0(K)}&&\ruTo{\pphibar_\infty(K)}&\\
&&\Hbb_1(K)&&
\end{diagram}
\end{equation}
are both exact.
\end{lem}
\begin{proof}
 The more general forms of exact triangles associated 
with pointed Heegaard diagrams are discussed
in Section 9 of \cite{AE}, using a generalization of Lemma 4.4 from 
\cite{OS-algebra}. The arguments are rather standard and are omitted from the paper. 
The only remark is that if the intersection pattern of $\lambda_0,\lambda_1$ and 
$\lambda_\infty$ follows the left-hand-side of Figure~\ref{fig:punctures}, 
the contributing holomorphic triangles for $(\Sig,\betas_0,\betas_1,\betas_\infty;u,v,w)$
come in cancelling pairs, allowing us to follow the standard arguments. For the Heegaard 
triple $(\Sig,\betas_0,\betas_1,\betas_\infty;\ubar,\vbar,\wbar)$, however, there is a unique 
contributing triangle class, which corresponds to the small triangle bounded between the 
three-curves, which implies that the corresponding triangle map takes 
$\Theta_{0,1}\otimes \Theta_{1,\infty}$ to $\Theta_{0,\infty}$. Nevertheless, the position 
of the punctures in this case implies that the map 
$\phi(\alphas,\betas_0,\betas_\infty;\ubar,\vbar,\wbar)$ which is defined using 
$\Theta_{0,1}$ is trivial (unlike $\phi(\alphas,\betas_\infty,\betas_0;\ubar,\vbar,\wbar)$).
From here, the rest of the argument is standard.
\end{proof}

The exactness of the triangles in (\ref{eq:lem-exact-triangles}) implies that 
$\Ker(\pphi_\infty(K))$ is isomorphic to $\Coker(\pphi_0(K))$, while 
$\Ker(\pphibar_\infty(K))$ 
is isomorphic to $\Coker(\pphibar_0(K))$. Furthermore, the first isomorphism is induced by
the natural chain map $\pphi_1(K)$ while the 
second isomorphism is induced by 
$\pphibar_1(K)$. Let $\theta(K) \co \Hbb_0(K)\lra \Hbb_\infty(K)$ 
denote a map which has the same 
rank as $\pphi_1(K)$ and induces the inverse of the isomorphism
$$\pphi_1(K) \co \Ker(\pphi_\infty(K))\lra \Coker(\pphi_1(K)),$$
while $\thetabar(K) \co \Hbb_0(K)\lra \Hbb_\infty(K)$ denotes a map which has the same 
rank as $\pphibar_1(K)$ and induces the inverse of the isomorphism
$$\pphibar_1(K) \co \Ker(\pphi_\infty(K))\lra \Coker(\pphi_1(K)).$$ 
The choice of the maps $\theta(K)$ and $\thetabar(K)$ are of course not unique.
If 
\[\phi_\infty=\phi(\alphas,\betas_0,\betas_1;u,v,w)
\ \ \text{and}\ \ \phi_0=\phi(\alphas,\betas_1,\betas_\infty;u,v,w)\] 
denote the triangle maps associated 
with the punctured Heegaard triples 
\[(\Sig,\alphas,\betas_0,\betas_1;u,v,w)\ \ \text{and}\ \ 
(\Sig,\alphas,\betas_1,\betas_\infty;u,v,w)\] as above,
the map $\theta(K)$  is in fact the correction term, 
in the sense of Lemma~\ref{lem:G-complex},  associated with the 
sequence (or in fact, graph of complexes)
\begin{displaymath}
\begin{diagram}
\ov\CFT(\Sig,\alphas,\betas_0;u,v,w)&\rTo{\ \ \phi_\infty\ \ } 
&\ov\CFT(\Sig,\alphas,\betas_1;u,v,w)&\rTo{\ \ \phi_0\ \ } 
&\ov\CFT(\Sig,\alphas,\betas_\infty;u,v,w).
\end{diagram}
\end{displaymath}
Similarly, $\thetabar(K)$ corresponds to the sequence 
\begin{displaymath}
\begin{diagram}
\ov\CFT(\Sig,\alphas,\betas_0;\ubar,\vbar,\wbar)&\rTo{\ \ \phibar_\infty\ \ } 
&\ov\CFT(\Sig,\alphas,\betas_1;\ubar,\vbar,\wbar)&\rTo{\ \ \phibar_0\ \ } 
&\ov\CFT(\Sig,\alphas,\betas_\infty;\ubar,\vbar,\wbar),
\end{diagram}
\end{displaymath}
where 
 \[\phibar_\infty=\phi(\alphas,\betas_0,\betas_1;\ubar,\vbar,\wbar)\ \ \text{and}
 \ \ \phibar_0=\phi(\alphas,\betas_1,\betas_\infty;\ubar,\vbar,\wbar).\] 

\subsection{Some properties of the maps $\pphi_\bullet(K)$ 
and $\pphibar_\bullet(K)$}\label{subsec:properties}
Our first observation is that changing the orientation of the knot $K$, and
 correspondingly that of $K_{1}$ and $K_{0}$, corresponds to changing  
 the markings $u,v,w$ with $\ubar,\vbar,\wbar$ in
 Figure~\ref{fig:punctures}.  Suppose that $(\Sig,\alphas,\betas;z_1,z_2)$ represents
$K_\bullet$, meaning that an oriented longitude for $K_\bullet$ is constructed from gluing
an oriented arc on $\Sig$ from $z_1$ to $z_2$ in the complement of $\alphas$ and an 
oriented arc on $\Sig$ from $z_2$ to $z_1$ in the complement of $\betas$. Then 
 $(\Sig,\alphas,\betas;z_2,z_1)$ is a Heegaard diagram for $-K_\bullet$
(the knot $K_\bullet$ with the reverse orientation) while  
 $(-\Sig,\betas,\alphas;z_2,z_1)$ is a Heegaard diagram for  $K_\bullet$.
  The chain complexes associated with the above three 
Heegaard diagrams are identical. Heegaard moves give chain homotopy equivalences
\begin{displaymath}
\tau_\bullet(K) \co \ov{\CFT}(\Sig,\alphas,\betas;z_1,z_2)\lra 
\ov{\CFT}(-\Sig,\betas,\alphas;z_2,z_1)=\ov{\CFT}(\Sig,\alphas,\betas;z_1,z_2).
\end{displaymath}
These chain homotopy equivalences induce the involutions
\begin{displaymath}
\tau_\bullet(K) \co \Hbb_\bullet(K)\lra \Hbb_\bullet(K),\ \ \ \bullet\in\{0,1,\infty\}.
\end{displaymath}
In terms of these isomorphisms  
\begin{equation}\label{eq:duality}
\begin{split}
&\pphibar_0(K)=\tau_{\infty}(K)\circ \pphi_0(K)\circ \tau_{1}(K),\\ 
&\pphibar_1(K)=\tau_{0}(K)\circ \pphi_1(K)\circ \tau_{\infty}(K)\ \  
 \text{and}\\
&\pphibar_\infty(K)=\tau_{1}(K)\circ \pphi_\infty(K)\circ \tau_{0}(K).
\end{split}
\end{equation}
Note however, that the equality
$\thetabar(K)=\tau_\infty(K)\theta(K)\tau_0(K)$ is only satisfied for the induced maps 
from $\Ker(\pphibar_\infty(K))$ to $\Coker(\pphibar_0(K))$.\\

The exactness of the sequences in (\ref{eq:lem-exact-triangles}) implies that  in appropriate 
decompositions 
\begin{equation}\label{eq:decompositions}
\begin{split}
&\Hbb_0(K)=\frac{\Hbb_0(K)}{\Ker(\pphi_\infty(K))}\oplus \Ker(\pphi_\infty(K))
=:\A_\infty(K)\oplus \A_1(K),\\ 
&\Hbb_1(K)=\frac{\Hbb_1(K)}{\Ker(\pphi_0(K))}\oplus \Ker(\pphi_0(K))
=:\A_0(K)\oplus \A_\infty(K)\ \ \ \text{and}\\
 &\Hbb_\infty(K)=\frac{\Hbb_\infty(K)}{\Ker(\pphi_1(K))}\oplus \Ker(\pphi_1(K))
=:\A_1(K)\oplus \A_0(K)
\end{split}
\end{equation}
we have $\pphi_\bullet(K)=\colvec{0&0\\ I_{a_\bullet(K)}&0}$,
where $a_\bullet(K)$ denotes the rank of $\A_\bullet(K)$ for every $\bullet\in \{0,1,\infty\}$.
In this basis we may present the matrices $\tau_\bullet(K)$  as
\begin{displaymath}
\tau_\bullet(K)=\left(\begin{array}{cc}
A_\bullet(K)& B_\bullet(K)\\ C_\bullet(K)& D_\bullet(K)
\end{array}\right),\ \ \ \
\bullet\in\{0,1,\infty\}.
\end{displaymath}

The map $B_0(K)$ corresponds to the induced map
$$\tau_0(K) \co \Ker(\pphi_\infty(K))\ra \frac{\Hbb_0(K)}{\Ker(\pphi_\infty(K))}.$$
The decomposition
$\Hbb_0(K)=\A_\infty(K)\oplus \A_1(K)$
may be modified using a change of basis  of the form
$P_X=\colvec{
I&0\\ -X&I}$, which does not change the block presentations of the 
maps $\pphi_\infty(K)$ and $\pphi_1(K)$. In the new 
basis $\tau_0(K)$ has the following presentation:
\begin{displaymath}
\begin{split}
\tau_0(K)&=\left(\begin{array}{cc}
I&0\\ -X&I
\end{array}\right)\left(\begin{array}{cc}
A_0(K)&B_0(K)\\
C_0(K)&D_0(K)
\end{array}\right)\left(\begin{array}{cc}
I&0\\ -X&I
\end{array}\right)
\\&
=\left(\begin{array}{cc}
A_0(K)-B_0(K)X&B_0(K)\\
\star & -XB_0(K)+D_0(K)
\end{array}\right)
\end{split}
\end{displaymath}
If $B_0(K)$ is injective, we may thus assume that $D_0(K)=0$, while if 
$B_0(K)$ is surjective, we may assume that $A_0(K)=0$.
With a similar reasoning, if $B_\bullet(K)$ is injective we may assume that $D_\bullet(K)=0$,
while if $B_\bullet(K)$ is surjective we may assume that $A_\bullet(K)=0$. \\

In the above decompositions for $\Hbb_\bullet(K)$, the map 
$\theta(K):\Hbb_0(K)\ra \Hbb_\infty(K)$ takes the form
\begin{displaymath}
\theta(K)=\left(\begin{array}{cc}
X&I\\Z&Y
\end{array}\right),
\end{displaymath}
since the induced map from $\A_1(K)\subset \Hbb_0(K)$ to 
$\A_1(K)\subset \Hbb_\infty(K)$ is the inverse of the map induced by 
$\pphi_1(K)$, i.e. the identity. Moreover, since the rank of $\theta(K)$ is the 
same as the rank of $\pphi_1(K)$, we conclude that $Z=YX$. Applying the 
change of basis $P_Y$ on $\Hbb_0(K)$ and the corresponding change of basis 
$P_X$ on $\Hbb_\infty(K)$, $\theta(K)$ takes the form
\begin{displaymath}
\left(\begin{array}{cc}I&0\\ -Y&I\end{array}\right)
\left(\begin{array}{cc}X&I\\ YX&Y\end{array}\right)
\left(\begin{array}{cc}I&0\\ -X&I\end{array}\right)=
\left(\begin{array}{cc}0&I\\ 0&0\end{array}\right).
\end{displaymath}
It is thus possible to choose the above decompositions so that 
$\theta(K)=\left(\begin{array}{cc}0&I\\ 0&0\end{array}\right)$.
If this is the case, the $2\times 2$ presentation of 
$\tau_\infty(K)\thetabar(K)\tau_0(K)$ would be of the form
\begin{displaymath}
\tau_\infty(K)\thetabar(K)\tau_0(K)=\left(\begin{array}{cc}
M&I\\Q&P
\end{array}\right),
\end{displaymath}
and since the ranks of $\theta(K)$ and $\thetabar(K)$ are the same, we find
$Q=PM$.
\subsection{Relative $\SpinC$ structures}\label{subsec:SpinC}
 The vector spaces $\Hbb_\infty(K)$ and 
$\Hbb_{1}(K)$ are naturally decomposed by relative $\SpinC$ classes 
in $$\RelSpinC(Y,K)= \RelSpinC(Y_{1}(K),K_1)=\Z,$$
where the identification with $\Z$ is made using the first Chern class (divided by $2$).
Similarly, the relative $\SpinC$ classes corresponding to $K_{0}$ are 
identified with $\half+\Z$. Thus
\begin{displaymath}
\Hbb_\bullet(K)=\bigoplus_{i\in\Z}\Hbb_\bullet(K,i),\ \ \bullet\in\{1,\infty\}\ \ \text{and}
\ \ \Hbb_0(K)=\bigoplus_{j\in\half+\Z}\Hbb_0(K,j).
\end{displaymath}
Note that $\tau_\bullet(K)$ takes $\Hbb_\bullet(K,i)$ isomorphically to 
$\Hbb_\bullet(K,-i)$ for $\bullet=0,1,\infty$.\\

Let
 $H_{0}=(\Sig,\alphas,\betas_1,\betas_\infty;u,v,w)$ be a Heegaard 
triple used for defining $\pphi_0(K)$. If
$\x\in\Ta\cap\mathbb{T}_{\beta_1}$ and  $\y\in\Ta\cap\mathbb{T}_{\beta_\infty}$ 
are two generators connected by a triangle class $\Delta\in\pi_2(\x,\Theta_{1,\infty},\y)$
with $n_u(\Delta)=n_w(\Delta)=0$ (as observed in  the
surgery exact sequences of \cite{OS-surgery}) $c_1(\relspincs_{u,w}(\x))=
c_1(\relspincs_{u,v}(\y))$. This observation, together with 
(\ref{eq:duality}) imply that the maps $\pphi_0(K)$ and $\pphibar_0(K)$ are 
decomposed as
\begin{displaymath}
\begin{split}
&\pphi_0(K)=\bigoplus_{i\in\Z} \pphi_0(K,i),\ \ \ \pphi_0(K,i) \co \Hbb_1(K,i)\lra \Hbb_\infty(K,i)
\ \ \text{and}\\
&\pphibar_0(K)=\bigoplus_{i\in\Z} \pphibar_)(K,i),\ \ \ \pphibar_0(K,i) \co \Hbb_1(K,i)
\lra \Hbb_\infty(K,i).
\end{split}
\end{displaymath} 

The map $\pphibar_\infty(K) \co \Hbb_0(K)\ra \Hbb_1(K)$ drops the  $\SpinC$ 
grading by $\frac{1}{2}$, while the map $\pphi_\infty(K) \co \Hbb_1(K)\ra \Hbb_0(K)$ increases
the  $\SpinC$ grading by $\frac{1}{2}$. The corresponding decompositions are
thus the following
\begin{displaymath}
\begin{split}
&\pphi_\infty(K)=\bigoplus_{i\in\Z} \pphi_\infty(K,i),\ \ \ 
\pphi_\infty(K,i) \co \Hbb_0\left(K,i-\half\right)\lra \Hbb_1(K,i)
\ \ \text{and}\\
&\pphibar_\infty(K)=\bigoplus_{i\in\Z} \pphibar_\infty(K,i),\ \ \ 
\pphibar_\infty(K,i) \co  \Hbb_0\left(K,i+\half\right) \lra\Hbb_1(K,i).
\end{split}
\end{displaymath} 

In particular,  for a knot $K$ of genus $g$ the maps
$\pphibar_\infty(K,g)$ and $\pphi_\infty(K,-g)$ are trivial, since 
$\Hbb_0(K,g+\half)=\Hbb_0(K,-g-\half)=0$ by the Theorem~3.2 of
\cite{Ef-Whitehead}. Moreover, 
\begin{displaymath}
\begin{split}
&\pphi_1(K)=\bigoplus_{i\in\Z} \pphi_1(K,i),\ \ \ \pphi_1(K,i) \co  \Hbb_\infty(K,i)
\lra\Hbb_0\left(K,i-\half\right)\ \ \text{and}\\
&\pphibar_1(K)=\bigoplus_{i\in\Z} \pphibar_1(K,i),\ \ \ 
\pphibar_1(K,i) \co \Hbb_\infty(K,i)\lra\Hbb_0\left(K,i+\half\right) .
\end{split}
\end{displaymath}

 Let us now assume that $(\Sig,\alphas,\betas_0,\betas_1,\betas_\infty;u,v,w)$
 is one of the Heegaard quadruples 
 illustrated in Figure~\ref{fig:general-triple}. If we drop the marked point $u$ 
 (respectively the marked point $w$) from the Heegaard diagram, 
associated with either of the two resulting punctured Heegaard quadruples we obtain a triangle 
of chain maps, 
 \begin{displaymath}
 \begin{split}
& \begin{diagram}
\ov\CFT(\Sig,\alphas,\betas_\infty;v,w)&&&&
\lTo{\phi(\alphas,\betas_1,\betas_\infty;v,w)}&&&&\ov\CFT(\Sig,\alphas,\betas_1;v,w)\\
&&\rdTo(3,2)_{\phi(\alphas,\betas_\infty,\betas_0;v,w)\ \ } &&&& 
\ruTo(3,2)_{\ \ \phi(\alphas,\betas_0,\betas_1;v,w)} && \\ 
&&&& \ov\CFT(\Sig,\alphas,\betas_0;v,w) &&&& 
\end{diagram} \\
&\begin{diagram}
\ov\CFT(\Sig,\alphas ,\betas_1;u,v)&&&&
\lTo_{\phi(\alphas,\betas_0,\betas_1;u,v)}&&&&\ov\CFT(\Sig,\alphas, \betas_0;u,v)\\
&&\rdTo(3,2)_{\phi(\alphas,\betas_1,\betas_\infty;u,v)\ \ } &&&& 
\ruTo(3,2)_{\ \ \phi(\alphas,\betas_\infty,\betas_0;u,v)} && \\ 
 &&&& \ov\CFT(\Sig,\alphas,\betas_\infty;u,v) &&&& 
\end{diagram}
 \end{split}
 \end{displaymath}
 The domain of any holomorphic triangle which contributes to 
 $\phi(\alphas,\betas_1,\betas_\infty ;v,w)$ has coefficient $1$ precisely 
 at one of the base points $u$ and $\ubar$, and coefficient $0$ at the other 
 one. In other words,
\begin{displaymath}
\begin{split}
\phi(\alphas,\betas_1,\betas_\infty ;v,w)&=
 \phi(\alphas,\betas_1,\betas_\infty ;u,v,w)+\phi(\alphas,\betas_1,\betas_\infty ;\ubar,v,w)\\
 &=\phi(\alphas,\betas_1,\betas_\infty ;u,v,w)+\phi(\alphas,\betas_1,\betas_\infty ;\ubar,
 \vbar,\wbar).
\end{split}
\end{displaymath}
A similar argument implies that 
\begin{displaymath}
\begin{split}
\phi(\alphas,\betas_0,\betas_1 ;u,v)
 &=\phi(\alphas,\betas_0,\betas_1 ;u,v,w)+\phi(\alphas,\betas_0,\betas_1 ;\ubar,
 \vbar,\wbar).
\end{split}
\end{displaymath} 
 We thus obtain the 
 following two exact triangles, respectively:
 \begin{equation}\label{eq:exactness-2}
\begin{diagram}
\Hbb_\infty(K)&&\lTo{\pphi_0+\pphibar_0}&&\Hbb_1(K)\ \ &
\ \ \Hbb_1(K)&&\lTo{\pphi_\infty+\pphibar_\infty}&&\Hbb_0(K)\\
&\rdTo && \ruTo & &   &\rdTo && \ruTo &\\
&& \ov\HFT(Y_0(K)) && & && \ov\HFT(Y) &&
\end{diagram}, 
 \end{equation}
 where $\pphi_\bullet=\pphi_\bullet(K)$ and $\pphibar_\bullet=\pphibar_\bullet(K)$.  
 The exact triangles in (\ref{eq:lem-exact-triangles}) and (\ref{eq:exactness-2})
 may be used to deduce the following conclusions regarding the ranks of the 
 chain maps
 \begin{equation}\label{eq:ranks}
 \begin{split}
 &\rank(\pphi_\bullet(K))=\rank(\pphibar_\bullet(K))=
 \frac{h_\infty(K)+h_1(K)+h_0(K)-2h_\bullet(K)}{2} \\
 &\rank(\pphi_\bullet(K)+\pphibar_\bullet(K))=
 \frac{h_\infty(K)+h_1(K)+h_0(K)-y_\bullet(K)-h_\bullet(K)}{2} 
 \end{split}
 \end{equation}
 where $h_\bullet(K)$ denotes the rank of $\Hbb_\bullet(K)$ and 
 $y_\bullet(K)$ denotes the rank of $\ov\HFT(Y_\bullet(K))$. \\

\newpage
\section{Combinatorial presentation of the exact triangles}\label{sec:complex}
\subsection{Heegaard diagrams for knot complements}\label{subsec:KinY}\label{sec:KinY}
The aim of this subsection is to construct Heegaard diagrams of particular
type associated with a knot $K$ inside a $3$-manifold $Y$, so that the 
chain complexes $C_\bullet(K)$ and the chain maps $\pphi_\bullet(K)$ and 
$\pphibar_\bullet(K)$ may all be described combinatorially.\\

Let us assume that a framed longitude $\ov\lambda$  for $K$ is given as a simple 
closed curve on the torus boundary of $Y\setminus\nd(K)$. Together with 
the meridian $\ov\mu$ of the knot $K$, $\ov\lambda$ gives a parametrization of the 
boundary of $Y\setminus \nd(K)$. It also determines the three-manifold 
$Y_{\ov\lambda}(K)$ obtained by surgery on $K$. 
The curves $\ov\mu$ and $\ov\lambda$ thus give $Y\setminus \nd(K)$ the structure of a bordered 
three-manifold. As such, we remind the reader that a {\emph{nice Heegaard diagram}}  
$$(\Sig,\alphas=\{\alpha_1,...,\alpha_g\},
\ov\betas=\{\beta_1,...,\beta_{g-1}\},\mu,\lambda;z)$$
for  the bordered three-manifold determined by $(Y,K)$ and $\ov\lambda$ 
 consists of a surface $\Sig$ of genus $g$, a $g$-tuple of disjoint simple
closed curves $\alphas$, a $(g-1)$-tuple of disjoint simple closed curves 
$\ov\betas$, a pair of simple closed curves $\mu$ and $\lambda$ disjoint from 
$\ov\betas$ which intersect in a single transverse point, and a marked point $z$ 
in the complement of all curves in $\Sig$. The data satisfies the following conditions.
\begin{itemize}
\item The diagram $(\Sig,\alphas,\ov\betas)$ corresponds to 
$Y\setminus \nd(K)$, while $(\Sig,\alphas,\ov\betas\cup\{\mu\})$ and 
$(\Sig,\alphas,\ov\betas\cup\{\lambda\})$ correspond to the three-manifolds $Y$ and 
$Y_{\ov\lambda}(K)$ respectively. 
\item All domains in 
$\Sig\setminus(\alphas\cup\ov\betas\cup\{\mu,\lambda\})$
are either bigons, triangles or rectangles, except for the domain $D_z$ containing the 
marked point $z$ which is a $(2N+1)$-gon for some integer $N$.
In particular, $D_z$ contains the single intersection point of 
$\mu$ and $\lambda$ as a corner. 
\item Every curve $\beta_i\in\ov\betas$ contains at least one of the 
$2N+1$ edges of $D_z$.
\end{itemize}
Nice Heegaard diagrams exist by Proposition 8.2 from \cite{LOT}. However, two 
remarks are necessary here. First, note that 
in the aforementioned proposition, the roles of 
$\alpha$ and $\beta$ curves  is the opposite of our convention. In particular, the 
curves $\mu$ and $\lambda$ are $\alpha$-curves in \cite{LOT}. 
The second point is that the third condition above is {\emph{a priori}} not 
guaranteed by Proposition 8.2 from \cite{LOT}. 
However, if $\beta_i$ does not contain any of 
the edges of $D_z$,  all neighboring regions of $\beta_i$ would be bigons or 
rectangles. Since $\beta_i$ is homotopically non-trivial, a computation of the euler 
characteristic for the neighborhood of $\beta_i$ 
(the union of all regions which are neighbors of $\beta_i)$ implies that
all neighboring regions of $\beta_i$ are rectangles. However, this in turn implies 
that for some $j\neq i$, $\beta_j$ is parallel (and thus homologous) to $\beta_i$,
a contradiction. Thus, the third condition is also  guaranteed by 
Proposition 8.2 from \cite{LOT}. 
\\

\begin{figure}\def\svgwidth{13cm}
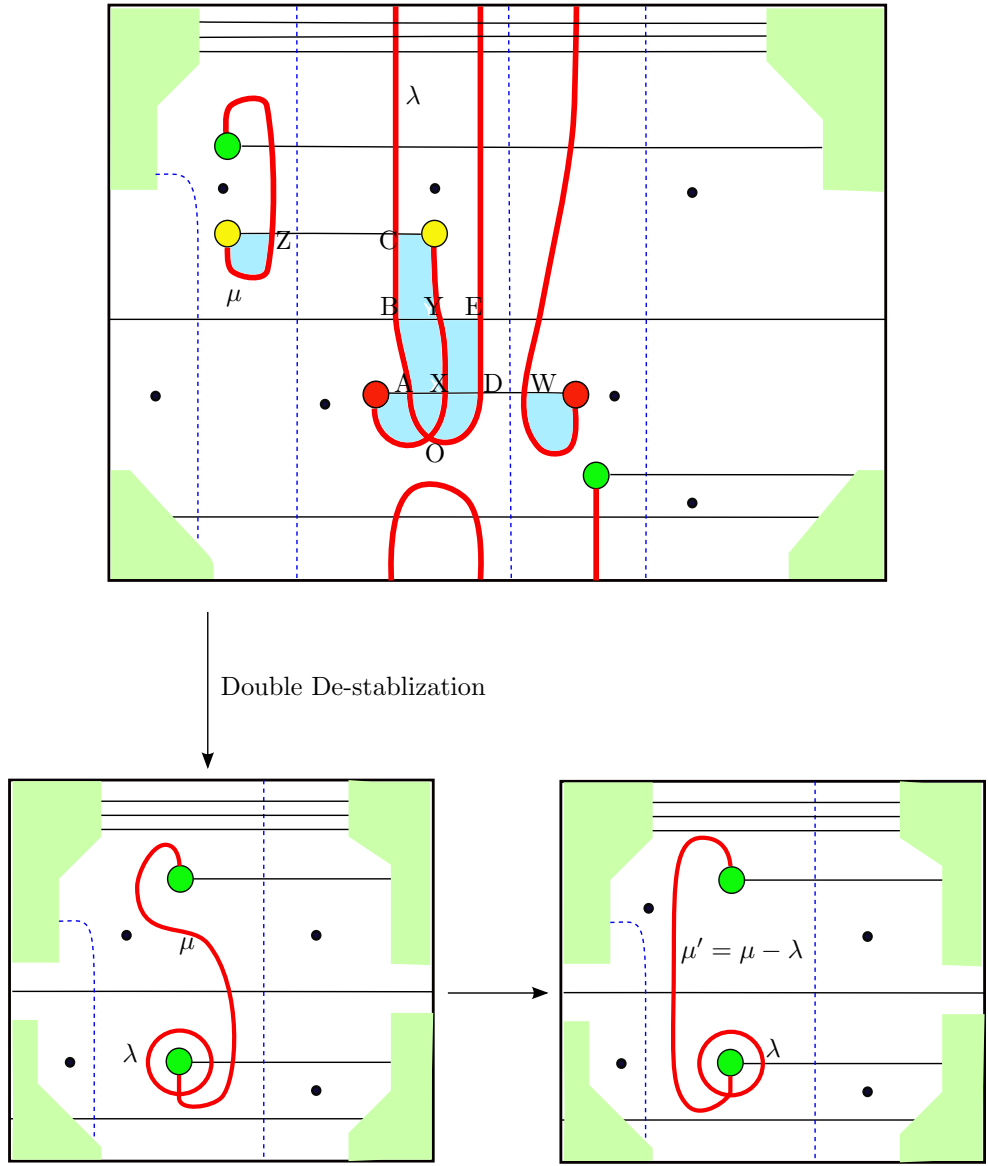
 \caption{\label{fig:Nice}
{Special Heegaard diagrams for knot complements  are 
the union of the genus $4$ surface $\Sig_1$ with boundary illustrated as the 
white part of the figure on top with another surface with boundary. 
The curves $\lambda$ and 
$\mu$ are illustrated as bold red curves, while $\alphas\cap \Sig_1$ and 
$\ov\betas\cap \Sig_1$ are denoted 
by black curves and dashed blue curves, respectively. The intersection of $\mu$ and 
$\lambda$ is denoted by $O$ and some of the intersection points in 
$\alphas\cap(\lambda\cup\mu)$ are labeled (by $A,B,C,D,E,X,Y,Z$ and $W$).
Double Destablization and a change in 
the framing (equivalently, in the parametrization of the boundary torus) gives the two 
Heegaard diagrams on the bottom of the figure. }}
\end{figure}

The picture on the top of Figure~\ref{fig:Nice} describes a surface $\ov\Sig_1$ of genus $4$.
The opposite edges of the rectangle are identified and the pairs of yellow and 
red circles are also glued together (using a horizontal reflection). The pair of green circles are 
identified using a vertical reflection. The solid red  curves are labelled 
$\mu$ and $\lambda$, which meet in a single transverse point $O$. 
The green domains glue together and form a disk $D$ on $\ov\Sig_1$. We set 
$\Sig_1=\ov\Sig_1\setminus \mathrm{Int}(D)$.
The dashed  
blue curves  in $\Sig_1$ correspond to the $\beta$ curves, 
while the solid black curves correspond to
the $\alpha$ curves.  The $\alpha$ and $\beta$ curves may 
have boundary in $\partial D$.

\begin{lem}\label{lem:nice-to-nicer}
Let $K$ be a knot inside a three-manifold $Y$ together with an arbitrary framing. 
Then there is a nice Heegaard diagram 
$(\Sig,\alphas,\ov\betas\cup\{\mu,\lambda\},z)$
for  the corresponding bordered three-manifold with the following properties:
\begin{itemize}
\item $\Sig=\Sig_1\coprod_{\partial \Sig_1=\partial \Sig_2}\Sig_2$, 
where $\ov\Sig_1$ is the surface of genus 
$4$ illustrated in Figure~\ref{fig:Nice} and $\Sig_2$ is a surface with one 
boundary component.
\item The arcs in $\alphas\cap\Sig_1$ are identified with the solid black curves 
in Figure~\ref{fig:Nice}, while the arcs in 
$\ov\betas\cap\Sig_1$ are identified with the dashed blue 
curves in Figure~\ref{fig:Nice}.
\item The curves $\mu$ and $\lambda$ correspond to  the 
bold red curves on $\Sig_1$.
\item The domains on $\Sigma_1$ which contain the bold markings belong to the connected 
component $D_z$ in $\Sig\setminus(\alphas\cup\ov\betas\cup\lambda\cup\mu)$ which contains $z$.
\end{itemize}
\end{lem}

\begin{proof}
Destabilization on $\Sig_1$ gives the equivalent Heegaard diagram which locally 
looks like the surface on the lower left part of Figure~\ref{fig:Nice}. 
Changing $\mu$ to 
$\mu'=\mu-\lambda$ in the aforementioned diagram corresponds to changing the 
parametrization of the boundary. It is thus enough to show that every bordered three-manifold 
with torus boundary admits a nice 
Heegaard diagram  which locally looks like the lower right 
side of   Figure~\ref{fig:Nice}, so that  every domain which meets the green region 
is either a bigon, a rectangle or contains the puncture. If this is the case, every 
domain in the Heegaard diagram illustrated on the upper side of Figure~\ref{fig:Nice}
is either a bigon, or a triangle, or a rectangle, or contains the puncture. In other words,
the diagram on the upper side of Figure~\ref{fig:Nice} is nice.\\

Start with a nice bordered Heegaard diagram for $Y\setminus\nd(K)$
with parametrization given by $\mu'$ and $\lambda$, which exists by 
Proposition 
8.2 from \cite{LOT}.
Denote the intersection point of $\mu'$ and $\lambda$ by $O$. 
Three of the four quadrants around $O$ are triangles, while the last 
quadrant contains the marked 
point $z$. There is thus some curve $\alpha_i$ in $\alphas$ which cuts $\mu'$ in 
the points $D$ and $A$ close to $O$, and the curve  $\lambda$ in $X$ and $W$ 
(close to $O$), so that the picture around $O$ on $\Sig$ is the one illustrated 
in part (a) of Figure~\ref{fig:nice-to-nicer}. We may assume for simplicity 
that $i=g$.
The three triangles are thus $[DOX], [XOA]$ and $[AOW]$. There is a path 
$\gamma$ disjoint from $\betas_0\cup\{\mu',\lambda\}$ which starts from 
the interior of the 
triangle $[AOW]$ and ends at the marked point $z$ and passes only through 
the rectangles. One may add a $1$-handle to $\Sig$ with attaching circles 
placed at the end points of $\gamma$. The core of this 
one handle may be added to $\alphas$ as the curve $\alpha_{g+1}$
and the arc $\gamma$ may be  completed to a simple closed curve $\beta_g$ by 
 attaching its endpoints with an arc going over the one-handle.
This gives a stabilization of the previous Heegaard diagram. 
We may then handle slide $\alpha_{g+1}$ over $\alpha_g$ to obtain the Heegaard 
diagram illustrated in part (b) of Figure~\ref{fig:nice-to-nicer}.\\

 Next, we may add  a one-handle to the Heegaard diagram with attaching circles 
 placed in the middle  of the arcs $[OX]$ and $[OW]$. Denote the arc connecting
  the above two midpoints by $\delta$.
The curve $\mu'$ will be re-named $\beta_{g+1}$, 
the core of this handle would be replaced for  $\mu'$, the curve $\lambda$ will be 
modified by deleting the arc $\delta$ from it, and replacing a corresponding 
arc which travels over the one-handle, and finally, the arc $\delta$ is completed 
to a simple closed curve $\alpha_{g+2}$ using the one-handle.
The new Heegaard diagram is illustrated in part (c) of Figure~\ref{fig:nice-to-nicer}.
This new Heegaard diagram corresponds to the same bordered three-manifold.\\

\begin{figure}\def\svgwidth{14cm}
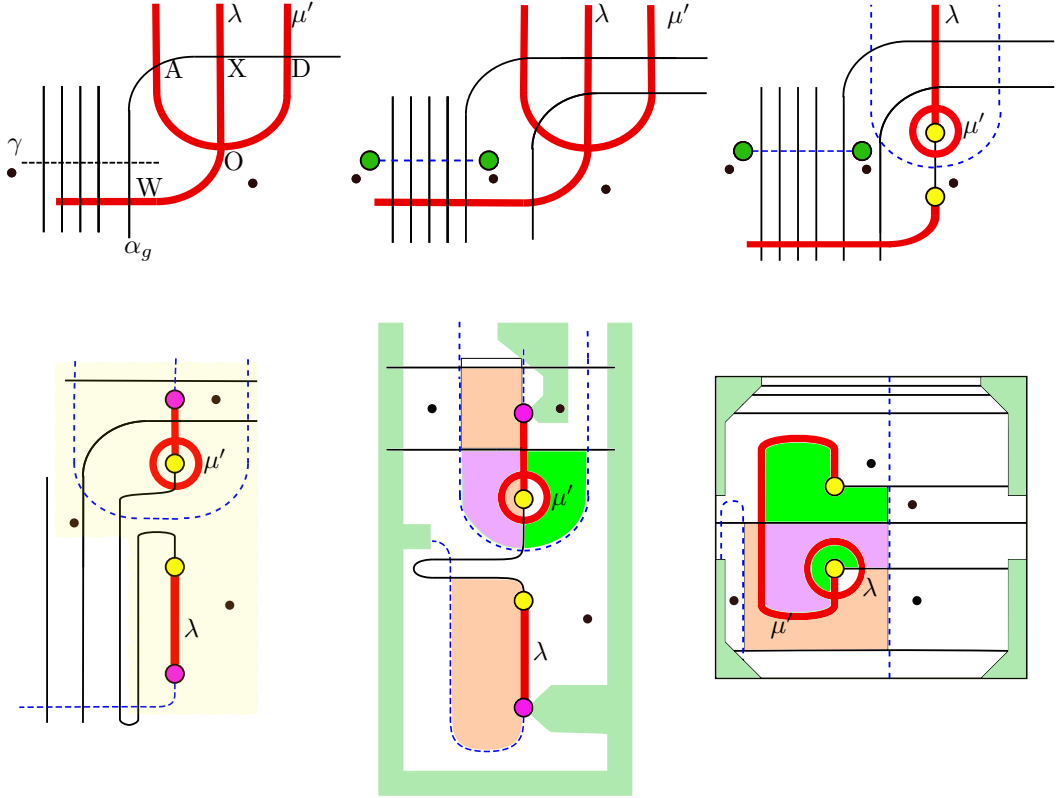
 \caption{\label{fig:nice-to-nicer}
{The $\alpha$ curves are denoted by solid black lines, the
$\beta$ curves are the dashed blue lines, and the curves $\mu'$ and 
$\lambda$ are denoted by bold red lines. 
$D_z$ is the domain containing  bold  circles.
\newline
(a) In a nice Heegaard diagram, three of the quadrants around 
 $O=\mu'\cap\lambda$ are triangles. Use an arc $\gamma$ disjoint 
from $\ov\betas\cup\{\mu',\lambda\}$  to connect the  triangle $[AOW]$ to $z$. 
The closest $\alpha$-curve to  $O$ is  $\alpha_g$.\newline
(b) Attach a handle at the end points of  $\gamma$, 
 complete $\gamma$ to a $\beta$-curve and  slide the core
of the handle  over $\alpha_g$ to produce a new $\alpha$-curve.\newline
(c) Attach a  handle  on $\lambda$ at the 
 two sides of $O$ (the attaching circles are painted yellow). 
 Rename $\mu'$ to  $\beta_{g+1}$ and replace 
the core of the handle for $\mu'$. Push 
 $\lambda$   above the handle and complete the segment  
 on $\mu'$ containing $O$ to  $\alpha_{g+2}$. \newline
(d) Attach a handle on $\lambda$ at the points illustrated by purple
circles. The arcs on $\lambda$ connecting the purple attaching circles to the
yellow attaching circles may be completed to a closed curve, which will 
be replaced for $\lambda$. The complement of these two arcs on initial $\lambda$
may be completed to a $\beta$-curve. The core of the $1$-handle slides over 
$\alpha_g$ to produce the new $\alpha$-curve. Finally, a finger move modifies 
$\alpha_{g+2}$.\newline 
(e),(f) Re-draw the  subsurface of genus $2$ around the intersection of $\mu'$ and 
$\lambda$ which was shaded in part (d).}}
\end{figure}

Next, we attach another one-handle to the Heegaard diagram. 
The attaching circles are placed on $\lambda$ on the two sides of the arc
bounded between the intersection of $\alpha_{g+1}$ and $\lambda$ and the 
intersection of $\mu'$ and $\lambda$. The aforementioned arc may be completed (by adding 
to it a segment which travels over the one-handle) to a simple closed curve, which
will be replaced for $\lambda$. The remainder of (the old) $\lambda$ may also be completed 
(again by adding to it a segment which travels over the one-handle) to a simple 
closed curve which will be denoted by $\beta_{g+2}$. One may slide the core of the 
new one-handle over $\alpha_{g+1}$ to obtain $\alpha_{g+3}$. Finally, we apply a 
finger move isotopy to $\alpha_{g+2}$ to create a pair of intersection points between 
$\alpha_{g+2}$ and $\beta_{g+2}$. The new Heegaard diagram (which still corresponds to 
the same bordered three-manifold) is illustrated in part (d) of   Figure~\ref{fig:nice-to-nicer}
and and a subset of diagram which lives on a sub-surface of genus $2$ 
is re-drawn in part (e) of the same picture, where a $7$-gon and a pair of pentagons are 
painted orange, green and purple, respectively.
 One may then identify the aforementioned subsurface of genus $2$ with  the 
 surface illustrated in part (f). To illustrate the correspondence, 
 the domains corresponding to the $7$-gon and the  two pentagons are painted in the new 
 picture with the relevant color. This completes the proof of the lemma.
\end{proof}

\begin{defn}\label{defn:Nice}
For every knot $K\subset Y$ and every framing $\lambda$ for $K$ the 
Heegaard diagrams of the 
type constructed in Lemma~\ref{lem:nice-to-nicer} are called {\emph{special}} Heegaard 
diagrams.
\end{defn}

\subsection{A combinatorial description of $\pphi_\bullet(K)$ and $\pphibar_\bullet(K)$}

Suppose that $(Y,K)$ denotes a knot $K$ inside a homology sphere $Y$.
Let us assume that 
$$(\Sig,\alphas=\{\alpha_1,...,\alpha_g\},\ov\betas=\{\beta_1,...,\beta_{g-1}\},
\mu,\lambda,z)$$
is  a special Heegaard diagram for the bordered three-manifold 
determined by a zero framed longitude for  $K$ inside 
$Y$. 
The picture around the intersection point $O$ of the simple closed curves 
$\mu$ and $\lambda$ is illustrated on the top of 
Figure~\ref{fig:Nice}. \\

\begin{figure}\def\svgwidth{13.5cm}
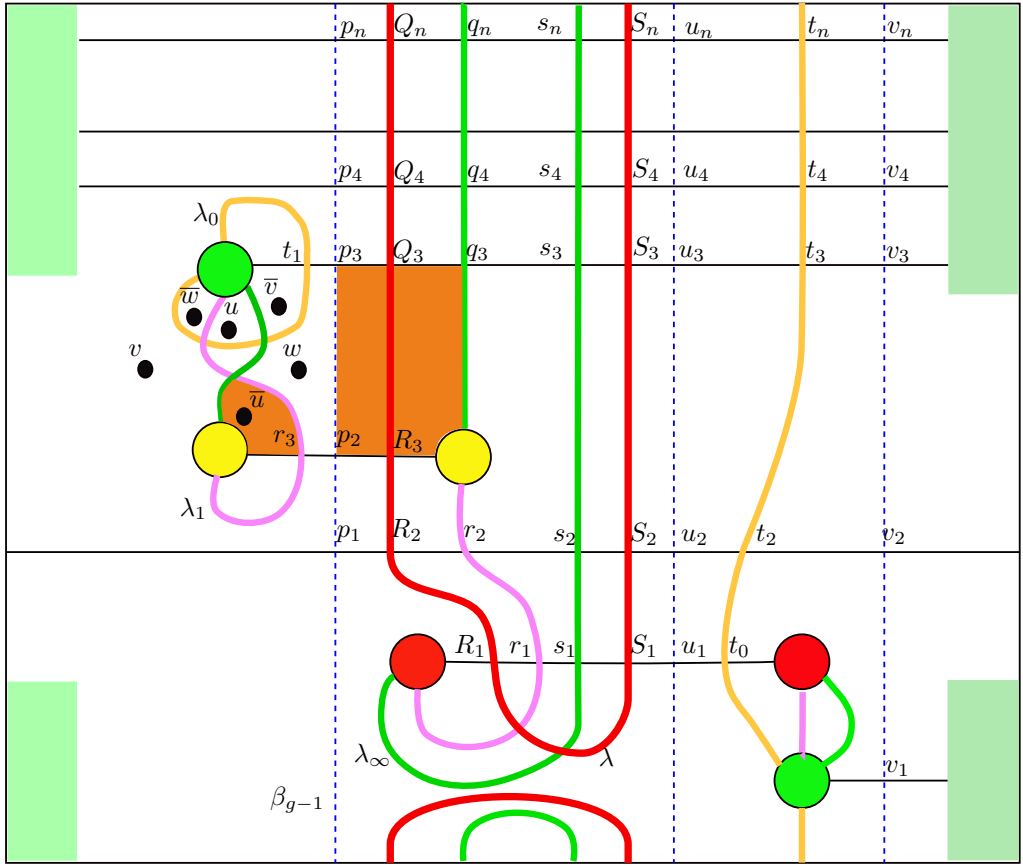
 \caption{\label{fig:H-0-1-inf}\label{fig:H1-and-H0}
{The curves in $\alphas$ are denoted by solid black lines while the curves 
in $\ov\betas$ are denoted by dashed blue lines. Three simple closed 
curves $\lambda_0$, $\lambda_1$ and $\lambda_\infty$ are 
denoted by bold red, purple and green lines respectively.
Six marked points $u,v,w,\ubar,\vbar$ and $\wbar$ are introduced close to the intersection 
points of these three curves. The intersection points on $\beta_{g-1}$, $\lambda_\infty$ 
and $\lambda_1$ are labelled. Associated with  $i\geq 3$ there is a pentagon which
 vertices at $P_0, r_3, p_2, p_i$ and $q_i$. For $i=3$ 
the pentagon is shaded orange in the picture. }}
\end{figure}

We introduce three auxiliary curves 
denoted by $\lambda_\infty,\lambda_0$ and $\lambda_1$ respectively, as in 
the Heegaard diagram illustrated in Figure~\ref{fig:H-0-1-inf}. 
The Heegaard diagram 
$$H_\bullet=(\Sig,\alphas,\ov\betas\cup\{\lambda_\bullet\};u,v,w)
\ \ \text{and}\ \ 
\ovl{H}_\bullet=(\Sig,\alphas,\ov\betas\cup\{\lambda_\bullet\};\ubar,\vbar,\wbar)
$$ 
are (triply punctured) diagrams which corresponds to the knot 
$K_\bullet\subset Y_\bullet(K)$ for  $\bullet\in\{0,1,\infty\}$ (note that 
two of the three punctures are placed in the same connected component of 
$\Sig\setminus(\alphas\cup\ov\betas\cup\lambda_\bullet)$ for 
$\bullet\in\{0,1,\infty\}$). The above claim is checked by computing the 
intersection numbers  of each $\lambda_\bullet$ with the simple closed curves 
$\mu$ and $\lambda$, since the curves are disjoint from  $\ov\betas$.
Each pair of these three curves intersect each-other exactly once.
Each one of the three diagrams $H_\bullet,\ovl{H}_\bullet\ \bullet\in\{0,1,\infty\}$ 
is a nice  Heegaard diagram and  they determine the chain complexes 
$C_\bullet=\ov\CFT(H_\bullet)=\ov\CFT(\ovl{H}_\bullet)$. 
Denote the differential of the
 complex $C_\bullet$  by  $d_\bullet$ for $\bullet\in\{0,1,\infty\}$.
The chain maps $\pphi_\bullet(K)$ and $\pphibar_\bullet(K)$ have a 
simple combinatorial  description, which is discussed in the remainder of 
this section. \\

Fix the labelling of the intersection points of $\lambda_0,\lambda_1,\lambda_\infty,
\beta_{g-1}$ and $\beta_{g-2}$ with the curves in $\alphas$ as in Figure~\ref{fig:H-0-1-inf}.
Let 
$$\{\Thet_0\}=\lambda_1\cap\lambda_\infty,\ \ \ 
\{\Thet_1\}=\lambda_0\cap\lambda_\infty,\ \ \text{and}\ \ 
\{\Thet_\infty\}=\lambda_0\cap\lambda_1. $$
The Heegaard triple 
$$(\Sig,\alphas,\ov\betas\cup\{\lambda_1\},\ov\betas\cup\{\lambda_\infty\};
\ubar,\vbar,\wbar)$$
 determines a {\emph{combinatorial triangle map}} 
 $\fbar_0 \co C_1\ra C_\infty$ as follows. 
Let $\beta_{g-1}$ be the $\beta$-curve which contains the intersection points 
$p_1,p_2,...,p_n$ in Figure~\ref{fig:H-0-1-inf}. 
Let $\x=(x_1,...,x_g)$ be a generator of  $C_1$ with 
$x_i\in \alpha_{\sig(i)}\cap \beta_i$ for some $\sig\in S_g$ and $i=1,...,g-1$, and 
$x_g\in \lambda_1$. Define
\begin{displaymath}
\fbar_0\left(\x\right):=\begin{cases}
(x_1,...,x_{g-1},s_i)\ \ \ &\text{if }x_g=r_i,\ \ i=1,2\\
0\ &\text{otherwise}
\end{cases}.
\end{displaymath}

Similarly, the Heegaard triple 
$$(\Sig,\alphas,\ov\betas\cup\{\lambda_1\},\ov\betas\cup\{\lambda_\infty\};u,v,w)$$ 
determines a combinatorial triangle map $f_0\co C_1\ra C_\infty$  defined by 
\begin{displaymath}
f_0\left(\x\right):=\begin{cases}
(x_1,...,x_{g-2},p_2,q_i)\ \ \ &\text{if } (x_{g-1},x_g)=(p_i,r_3)\ \ i\geq 3,\\
(x_1,...,x_{g-2},p_1,q_i)\ \ \ &\text{if } (x_{g-1},x_g)=(p_i,r_2)\ \ i\geq 3,\\
0&\text{otherwise.}
\end{cases}
\end{displaymath}

The Heegaard triples 
\begin{displaymath}
(\Sig,\alphas,\ov\betas\cup\{\lambda_0\},\ov\betas\cup\{\lambda_1\};u,v,w)
\ \ \text{and}\ \ (\Sig,\alphas,\ov\betas\cup\{\lambda_0\},\ov\betas\cup\{\lambda_1\};
\ubar,\vbar,\wbar),
\end{displaymath}
correspond to  the combinatorial triangle maps $f_\infty,\fbar_{\infty} \co C_0\ra C_1$.
For a generator $\x=(x_1,...,x_g)$ these two maps are defined by
which are defined by setting 
\begin{displaymath}
\begin{split}
&\fbar_\infty\left(\x\right)=\begin{cases}
(x_1,...,x_{g-1},r_1)\ \ \ &\text{if }\ x_g=t_0\\
0\ \ \ &\text{otherwise}
\end{cases}\ \ \ \ \ \text{and}\\
&f_\infty\left(\x\right)=\begin{cases}
(x_1,...,x_{g-2},p_3,r_3)\ \ \ &\text{if }\ (x_{g-1},x_g)=(p_2,t_1)\\
0\ \ \ &\text{otherwise}
\end{cases}.
\end{split}
\end{displaymath}


\begin{lem}{\label{lem:1}}
With the above notation fixed
$f_0\circ f_\infty=\fbar_0\circ \fbar_\infty=0$.
\end{lem}
\begin{proof}
This is trivial from the combinatorial definitions of 
$f_0,\fbar_0, f_\infty$ and $\fbar_\infty$.
\end{proof}

 Let 
$$\Sig\setminus(\alphas\cup\ov\betas\cup \lambda_0\cup \lambda_1)=\big(
\coprod_{i=1}^N D_i)\cup D_\ubar\cup D_\vbar \cup D_\wbar,$$
where $D_\bullet$ are the regions in the complement of these curves, with $D_\ubar,D_\vbar$ 
and $D_\wbar$ the regions containing the
marked points $\ubar,\vbar$ and $\wbar$, respectively. We  set 
$$\beta^0_i=\beta_i,\ \ i=1,...,g-1,\ \ \ \text{and}\ \ \ 
\betas=\{\beta_1,...,\beta_g\}=\ov\betas\cup\{\lambda_0\}.$$ 
The construction of  the Heegaard diagram implies  the following properties:
\begin{itemize}
\item
 The regions $D_2,...,D_N$ are  rectangles or bigons, while $D_1$ 
is a pentagon.
\item
 One of the corners of the pentagon $D_1$ is the unique intersection point 
$P=\Thet_\infty=\lambda_0\cap\lambda_1$, and the three punctures $\ubar,\vbar$ 
and $\wbar$  are placed on three of the quadrants around $P$ (other than the quadrant 
corresponding to $D_1$).
\item
 All the neighbours of $D_1$ (the regions having an edge in common with 
$D_1$) are punctured.
\item
 Each $\beta$-curve is adjacent to at least one of the punctured domains.
\end{itemize}

The edges of the pentagon are five arcs: $2$ of them are on $\lambda_0$ and 
$\lambda_1$, two of them are on the $\alpha$-curves, and one of them is on a 
$\beta$-curve which is assumed to be $\beta_1$. The $\alpha$-curve which 
cuts $\lambda_0$ in a corner of the pentagon is assumed to be $\alpha_1$, 
and the other one is assumed to be $\alpha_2$. Denote the vertices of the 
pentagon by $P=Q_1,Q_2,Q_3,Q_8$ and $Q_6$ in counter-clockwise order, so that 
$Q_1$ is the intersection point of $\lambda_0$ and $\lambda_1$,  $Q_2$ is 
on the intersection of  $\alpha_1$ with $\lambda_0$, and so that $Q_6$ is the 
intersection point of $\lambda_1$ with  $\alpha_2$. \\

For $i=2,...,g-1$
let $\beta^1_i=\gamma_i$ be a parallel copy of $\beta_i$ which is drawn very
close to $\beta_i$, and is slightly pushed to one of the punctured domains adjacent 
to $\beta_i$ by a finger move, so that a pair of intersection points 
(denoted by $X_i$ and $Y_i$) is created between these two curves (see the right hand 
side picture in Figure~\ref{fig:D1}).  Let us assume that the 
 small positively oriented disk connecting these two intersection points (with
  $\beta_i$ on the left and $\gamma_i$ on the right) goes from $X_i$ to $Y_i$. 
  In order to define $\gamma_1$, choose a parallel copy of $\beta_1$ and push it 
  slightly over the intersection point of $\beta_1$ with $\alpha_1$ to obtain $\gamma_1$,
so that a pair of cancelling intersection points $X_1$ and $Y_1$ is created between 
$\gamma_1$ and $\beta_1$ on the two sides of the intersection point $Q_3$ of
 $\alpha_1$ and $\beta_1$, and so that $\gamma_1$ slightly enters the punctured 
 domain next to the $\beta$-edge of the pentagon. The local picture around 
 $D_1$  looks like Figure~\ref{fig:D1} where this procedure is pictured.
Let $\gamma_g$ be the curve $\lambda_1$ and set $\betas_1=\gammas
=\{\gamma_1,...,\gamma_g\}$.\\

In order to construct 
$\beta^\infty_i$, for $i=2,...,g-1$ choose a parallel copy of $\gamma^i=\beta^1_i$ and
as this parallel copy enters the bigon $T_i$ push it into the neighbouring punctured
domain by a finger move.  The curve $\beta^\infty_1$ is constructed as illustrated in 
Figure~\ref{fig:pentagon}. We set
$$\betas_\infty=\{\beta^\infty_1,...,\beta^\infty_{g-1},\lambda_\infty\}.$$

\begin{lem}\label{lem:priodc-triangles}
The punctured Heegaard diagrams 
$$(\Sig,\alphas,\betas_\star,\betas_\bullet;u,v,w)\ \ \text{and}\ \ 
(\Sig,\alphas,\betas_\star,\betas_\bullet;\ubar,\vbar,\wbar)$$
for $(\star,\bullet)\in\{(0,1),(1,\infty)\}$
do not contain any non-trivial positive triply periodic domains.
\end{lem}
 \begin{proof}
Let $\Dcal$ denote a positive triply periodic domain in the Heegaard diagram
 $(\Sig,\alphas,\betas_1,\betas_\infty;u,v,w)$. Thus 
 $$\partial \Dcal=\sum_{i=1}^g a_i \alpha_i+\sum_{i=1}^{g-1} b_i \beta^1_i+ 
 \sum_{i=1}^{g-1} c_i \beta^\infty_i+b\lambda_1+c\lambda_\infty.$$
Let $\Dcal_i$ denote the doubly periodic domain with $\partial \Dcal_i=\beta^1_i-\beta^\infty_i$
for $i=1,...,g-1$.  Setting 
 $\Dcal'=\Dcal-\sum_{i=1}^{g-1}b_i\Dcal_i$ we find
 $$\partial \Dcal'=\sum_{i=1}^g a_i \alpha_i+\sum_{i=1}^{g-1} (c_i-b_i) \beta^\infty_i
 +b\lambda_1+c\lambda_\infty.$$
 Since the left-hand-side is trivial in $H_1(Y\setminus \nd(K);\Z)$, so is the right-hand-side.
 This implies that $c=-b$. Let $\Dcal_0$ denote the triply periodic domain in  the punctured 
 Heegaard triple $(\Sig,\betas_0,\betas_1,\betas_\infty;u,v,w)$ with 
 $\partial \Dcal_0=\lambda_1-\lambda_0-\lambda_\infty$. For $\Dcal''=\Dcal'-b\Dcal_0$ 
 we thus obtain
 $$\partial \Dcal''=\sum_{i=1}^g a_i \alpha_i
 +\sum_{i=1}^{g-1} (c_i-b_i) \beta^\infty_i
 +b\lambda_0.$$
 In other words, $\Dcal''$ is a doubly periodic domain for the nice 
 (and hence weakly admissible) 
 Heegaard diagram
 $$(\Sig,\alphas,\betas_\infty\cup\{\lambda_0\}\setminus\{\lambda_\infty\};u,v,w).$$
The coefficients of $\Dcal''$ and all $\Dcal_i, i=1,...,g-1$ 
over the small triangle bounded between $\lambda_0,\lambda_1$ and 
$\lambda_\infty$ is zero. In other words, the coefficient of 
$$\Dcal=\Dcal''+b\Dcal_0+\sum_{i=1}^{g-1} b_i\Dcal_i$$ 
  over this small triangle is $b$, which should thus be non-negative. 
  Choosing this triangle sufficiently small 
  we may thus assume that the total area of $b\Dcal_0$ is negative, unless $b=0$. 
   
 One may choose the area form on the surface $\Sig$ so that all doubly periodic domains 
 for the punctured Heegaard diagram $(\Sig,\alphas,\betas_\infty\cup\{\lambda_0\}
\setminus \{\lambda_\infty\};u,v,w)$ and all $\Dcal_i, i=1,\dots,g-1$ have zero total area. 
 However, this implies that the total area of 
 $\Dcal$ is the same as the total area of $b\Dcal_0$, which is at most zero. 
 Since $\Dcal$ is a positive 
 domain, we conclude $\Dcal=0$. This completes the proof for the triple
 $(\Sig,\alphas,\betas_1,\betas_\infty;u,v,w)$. The proof for the other triples is 
 completely similar.
 \end{proof}
 
The Heegaard diagrams
\begin{displaymath}
(\Sig,\alphas,\betas_\bullet;u,v,w)\ \ \text{and}\ \ 
(\Sig,\alphas,\betas_\bullet; \ubar,\vbar,\wbar)
\end{displaymath}
are nice, so by Sarkar and Wang \cite{SW}, the differentials of the complexes
\begin{displaymath}
\begin{split}
&\ov{\CFT}(\Sig,\alphas,\betas_\bullet;u,v,w)\ \ \text{and}\ \ 
\ov{\CFT}(\Sig,\alphas,\betas_\bullet; \ubar,\vbar,\wbar)
\end{split}
\end{displaymath}
are given by counts of bigons and rectangles. \\
 
\begin{thm}\label{thm:comb-to-hol-all}
Under the above identification of the chain complexes $(C_\bullet,d_\bullet)$ 
\begin{displaymath}
\begin{split}
f_0 &= \phi\left(\Sig,\alphas,\betas_1,\betas_\infty;u,v,w\right),
\ \ \ \ \ \ \ \ 
\fbar_0 =\phi\left(\Sig,\alphas,\betas_1,\betas_\infty;
\ubar,\vbar,\wbar\right),\\
f_\infty &= \phi\left(\Sig,\alphas,\betas_0,\betas_1;u,v,w\right)
\ \ \text{and}\ \ \ 
\fbar_\infty =\phi\left(\Sig,\alphas,\betas_0,\betas_1;
\ubar,\vbar,\wbar\right).
\end{split}
\end{displaymath}
\end{thm}
\subsection{Proof of Theorem~\ref{thm:comb-to-hol-all}}
 A similar discussion is carried over in \cite{Ef-Csurgery} (and in particular 
Theorem~2.3 from that paper). We repeat 
the proof, in most parts with more details, to keep the paper easier to read.\\

\begin{proof}
We start by proving the statement for $\fbar_\infty$.
Note that the
 top generator $\Theta$ of the Heegaard Floer homology group
$\widehat{\mathrm{HF}}(\#^{g-1}S^1\times S^2)$ coming from the Heegaard 
diagram $(\Sig,\betas,\gammas;\ubar,\vbar,\wbar)$ is 
 the generator $\{P,X_1,...,X_{g-1}\}$. \\

Let $\x=(x_1,...,x_g)$ and $\y=(y_1,...,y_g)$ be generators with 
$x_i\in \alpha_{\sigma(i)}\cap \beta_i$ and $y_i\in \alpha_{\tau(i)}\cap \gamma_{i}$,
with $\sigma,\tau\in S_g$. Let 
$\Delta \co \mathbb{D}\ra \mathrm{Sym}^g(\Sig)$ be the homotopy class of a triangle
in $\pi_2(\x,\Theta,\y)$, with Maslov index zero, such that it supports a holomorphic 
representative, and remains disjoint from the punctures.\\

There are two types of domains in the complement 
$\Sig\setminus(\alphas\cup\betas\cup\gammas)$ of 
the curves, the \emph{large} domains and the \emph{small} domain.
The small domains are those created between the parallel pairs of curves $\gamma_i$ and 
$\beta_i$ ($i=1,...,g-1$),
and their area may be chosen arbitrarily small by choosing $\gamma_i$ close enough to 
$\beta_i$. The large domains are the rest of the domains which are in
correspondence with the domains $D_\bullet,\ \bullet\in\{\ubar,\vbar,\wbar,1,...,N\}$
introduced above. We  abuse the notation and still denote these new 
regions by $D_\bullet$. \\

 Let us assume that the small bigon connecting $X_i$ to $Y_i$ is  denoted by
 $T_i$, and the region having the small interval $[X_i,Y_i]$ on $\beta_i$ in common
with $T_i$ is $D_{i},\ i=2,...,g-1$.
Then there are two triangles with corners $X_i$ and $Y_i$, which have an edge in
common with $D_{i}$, which will be denoted by $R_{i}$ and $L_{i}$ respectively.
For $i=1$, instead of these three regions we have $4$ triangles with one
corner being $X_1$ or $Y_1$, which will be denoted by $R_1,T_1,S_1$ and $L_1$ 
respectively (as they appear while we travel on $\beta_1$ from $X_1$ to $Y_1$, see 
Figure~\ref{fig:D1}). We are implicitly assuming that the regions $D_{i}$ for 
$i=1,...,g-1$ (as described above) are different, while it may happen that this is not 
the case. However, the argument we give below remains true in general, and only 
needs notational corrections.\\

\begin{figure}
\def\svgwidth{12cm}\begin{center}
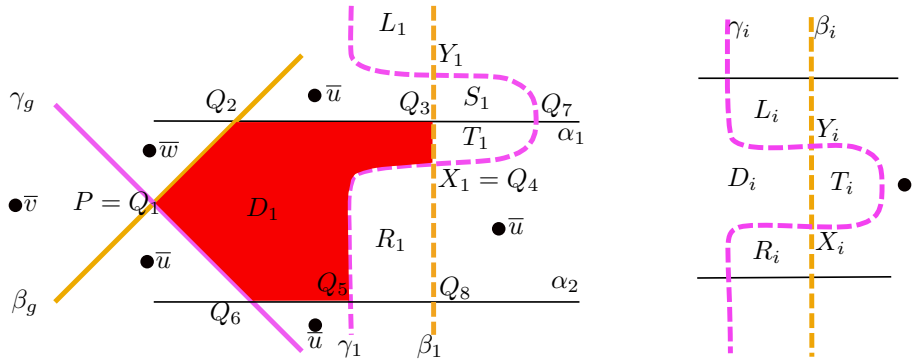
\caption{\label{fig:D1}
{The region around the \emph{pentagon} $D_1$ is illustrated on the left hand side.
 The punctured domains are marked 
by solid circles inside them. The curves in $\betas=\betas_0$, 
$\gammas=\betas_1$ and 
$\alphas$  have colors orange, pink and black respectively. 
The pentagon is changed to a hexagon in the new 
Heegaard diagram which is colored red. The initial pentagon is the union of the 
hexagon $D_1$ with the triangle $R_1$. The right hand side illustrates the labelling 
near the intersection of $\beta_i$ with its  Hamiltonian isotope $\gamma_i$.
}}\end{center}
\end{figure}
Let $\Dcal=\Dcal(\Delta)$
denote the domain (i.e. the $2$-chain on $\Sig$) associated with the triangle class $\Delta$.
Let $d_i\geq 0$ denote the coefficient of $D_i$ in $\Dcal$. Similarly, denote
the coefficients of $T_i, R_i$ and $L_i$ by $t_i,r_i$ and $l_i$ respectively.
The coefficient of $S_1$ will be denoted by $s_1$. Of course, there are other regions
which may appear in $\Dcal$ with positive coefficient, but all such regions are
bigons or rectangles. Since $P$ appears in $\Theta$ and three of the corners around
$P$ are punctured, the coefficient $d_1$ is equal to $1$. \\

Let $P=Q_1,Q_2,...,Q_6$
denote the corners  of $D_1$ (now a hexagon) in counter-clockwise order 
(so $Q_4=X_1$). Since two opposite quadrants around each one of $Q_2$ and 
$Q_6$ are punctured, we have $x_g=Q_2$ and $y_g=Q_6$.  Thus, 
$Q_3$ is not one of $x_1,...,x_g$ and $Q_5$
is not one of $y_1,...,y_g$. Considering the local coefficients around
$Q_3$ we  conclude that $t_1=1+s_1$. If $Q_7$ is the third corner of $T_1$
(other than $Q_3$ and $Q_4$),  in order for $\Dcal$ to be a non-negative
domain, we need $x_1=Q_7$, and the $4$ local coefficients around $Q_7$ are
forced to be $t_1=1+s_1, s_1, 0$ and $0$ in the counter-clockwise order.
Two opposite quadrants around $Y_1$  have zero coefficients in $\Dcal$. Since
 $Y_1$  does not appear in $\Theta$, this implies that $s_1=l_1=0$ (thus $t_1=1$).
Similarly, considering the local coefficients around $p_1$ we conclude $r_1=1$.
Since $Q_5$ is not among $y_1,...,y_g$, the local coefficients around $Q_5$ are
$1,r_1=1,0$ and $0$ in the counter-clockwise order.
Let $Q_8$ be the third corner of $R_1$ other than $Q_4$ and $Q_5$. Since
two opposite corners around $Q_8$ have zero coefficient and $r_1=1$ we  
have $x_1=Q_8$. Thus $\Dcal=\Dcal'+\Dcal_1=\Dcal'+(R_1+D_1+T_1)$ 
where $\Dcal'$ is a non-negative $2$-chain which is disjoint from $\Dcal_1$, and 
$\Dcal_1$ is a hexagon with $5$ acute angles and one obtuse angle, and with vertices 
$\{P,y_g, x_1, X_1,y_1,x_g\}$. The contribution of $\Dcal_1$ to the index of $\Delta$ 
is zero, by Sarkar's formula from \cite{Sar}.\\

By Sarkar's formula for the index of triangles from \cite{Sar}
\begin{equation}\label{eq:Index-Formula}
\mu(\Delta)=e(\Dcal)+\mu_{\x}(\Dcal)+\mu_{\y}(\Dcal)+\bb(\Dcal).\cc(\Dcal)-\frac{g}{2}.
\end{equation}
Here $e(\Dcal)$ is the Euler measure of the domain $\Dcal$, $\bb(\Dcal)$ is the part 
of $\partial \Dcal$ on the $\beta$-curves, and
$\cc(\Dcal)$ is the part of $\partial \Dcal$ on the $\gamma$-curves. 
Furthermore, $\mu_\x(\Dcal)$ and $\mu_\y(\Dcal)$ denote the local contributions
of the intersection points included in $\x$ and $\y$, respectively, to the corners of 
$\Dcal$. We refer to \cite{Sar} for more detailed definitions. Separating 
$\Dcal_1$, which has Maslov index $0$, from $\Dcal$ we obtain the
following equality
$$\mu(\Delta)=e(\Dcal_{s})+e(\Dcal_{l})+\mu_{\x}(\Dcal')+\mu_{\y}
(\Dcal')+\bb(\Dcal').\cc(\Dcal')-\frac{g-2}{2}.$$
Here $\Dcal_{s}$ denotes the part of $\Dcal'$ which uses the regions $D_i,R_i,T_i$ and 
$L_i$ for $i=2,...,g-1$ and $\Dcal_l=\Dcal'-\Dcal_s$. Clearly, $e(\Dcal_l)\geq 0$ and 
$$\Dcal_s=\sum_{i=2}^{g-1}(d_iD_i+t_iT_i+r_iR_i+l_iL_i).$$

Considering the local coefficients around $X_i$ and $Y_i$ we conclude 
$r_i=l_i+1$ and $d_i=t_i+l_i$. Having in mind that $T_i$ are bigons, $R_i$ and $L_i$ 
are triangles and $D_i$ are hexagons, this implies the following computation:
\begin{equation}\label{eq:1}
 \begin{split}
e(\Dcal_s)&=\sum_{i=2}^{g-1}\left((t_i+l_i)e(D_i)+t_ie(T_i)+(l_i+1)e(R_i)+l_ie(L_i)\right)\\
&=\sum_{i=2}^{g-1}\left((t_i+l_i)(-\frac{1}{2})+t_i(\frac{1}{2})+(l_i+1)
(\frac{1}{4})+l_i(\frac{1}{4})\right)=\frac{g-2}{4}.
 \end{split}
\end{equation}

The $1$-chain $\bb(\Dcal')$ is a union of $1$-chains on $\beta_i, i=2,...,g-1$, 
denoted by $\bb_i(\Dcal')$. Similarly we have 
$\cc(\Dcal')=\sum_{i=2}^{g-1}\cc_i(\Dcal')$.
It is clear that $\bb_i(\Dcal')$ and $\cc_j(\Dcal')$ are disjoint unless $i=j$. In this 
latter case, the only possible geometric intersections are at $X_i$ and $Y_i$, where the
intersection numbers are $(l_i+\frac{1}{2})(t_i-\frac{1}{2})$ and $-l_it_i$ respectively. 
Thus
\begin{equation}\label{eq:2}
 \bb(\Dcal').\cc(\Dcal')=\sum_{i=2}^{g-1}\left((l_i+\frac{1}{2})(t_i-\frac{1}{2})
 -l_i.t_i\right)=-\frac{g-2}{4}+\frac{1}{2}\sum_{i=2}^{g-1}(t_i-l_i).
\end{equation}

Let us now consider the coefficients around the intersection points $x_i$ and $y_i$ for 
$i=2,...,g-1$. Since $x_i$ is on $\beta_i$ there are
non-negative integers $a_i,b_i,c_i$ and $e_i$ such that the local coefficients around 
$x_i$ are $a_i,b_i, b_i+l_i+1$ and $a_i+l_i$, and the local coefficients
around $y_i$ are $c_i,e_i,e_i+t_i-1$ and $c_i+t_i$. Thus
\begin{equation}\label{eq:3}
 \mu_\x(\Dcal')+\mu_\y(\Dcal')=
 \frac{1}{2}\sum_{i=2}^{g-1}\big((a_i+b_i+c_i+e_i)+(l_i+t_i)\big).
\end{equation}

Combining (\ref{eq:1}), (\ref{eq:2}) and (\ref{eq:3}) and replacing for the terms 
in the definition of
$\mu(\Delta)$  we obtain
\begin{displaymath}
 \begin{split}
  0=\mu(\Delta)&=e(\Dcal_{s})+e(\Dcal_{l})+\mu_{\x}(\Dcal')
  +\mu_{\y}(\Dcal')+\bb(\Dcal').\cc(\Dcal')-\frac{g-2}{2}\\
&=e(\Dcal_l)-\frac{g-2}{2}+\frac{1}{2}\sum_{i=2}^{g-1}(a_i+b_i+c_i+e_i+2t_i)\\
&\geq \frac{1}{2}\sum_{i=2}^{g-1}(a_i+b_i+c_i+(e_i+t_i-1)+t_i).
 \end{split}
\end{displaymath}

Note that $e_i+t_i-1$ is the coefficient of one of the domains around $y_i$ and is
 thus non-negative. The above inequality thus implies that $a_i=b_i=c_i=t_i=0$ and
$e_i=1$ for $i=2,...,g-1$. Thus, the coefficients on the two sides of $\gamma_i$ 
 either agree, or differ by $1$, and the coefficients on the two sides of $\beta_i$ differ
either by $l_i$ or by $l_i+1$.  If we start from $y_i$, where on the left (or right) side of 
$y_i$ the coefficients on the two sides of $\gamma_i$ are zero, and travel on the $\alpha$ 
curve intersecting $\gamma_i$ (i.e. orthogonal to $\gamma_i$) 
until we get to an intersection point 
with $\beta_i$, as we pass $\beta_i$ the coefficient  
changes either to $-l_i$ or to $-l_i-1$. Since the latter is negative, the former 
happens and $l_i=0$. It is easy to see from here that $x_i$ and 
$y_i$ are the corresponding intersection points of $\beta_i$ and $\gamma_i$ with 
the same $\alpha$ curve, and that the domain $\Dcal'$ is a union of obvious triangles 
which are disjoint from each-other.  \\

We conclude that the domain of $\Delta$ 
is the disjoint union of $g-2$ simple triangles with a hexagon with $5$ acute angles 
and one obtuse angle. It is quite well-known that the moduli space corresponding to 
this homotopy class contributes $1$ to the triangle map, for a generic path of almost 
complex structures. These are thus the only holomorphic triangles which contribute to 
the chain map $\fbar_\infty$ defined using the Heegaard triple 
$(\Sig,\alphas,\betas,\gammas;\ubar,\vbar,\wbar)$. Under the obvious identification of 
$\widehat{\mathrm{CF}}(\Sig,\alphas,\gammas;\ubar,\vbar,\wbar)$ with
$\widehat{\mathrm{CF}}(\Sig,\alphas,\betas_1;\ubar,\vbar,\wbar)$, this is just the map
 which replaces the pair $\{Q_2,Q_8\}$ with $\{Q_6,Q_7\}$.
 This completes the proof of Theorem~\ref{thm:comb-to-hol-all} for $\fbar_\infty$.\\

\begin{figure}
\def\svgwidth{9cm}\begin{center}
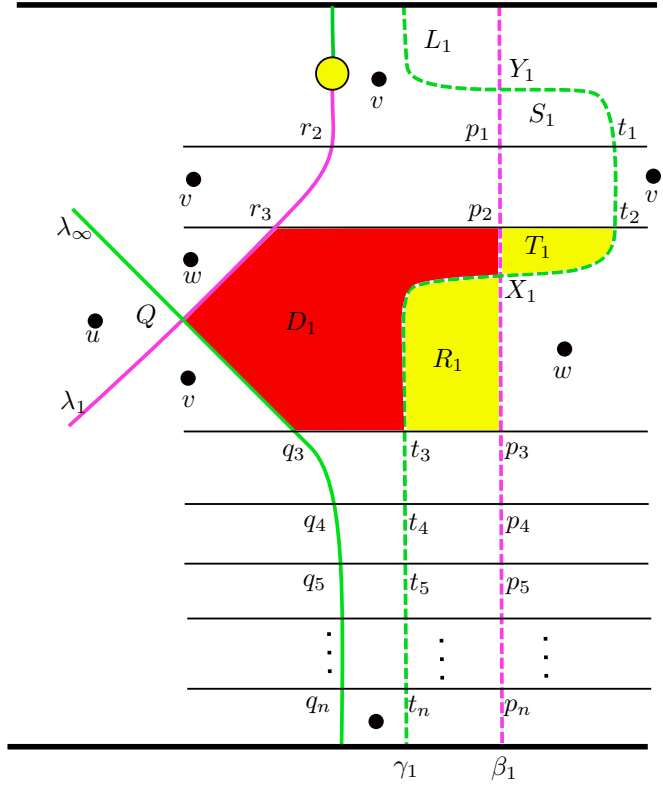
\caption{\label{fig:pentagon}
{
The region around the {hexagon} $D_1$ is illustrated.
The labelling of the intersection points  in the Heegaard diagram,
as well as the labelling of some of the connected components in the complement of the 
curves, is illustrated. The curves in $\betas_1$, $\betas_\infty$ and $\alphas$ have
colors pink,  green and  black, respectively.
}}\end{center}
\end{figure}

The proofs of 
the other three claims are completely similar. In fact, the proofs of the statement of 
the theorem for $\fbar_0$ and $f_\infty$ are even easier, 
since the domains which are not punctured in 
the corresponding Heegaard triple are all bigons, rectangles or triangles. We thus only 
 need  to use the second part of the above argument
in these two cases (and the study of the neighbourhood of the hexagon is not needed).
The proof of the claim for $f_0$ 
requires some more serious modification, which will be outlined below.\\ 

Note that the Heegaard triples $(\Sig,\alphas,\betas_1,\betas_\infty;u,v,w)$ and 
$$(\Sig,\alphas,\betas'=\ov\betas\cup\{\lambda_1\},
\gammas'=\{\beta^\infty_1,\gamma_2,...,\gamma_{g-1},\lambda_\infty\};u,v,w)$$
may be identified using a diffeomorphism of the surface $\Sig$. It is thus enough 
to show that $f_0=\phi(\Sig,\alphas,\betas',\gammas';u,v,w)$. This  allows us  
keep the same labelling for the points $X_i,Y_i,\ i=2,...,g-2$. For the intersection
points on $\gamma'_1=\beta^\infty_1$ and $\beta_1$, as well as some of the 
intersection points on $\lambda_\infty$ and $\lambda_1$ we  use the labelling of 
Figure~\ref{fig:pentagon}.  We abuse the notation and denote the two intersection 
points between $\beta_{1}$ and $\gamma'_{1}$ by $X_1$ and $Y_1$. 
Moreover, some of the regions in the neighbourhood of $X_1$ and $Y_1$ are 
labelled: again by abuse of notation we  denote these regions by $D_1, R_1, L_1,
S_1$ and $T_1$ (see Figure~\ref{fig:pentagon}). Let us use $d_i, r_i,s_i, t_i$ and $l_i$
to denote the coefficients of the domains $D_i,R_i,S_i,T_i$ and $L_i$ in the $2$-chain
$\Dcal$ associated with a holomorphic triangle connecting  $\x=(x_1,...,x_g),
\y=(y_1,...,y_g)$ and $\Theta$ which contributes to 
$\phi(\Sig,\alphas,\betas',\gammas';u,v,w)$. We  assume that for some elements 
$\sig,\tau\in S_g$ 
\begin{displaymath}
x_i\in \begin{cases}\beta_i\cap\alpha_{\sigma(i)} \ \ &\text{if } i=1,...,g-1\\
\lambda_1\cap \alpha_{\sig(g)}\ \ &\text{if } i=g\end{cases}
\ \ \text{and}\ \ 
y_i\in \begin{cases}\gamma_i\cap\alpha_{\tau(i)} \ \ &\text{if } i=2,...,g-1\\
\gamma'_1\cap\alpha_{\tau(1)}\ \ &\text{if }i=1\\
\lambda_0\cap \alpha_{\tau(g)}\ \ &\text{if } i=g\end{cases}
\end{displaymath}

The examination of the coefficients in 
Figure~\ref{fig:pentagon} implies the following:
\begin{itemize}
\item   We have $d_1=r_1=t_1=1$ and $s_1=l_1=0$.
\item   Either $x_g=r_2, y_1=t_1$, or $x_g=r_3, y_1=t_2$.
\item   There are some $j,k\in\{3,4,...,n\}$ such that $y_g=q_k$ and $x_1=p_j$.
\end{itemize}

Let us write $\Dcal=\Dcal_s+\Dcal_l$, where 
$$\Dcal_s:=s_1S_1+\sum_{i=1}^{g-1}\left(d_iD_i+t_iT_i+r_iR_i+l_iL_i\right)\ \ 
\text{and}\ \ \Dcal_l:=\Dcal-\Dcal_s.$$
Considering the local coefficients at $X_i$ and $Y_i$ we find 
$r_i=l_i+1$ and $d_i=t_i+l_i$.
Applying the index formula in (\ref{eq:Index-Formula}) we obtain
\begin{equation}\label{eq:index-step-1}
\begin{split}
0&=e(\Dcal)+\mu_\x(\Dcal)+\mu_\y(\Dcal)+\bfrak(\Dcal).\cfrak(\Dcal)-\frac{g}{2}\\
&=\left(e(\Dcal_l)+\left(-\frac{1}{2}+\frac{1}{4}+\frac{1}{4}\right)+\frac{g-2}{4}\right)
+\mu_\x(\Dcal)+\mu_\y(\Dcal)+\bfrak(\Dcal).\cfrak(\Dcal)-\frac{g}{2}\\
&\geq \mu_\x(\Dcal)+\mu_\y(\Dcal)+\bfrak(\Dcal).\cfrak(\Dcal)-\frac{g+2}{4}.
\end{split}
\end{equation}

The $1$-chains $\bfrak(\Dcal)$ and $\cfrak(\Dcal)$ may be written as 
$$\bfrak(\Dcal)=\sum_{i=1}^g \bfrak_i(\Dcal)\ \ \text{and}\ \ \cfrak(\Dcal)
=\sum_{i=1}^g \cfrak_i(\Dcal)$$
 as before. Note that $\bfrak_1(\Dcal)$ is the arc on $\beta_1$ 
 from $X_1$ to $p_j$, while  $\cfrak_1(\Dcal)$ is the arc from one of $t_1$ and 
 $t_2$ to $X_1$. Moreover, $\bfrak_g(\Dcal)$ is the arc on $\lambda_1$ from 
 $Q$ to one of $r_2$ and $r_3$, while $\cfrak_g(\Dcal)$ is the arc on 
 $\lambda_\infty$ from $q_i$ to $Q$.  Thus
\begin{equation}\label{eq:index-step-2}
\begin{split} 
 \bfrak(\Dcal).\cfrak(\Dcal)&=\left(\frac{-1}{4}+\frac{1}{4}\right)+
 \sum_{i=2}^{g-1}\left((l_i+\frac{1}{2})(t_i-\frac{1}{2})
 -l_i.t_i\right)\\
 &=-\frac{g-2}{4}+\frac{1}{2}\sum_{i=2}^{g-1}(t_i-l_i).
 \end{split}
 \end{equation}
 
 Let us now assume that
 the local coefficients around $x_i$ are $a_i,b_i,b_i+l_i+1$ and $a_i+l_i$,
 while the local coefficients around $y_i$ are $c_i,e_i,e_i+t_i-1$ and $c_i+t_i$
 for $i=2,...,g-1$. The corresponding local coefficients around $x_1,y_1,x_g$ and 
 $y_g$  would be $ (a_1,b_1, b_1+1,a_1)$, $(c_1,e_1,e_1,c_1+1)$, $(0,0,1,0) $ and 
 $(0,0,0,1)$ respectively for some non-negative integers $a_i,b_i,c_i$ and $e_i$ 
 $i=1,...,g-1$. Thus
 \begin{equation}\label{eq:index-step-3}
 \begin{split}
 \mu_\x(\Dcal)+\mu_\y(\Dcal)&=\left(\frac{1}{2}\right)
 +\frac{1}{2}\sum_{i=1}^{g-1}\left((a_i+b_i+c_i+e_i)+(l_i+t_i)\right)\\
 \end{split}
 \end{equation}  
 If we combine (\ref{eq:index-step-1}), (\ref{eq:index-step-2}) and 
 (\ref{eq:index-step-3}) we find
 \begin{displaymath}
 \begin{split}
0&\geq -\frac{g}{2}+\left(\frac{1}{2}\right)+\frac{1}{2}
\sum_{i=1}^{g-1}\left( a_i+b_i+c_i+e_i+2t_i\right)\\
&=\frac{1}{2}\sum_{i=1}^{g-1}\left(a_i+b_i+c_i+(e_i+t_i-1)+t_i\right).
 \end{split}
 \end{displaymath}
 
 As in the proof of the theorem for $\fbar_\infty$ this implies that 
 $a_i=b_i=c_i=t_i=0$, while $e_i=1$ for $i=1,...,g-1$. It is easy to see from 
 here that $j=k$, and complete the proof as before.
\end{proof}

\subsection{The maps $\theta(K)$ and $\thetabar(K)$}\label{subsec:theta}
Let $\Hbb_\bullet$ denote the homology of the chain complex $C_\bullet$ for 
$\bullet\in \{\infty,1,0\}$. 
If we choose a representative $a\in C_0$ of a class 
$$[a]\in\Ker((f_\infty)_*)\subset \Hbb_0$$  
there exists some $b\in C_1$ such that $f_\infty(a)=d_1(b)$. Then
$d_\infty(f_0(b))=f_0(d_1(b))=f_0(f_\infty(a))=0$, so $f_0(b)$ is closed and 
represents a class in $\Hbb_\infty$. If we replace $b$ with another element 
$b'=b+\Delta b$ so that $d_1(b')=f_\infty(a)$,
$\Delta b$ is closed (i.e. it represents an element in $\Hbb_1$). 
The difference $f_0(b')-f_0(b)=f_0(\Delta b)$ is an element in 
$\Image((f_0)_*)$. Thus the class 
$$\theta([a])=[f_0(b)]\in\Coker((f_0)_*)$$
is well-defined. This gives a homomorphism $$\theta=\theta(K) \co 
\Ker((f_\infty)_*)\lra \Coker((f_0)_*).$$
Similarly, we define the map 
$\thetabar=\thetabar(K) \co \Ker((\fbar_\infty)_*)\ra\Coker((\fbar_0)_*)$ from
$$\fbar_\infty \co C_0\ra \C_1\ \ \ \text{and}\ \ \ \fbar_0 \co 
C_1\ra C_\infty.$$

\begin{prop}\label{prop:theta}
The maps $$\theta(K) \co \Ker(\pphi_\infty(K))\ra \Coker(\pphi_0(K))\ \ \ \text{and}\ \ \  
\thetabar(K) \co \Ker(\pphibar_\infty(K))\ra \Coker(\pphibar_0(K))$$ are the inverses
of the maps induced by 
$\pphi_1(K),\pphibar_1(K) \co \Hbb_\infty(K)\ra \Hbb_0(K)$ which sit in the  
exact sequences
\begin{displaymath}
\begin{diagram}
\Hbb_0(K)&&\lTo{\pphi_1(K)} & & \Hbb_\infty(K) &&\Hbb_0(K)&&\lTo{\pphibar_1(K)}
 && \Hbb_\infty(K)\\
&\rdTo{\pphi_\infty(K)} &&\ruTo{\pphi_0(K)}&& \text{and} &&\rdTo{\pphibar_\infty(K)} &&\ruTo{\pphibar_0(K)}&\\
&&\Hbb_1(K)&&&&&&\Hbb_1(K)&&
\end{diagram} 
\end{displaymath}
\end{prop}
\begin{proof}
For this purpose, let us assume that the  Heegaard diagram 
$$(\Sig,\alphas,\ov\betas,\{\lambda_0,\lambda_1,\lambda_\infty\};u,v,w,\ubar,\vbar,\wbar)$$
is constructed from a special Heegaard diagram as before. Let $\betas_\bullet$, for
$\bullet\in\{0,1,\infty\}$ denote the set 
$\betas_\bullet=\{\beta_1^\bullet,...,\beta_{g-1}^\bullet,\lambda_\bullet\}$
constructed before. Let us furthermore assume that $\betas_1'$ is a set of 
$g$ simple closed curves, where the first $g-1$ of them are small Hamiltonian
isotopes of the first $g-1$ curves in $\betas_1$ (with two transverse intersection 
points with the corresponding simple closed curve in $\betas_1$)
while the last ($g$th) curve
is denoted by $\lambda_1'$. We  assume that  $\lambda_1'$ is a Hamiltonian 
isotope of  $\lambda_1$, which is very close to the juxtaposition of the curves 
$\lambda_0$  and $\lambda_\infty$.\\
 
 Consider the two Heegaard quadruples 
 $$H=(\Sig,\alphas,\betas_0,\betas_1,\betas_\infty;u,v,w)\ \ \text{and}\ \ 
 H'=(\Sig,\alphas,\betas_0,\betas_1',\betas_\infty;u,v,w).$$      
Let us denote the triangle maps associated with the first Heegaard diagram by 
\begin{displaymath}
\begin{split}
&\pphi_0(H)=\phi(H\setminus \betas_0) \co C_1(K;H)\ra C_\infty(K;H)\ \ \text{and}\\ 
&\pphi_\infty(H)=\phi(H\setminus \betas_\infty) \co C_0(K;H)\ra C_1(K;H),
\end{split}
\end{displaymath}
while the triangle maps associated with the Heegaard quadruple $H'$ are 
denoted by 
\begin{displaymath}
\begin{split}
&\pphi_0(H')=\phi(H'\setminus \betas_0) \co 
C_1(K;H')\ra C_\infty(K;H')=C_\infty(K;H)\ \ \text{and}\\ 
&\pphi_\infty(H')=\phi(H'\setminus \betas_\infty) \co C_0(K;H')=C_0(K;H)\ra C_1(K;H').
\end{split}
\end{displaymath}

The  holomorphic triangle map 
$\pphi_1(H)=\pphi_1(H') \co C_\infty(K;H)\ra C_0(K;H)$ may be defined 
using the Heegaard triple $(\Sig,\alphas,\betas_\infty,\betas_0;u,v,w)$. 
Count of holomorphic rectangles in $H$ and $H'$, respectively, that avoid the punctures 
$u,v$ and $w$ gives the homomorphisms
\begin{displaymath}
\begin{split}
&\Phi_1  \co C_0(K;H)\ra C_\infty(K;H),\ \ \text{and}\ \ \Phi_1' \co C_0(K;H)\ra C_\infty(K;H)
\end{split}
\end{displaymath}
such that 
\begin{displaymath}
\begin{split}
&d_\infty\circ \Phi_1+\Phi_1\circ d_0=\pphi_0(H)\circ \pphi_\infty(H)\ \  \text{and}\ \
d_\infty\circ \Phi_1'+\Phi_1'\circ d_0=\pphi_0(H')\circ \pphi_\infty(H').\\
\end{split}
\end{displaymath}
The interesting observation is that both $\Phi_1$ and $\Phi_1'$ vanish
when the Heegaard diagram is chosen as above. The reason for the first vanishing 
result is that there are no  positive squares connecting the four intersection points
\begin{displaymath}
\begin{split}
&\x\in\Ta\cap\mathbb{T}_{\beta_0},\ \ \Theta_{0,1}\in 
\mathbb{T}_{\beta_0}\cap\mathbb{T}_{\beta_1},\ \ 
\Theta_{1,\infty}\in \mathbb{T}_{\beta_1}\cap\mathbb{T}_{\beta_\infty}
\ \ \text{and}\ \  \y\in\mathbb{T}_{\beta_\infty}\cap\Ta
\end{split}
\end{displaymath}
In fact, for every square class $\square\in\pi_2^+(\x,\Theta_{0,1},\Theta_{1,\infty},\y)$,
 $n_\ubar(\square)=n_\wbar(\square)=1$. Thus, two opposite quarters
around the intersection point $r_1$ have zero coefficient, while one other quadrant 
has coefficient $1$. Since $r_1$ is not among the intersection points in any of 
$\x,\y,\Theta_{0,1}$ and $\Theta_{1,\infty}$ the coefficient of the 
last quadrant around $r_1$ is $-1$, and the contribution of $\square$ is thus trivial. 
A similar argument implies that $\Phi_1'$ is zero.\\

For $\bullet\in\{0,1,\infty\}$ the 
Heegaard triple $H_\bullet=(\Sig,\alphas,\betas_\bullet',\betas_\bullet;u,v,w)$
gives 
$$\imath_\bullet=\imath(H_\bullet) \co C_\bullet(K;H')\ra C_\bullet(K;H).$$
The homomorphisms $\imath_0$ and $\imath_\infty$ are the identity maps of 
$C_0(K;H)$ and $C_\infty(K;H)$ respectively. The Heegaard quadruple 
$$(\Sig,\alphas,\betas_0,\betas_1',\betas_1;u,v,w)$$
determines a holomorphic square map
$$\Psi_\infty \co C_0(K;H)\lra C_1(K;H).$$
Considering different possible degenerations of a holomorphic square of Maslov index 
zero, one finds the relation
\begin{equation}\label{eq:homotopy-3}
d_1\circ \Psi_\infty+\Psi_\infty\circ d_0=\imath_1\circ \pphi_\infty(H')+\pphi_\infty(H).
\end{equation}

Finally, one may consider the Heegaard $5$-tuple 
$$(\Sig,\alphas,\betas_0,\betas_1',\betas_1,\betas_\infty;u,v,w)$$
which may be used to construct a pentagon map
$Q \co C_0(K;H)\ra C_\infty(K;H).$
Considering all possible degenerations of a holomorphic pentagon of Maslov index $-1$,
one obtains the relation
\begin{equation}\label{eq:homotopy-5}
d_\infty\circ Q+Q\circ d_0=\Psi_0\circ \pphi_\infty(H')+
\pphi_0(H)\circ \Psi_\infty,
\end{equation}
where $\Psi_0 \co C_1(K;H')\ra C_\infty(K;H)$ is the holomorphic square map associated with 
$(\Sig,\alphas,\betas_1',\betas_1,\betas_\infty;u,v,w)$.
The reason for the above equality is that the contributing holomorphic squares in 
the Heegaard quadruple $(\Sig,\betas_0,\betas_1',\betas_1,\betas_\infty;u,v,w)$
come in cancelling pairs, while there is a single contributing holomorphic triangle 
corresponding to either of the Heegaard triples 
$$(\Sig,\betas_1',\betas_1,\betas_\infty;u,v,w)\ \ \text{and}\ \ 
(\Sig,\betas_0,\betas_1',\betas_1;u,v,w).$$
Figure~\ref{fig:quadrangle} illustrates the domain for a cancelling pair of 
contributing squares. Moreover, the maps $\Phi_1$ and $\Phi_1'$ which may 
potentially contribute, are trivial.\\

\begin{figure}
\def\svgwidth{15cm}\begin{center}
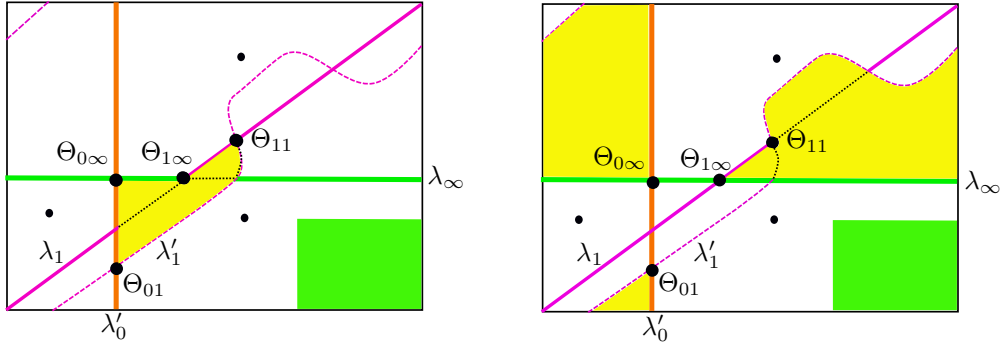
\caption{\label{fig:quadrangle}
{The domains for a cancelling pair of contributing squares connecting 
$\Theta_{01},\Theta_{1,1},\Theta_{1,\infty}$ and $\Theta_{0,\infty}$.
}}\end{center}
\end{figure}
Our choice of $\lambda_1'$ and the fact that the Heegaard diagram is nice 
imply that we have a short exact sequence
\begin{displaymath}
\begin{diagram}
0&\rTo{}&C_0(K;H)&\rTo{\pphi_\infty(H')}&C_1(K;H')&\rTo{\pphi_0(H')}
&C_\infty(K;H)&\rTo &0.
\end{diagram}
\end{displaymath}
Correspondingly, an isomorphism $\theta' \co \Ker(\pphi_\infty(K))\ra \Coker(\pphi_0(K))$
may be constructed.
Choose some closed element $a\in C_0(K;H)$, and let $\pphi_\infty(H')(a)=d_1'(b')$
for some $b'\in C_1(K;H')$. By (\ref{eq:homotopy-3})
\begin{displaymath}
\begin{split}
\pphi_\infty(H)(a)&=\left(\imath_1\circ \pphi_\infty(H')+d_1\circ \Psi_\infty\right)(a)\\
&=d_1\left(\imath_1(b')+\Psi_\infty(a)\right)=:d_1(b)
\end{split}
\end{displaymath} 
Using  (\ref{eq:homotopy-3}) and (\ref{eq:homotopy-5}) we compute
\begin{displaymath}
\begin{split}
\pphi_0(H)(b)&=\pphi_0(H)\left(\imath_1(b')+\Psi_\infty(a)\right)\\
\text{by }(\ref{eq:homotopy-3})\ \ \ \ \ \ \ \ \ \ \ 
&=\pphi_0(H')(b')+\left(d_\infty\circ \Psi_0+\Psi_0\circ d_1'\right)(b')+
\left(\pphi_0(H)\circ \Psi_\infty\right)(a)\\
&=\pphi_0(H')(b')+d_\infty\circ \Psi_0(b')+\left(\Psi_0\circ \pphi_\infty(H')+
\pphi_0(H)\circ \Psi_\infty\right)(a)\\
\text{by }(\ref{eq:homotopy-5})\ \ \ \ \ \ \ \ \ \ \ 
&= \pphi_0(H')(b')+d_\infty\left(Q(a)+\Psi_0(b')\right).
\end{split}
\end{displaymath}
This means that the maps $\theta(K)$ and $\theta'$, as maps from 
$\Ker(\pphi_\infty(K))$ to $\Coker(\pphi_0(K))$, are the same.
However, the map $\theta'$ is the inverse of the connecting homomorphism 
$\delta \co \Coker(\pphi_0(K))\ra \Ker(\pphi_\infty(K))$
resulting from the short exact sequence
\begin{displaymath}
\begin{diagram}
0&\rTo{}&C_0(K;H)&\rTo{\pphi_\infty(H')}&C_1(K;H')&\rTo{\pphi_0(H')}
&C_\infty(K;H)&\rTo &0.
\end{diagram}
\end{displaymath}

The above observations imply the claim for $\theta(K)$. The proof for the map 
$\thetabar(K)$ is similarly reduced  to showing the triviality of the holomorphic 
square map corresponding to the Heegaard quadruple
$(\Sig,\alphas,\betas_0,\betas_1,\betas_\infty;\ubar,\vbar,\wbar)$.\\

The domain of every contributing holomorphic square corresponding 
to the aforementioned punctured Heegaard diagram has coefficient zero at 
$\ubar,\vbar,\wbar$ and $w$, and coefficient $1$ at $u$ and $v$. 
This implies that two opposite 
quadrants around $r_3$  have coefficient zero, while a third quadrant has coefficient 
$1$. Since $r_3$ can not be among the intersection points on the vertices of the 
square, the fourth quadrant around $r_3$ has coefficient $-1$. This 
contradiction gives the triviality of the holomorphic square map corresponding to 
$(\Sig,\alphas,\betas_0,\betas_1,\betas_\infty;\ubar,\vbar,\wbar)$,
and completes the proof.
\end{proof}

\newpage
\section{Gluing the knot complements}\label{sec:gluing}
\subsection{Extracting a chain complex for splicing}\label{subsec:chain-complex}
Given two Heegaard diagrams for the complements of the knots $K_1$ and $K_2$, one may
 construct a Heegaard diagram for $Y(K_1,K_2)$ as follows, similar to the construction of 
 \cite{Ef-Whitehead}. Let
$$H_i=(\Sigma_i,\alphas^i,\ov\betas^i\cup\{\mu_i,\lambda_i\})$$
denote the Heegaard diagram for $K_i$ with Heegaard surface $\Sig_i$, and with
$\mu_i$ the meridian for $K_i$ and $\lambda_i$ a zero framed longitude for it which
cuts $\mu_i$ in a single point $O_i$. Then the Heegaard diagram for 
the three-manifold $Y=Y(K_1,K_2)$ obtained by  splicing the complement of 
$K_1\subset Y_1$
and the complement of ${K_2}\subset {Y_2}$  is
constructed as follows. Attach a one-handle to $\Sigma_1\cup\Sigma_2$, with 
attaching circles placed at the intersections $O_1$ and $O_2$. Use four 
parallel segments on this one-handle to connect the four intersections of 
$\mu_1\cup\lambda_1$ with one of the attaching circles to the four intersections 
of $\mu_2\cup\lambda_2$ with the other attaching circle, so that intersection 
points on $\mu_1$ are joined to the intersection points on $\lambda_2$. 
The union of the remaining parts from $\mu_1$ and $\lambda_2$ with two of the 
four parallel line segments gives a simple closed curve on $\Sig$ which will be denoted 
by $\mu_1\#\lambda_2$. The simple closed curve $\lambda_1\#\mu_2$ is 
constructed in a similar way. Let 
$$\alphas=\alphas^1\cup \alphas^2\ \ \ \text{and}\ \ \ 
\betas=\ov\betas^1\cup\ov\betas^2\cup\{\mu_1\#\lambda_2,\lambda_1\#\mu_2\}.$$
The resulting Heegaard diagram $H=(\Sig,\alphas,\betas)$
is a  Heegaard diagram for the three-manifold 
obtained by splicing the two knot complements.\\

If the initial Heegaard diagrams $H_i$ are special (c.f.
Definition~\ref{defn:Nice}) one may assume that the Heegaard diagram $H$ will have 
one bad region, and the rest 
of the regions are either bigons or rectangles. Thus the combinatorial algorithm of 
\cite{SW} may be used to compute its (hat) Heegaard Floer homology with $\Fbb$ 
coefficients. Let $z$ denote a marked point which is placed in the aforementioned bad region.
The marked point $z$ corresponds to the marked points $z_i\in \Sig_i$, $i=1,2$.
We may also choose a second marked point $z_i'$ for the Heegaard diagram $H_i$
which is placed next to $O_i$ and 
in the quadrant opposite to the quadrant containing $z_i$.  \\

\begin{figure}\def\svgwidth{11cm}
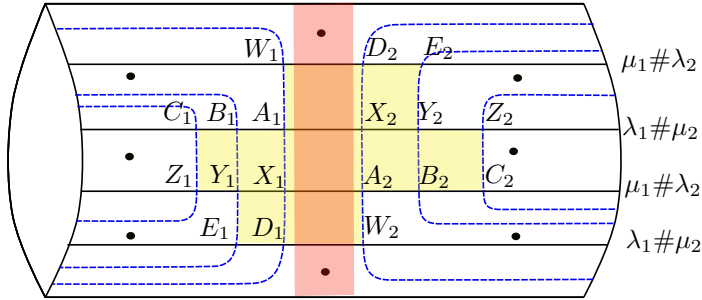
 \caption{\label{fig:local}
{The cylinder illustrates a neighbourhood of the $1$-handle used for 
attaching the two Heegaard diagrams. The union of the domains of the disks 
intersecting the one-handle and contributing to the differential 
is shaded yellow.  }}
\end{figure}

Define the chain complexes $M^i$ and $L^i$ associated with 
$K_i\subset Y_i$ using the Heegaard diagrams 
$$(\Sig_i,\alphas^i,\ov\betas^i\cup \{\mu_i\};z_i,z_i')\ \ \text{and}\ \ 
(\Sig_i,\alphas^i,\ov\betas^i\cup \{\lambda_i\};z_i,z_i'),$$
respectively. Note that the generators of the complex $C$ associated with the 
Heegaard diagram 
$H$ are in correspondence, either with the generators of $M=M^1\otimes M^2$ or the 
generators of $L=L^1\otimes L^2$, i.e. the $\Fbb$-module $C$ may be identified 
with $M\oplus L$.  
Denote the differential of $M$ by $d_M$ and the differential of $L$ by $d_L$.
The domain of every disk which contributes to the differential of $C$ is then a 
 rectangle or a bigon in the diagram. Such a disk may either stay in one of 
$\Sigma_i$, or intersect both $\Sigma_1$ and $\Sigma_2$. 
The disks that stay in one of $\Sigma_i$ correspond to the differentials 
$d_M$ and $d_L$ of the complexes $M$ and $L$. 
Only a few rectangles can intersect
both $\Sig_i$ and miss the marked point $z$ (see Figure~\ref{fig:local}), while no
bigons can intersect both $\Sig_1$ and $\Sig_2$. 
 Because of the way the bad region (the region containing the marked
point) enters the neighbourhood of the one-handle, the rectangles which intersect 
both $\Sig_1$ and $\Sig_2$ stay in the
neighbourhood of the one-handle. The contribution of such rectangles may be described 
after introducing some extra notation.\\

The assumption on the Heegaard diagrams $H_1$ and $H_2$ from 
Lemma~\ref{lem:nice-to-nicer}
implies that the local picture around $O_i$ looks like 
the genus $4$ surface illustrated on the top of Figure~\ref{fig:Nice}.
Denote the intersection points on $H_i$ which correspond to $A,B,C,D,E,X,Y,Z$ and $W$ by 
$A_i,B_i,C_i,D_i,E_i,X_i,Y_i,Z_i$ and $W_i$, respectively.\\

The generators of $M^i\oplus  L^i$ are the tuples 
$\x=(x_1,...,x_{g_i})$ such that for a permutation $\sigma \co \{1,...,g_i\}\ra\{1,...,g_i\}$,
$x_j\in \alpha_{\sigma(j)}\cap \beta_{j}$, $j=1,...,g_i-1$ and
$x_{g_i}\in \alpha_{\sigma(g_i)}\cap(\mu_i\cup\lambda_i)$. 
The complex $M^i$ is generated by those $\x$ such that $x_{g_i}\in \mu_i$, 
and the complex $L^i$ is generated by the $g_i$-tuples $\x=(x_1,...,x_{g_i})$ with
$x_{g_i}\in \lambda_i$. The homology of
the complex $M^i$ is the knot Floer homology $\widehat{\text{HFK}}(K_i)$, 
and the homology of the complex $L^i$ is the 
longitude Floer homology $\widehat{\text{HFL}}(K_i)$. 
The homomorphisms $\Phi^i \co M^i\rightarrow L^i$ over 
$\x=(x_1,\dots,x_{g_i})\in M^i$ are
defined by 
\begin{displaymath}
\Phi^i\left(\x\right)=\begin{cases}
(x_1,...,x_{g_i-1},X_i)\ \ \ &\text{if }x_{g_i}=A_i\\
(x_1,...,x_{g_i-1},Y_i)\ \ \ &\text{if }x_{g_i}=B_i\\
(x_1,...,x_{g_i-1},Z_i)\ \ \ &\text{if }x_{g_i}=C_i\\
0 &\text{otherwise}
\end{cases}.
\end{displaymath}
The corresponding contributing triangles  are  $[A_iO_iX_i],[B_iO_iY_i]$, and $[C_iO_iZ_i]$. 
The map $\Phi$ thus  corresponds to the  changes $x_{g_i}\rightarrow y_{g_i}$
which are one of the following: $A_i\ra X_i$, $B_i\ra Y_i$ or $C_i\ra Z_i$.
Similarly, the homomorphisms  $\Psi^i_1 \co L^i\rightarrow M^i$ 
 correspond to the triangles
$[W_iO_iA_i]$, and over $\x=(x_1,...,x_{g_i})\in L^i$ are defined by
\begin{displaymath}
\Psi^i_1\left(\x\right)=\begin{cases}
(x_1,...,x_{g_i-1},A_i)\ \ \ &\text{if }x_{g_i}=W_i\\
0 &\text{otherwise}
\end{cases}.
\end{displaymath}
Define the  maps $\Psi_2^i,\Psi_3^i \co L^i\ra M^i$, where $\Psi^i_2$ corresponds to the
changes $X_i\ra D_i$ and $Y_i\ra E_i$, and $\Psi^i_3$ corresponds to $W_i\ra D_i$. 
Thus the triangles contributing to $\Psi^i_2$ are $[X_iO_iD_i]$ and $[Y_iO_iE_i]$,
while the only triangle contributing to $\Psi^i_3$ is 
$$[W_iO_iD_i]=[W_iO_iA_i]\cup[A_iO_iX_i]\cup[X_iO_iD_i].$$

The contribution of the rectangles which intersect both $\Sig_1$ and $\Sig_2$  
to the differential of the complex 
$C=M\oplus L$ may thus be described by the maps 
\begin{displaymath}
\begin{split}
\Phi=\Phi^1\otimes \Phi^2 & \co L^1\otimes L^2\lra M^1\otimes M^2,\\
\begin{array}{c}
\Psi_1=\Psi^1_1\otimes \Psi^2_2\\
\Psi_2=\Psi^1_2 \otimes \Psi^2_1\\
\Psi_3=\Psi^1_3\otimes \Psi^2_3
\end{array}
 & \co M^1\otimes M^2 \lra L^1\otimes L^2.\\
\end{split}
\end{displaymath}
In other words, the differential of the complex $C=M\oplus L$ is   the homomorphism
$$d=d_C=\left(\begin{array}{cc}d_M&\Phi \\
\sum_{i=1}^3 \Psi_i& d_L\end{array}\right).$$

\begin{prop}\label{prop:ML-mapping-cone}
The complexes $M^i$ and $L^i$ are identified with the mapping cones of 
$\fbar_\infty^i=\fbar_\infty(K_i)$
and $f_0^i=f_0(K_i)$ respectively. 
More precisely, the $\Fbb$-module $M^i$ is isomorphic to  the direct 
sum of $C_1(K_i)$ and $C_0(K_i)$, while $L^i$ is isomorphic to the direct sum of 
$C_\infty(K_i)$ and $C_1(K_i)$. 
Moreover, the differentials $d_{M^i}$ and $d_{L^i}$ of $M^i$ and $L^i$ are identified as
\begin{displaymath}
\begin{split}
&d_{L_i}(c_1,c_\infty)=(d_1^i(c_1),d_\infty^i(c_\infty)+f_0^i(c_1))\ \ \ \ \ \ 
\forall\ (c_1,c_\infty)\in C_1(K_i)\oplus C_\infty(K_i)\ \ \text{and}\\
&d_{M_i}(c_0,c_1)=(d_0^i(c_0),d_1^i(c_1)+\fbar_\infty^i(c_0))\ \ \ \ \ \ \ \ \ 
\forall\ (c_0,c_1)\in C_0(K_i)\oplus C_1(K_i).
\end{split}
\end{displaymath}
\end{prop}

\begin{proof}
We sketch the proof of the claim for $L^i$. The corresponding claim for $M^i$ is proved 
in a completely similar way.
Consider the labelling of the intersection points of the $\alpha$-curves 
with  the curves $\lambda_\infty(K_i),\lambda_1(K_i)$ and $\lambda(K_i)$ as in 
Figure~\ref{fig:H1-and-H0}. The intersection points with the 
$\alpha$-curves on $\lambda_1(K_i)$ are  $r_1,r_2$ and $r_3$. The intersection 
points with the $\alpha$-curves on $\lambda_\infty(K_i)$ are
$s_1,s_2,...,s_n,q_3,q_4,...,q_n$, and the intersection points with the $\alpha$-curves
 on $\lambda(K_i)$ are $S_1,S_2,...,S_n,Q_3,Q_4,...,Q_n$, and $R_1,R_2,R_3$. Define the
 $\Fbb$-module  isomorphism 
\begin{displaymath}
 I_i \co C_1(K_i)\oplus C_\infty(K_i)\ra L^i,\ \ \ 
  I_i\left(\x=(x_1,...,x_{g_i})\right):=(x_1,...,x_{g_i-1},I_i(x_{g_i})),
\end{displaymath}  
  where $I_i$  changes the letter in the 
  labeling of an intersection point to a capital letter (so 
$I_i(r_j)=R_j, I_i(s_j)=S_j$ and $I_i(q_j)=Q_j$). Straightforward combinatorics may be used to 
verify $d_{L^i}(I_i(\x))=I_i(d_\infty^i(\x))$ 
for every generator $\x$ of $C_\infty(K_i)$, and   $d_{L^i}(I_i(\x))=I_i(d_1^i(\x))+I_i(f_0^i(\x))$
for every generator $\x$ of $C_1(K_i)$. 
\end{proof}

Under the identification of $M^i$ with the mapping cone of $\fbar_\infty^i$, and the
identification of $L^i$ with the mapping cone of $f_0^i$,  the map
$\Phi$ has a simple description; it is the map that takes $C_1(K_i)$ in the mapping
cone of $f_0 \co C_1(K_i)\ra C_\infty(K_i)$ to the complex $C_1(K_i)$ in the mapping cone of
$\fbar_\infty^i \co C_0(K_i) \ra C_1(K_i)$ via the identity map of $C_1(K_i)$. 
Furthermore, the map 
$f_\infty^i$ from $C_0(K_i)$
in $M^i$ to $C_1(K_i)$ in $L^i$ is  identified with the triangle map $\Psi_1^i$. 
The induced map $\fbar_0^i$ from the copy of $C_1(K_i)$ in $M$ to the copy 
of $C_\infty(K_i)$ in $L^i$ is the 
triangle map $\Psi_2^i$. The map $\Psi_3^i$ is obtained from the composition map
$\fbar_0^i\circ f_\infty^i \co  C_0(K_i)\ra C_\infty(K_i)$.
Set $C_\bullet^i=C_\bullet(K_i)$.
If we replace the mapping cone of $C_1^i\xra{f_0^i}C_\infty^i$ for $L^i$,
replace the mapping cone $C_0^i\xra{\fbar_\infty^i}C_1^i$ for  $M^i$, and also replace 
$\Phi^i$ and $\Psi^i_j$ with the 
appropriate descriptions in terms of $\fbar_0^i$ and $f_\infty^i$,
we obtain an alternative description of the complex $C$.\\

 The {\emph{cube}} 
$\M=\M(f_\bullet^i,\fbar_\bullet^i\ |\  \bullet=0,\infty,\ i=1,2)$
associated with the knots 
$K_1$ and $K_2$, the corresponding complexes $C^i_\bullet$, $i=1,2,\ 
\bullet\in\{0,1,\infty\}$ and the maps $f_0^i,\fbar_0^i \co C_1^i\ra C_\infty^i$ and 
the maps $f_\infty^i,\fbar_\infty^i \co C_0^i\ra C_1^i$  is the chain complex 
$(\M,d_\M)$ associated with the graph of complexes represented by the 
 following cube:
\begin{equation}\label{eq:complex-cube}
\begin{diagram}
\Cinfinf   &              & &\lTo^{I\otimes \f_0^2}&             & \Cinfone       &               &          \\
	   &\luTo^{\f_0^1\otimes I}& && & \uTo        &\luTo^{\f_0^1\otimes I} &          \\
&& \Coneinf     &&\lTo^{I\otimes \f_0^2}& \HonV &&\Coneone  \\
      \uTo^{\fbar_0^1\f_\infty^1\otimes \fbar_0^2\f_\infty^2}  
      &              &              &&             & \vLine_{\fbar_0^1\otimes\f_\infty^2}   &               &          \\
	   &              & \uTo^{\f_\infty^1\otimes \fbar_0^2}  &&             &                &          &\dTo^{I}\\
\Czerozero & \hLine & \VonH &\rTo{\fbar_\infty^1\otimes I}&& \Conezero      && \\
	  & \rdTo_{I\otimes \fbar_\infty^2}&              && &  &\rdTo_{I\otimes \fbar_\infty^2} &          \\
	   &              &\Czeroone     &&\rTo^{\fbar_\infty^1\otimes I}  &  &   &\Coneone  \\
\end{diagram}.
\end{equation}

\begin{prop}\label{prop:complex-identification}
With the above notation fixed, the 
complex $(C,d)$  is  identified, as a chain complex, with the cube 
$$\left(\M=\M(f_\bullet^i,\fbar_\bullet^i\ |\  \bullet=0,\infty, \ i=1,2),d_\M\right).$$ 
\end{prop}

\subsection{The linear algebra of the cubes}\label{subsec:linear-algebra}
Let  $\Hbb^i_\bullet$ denote the homology of the chain complex 
$(C^i_\bullet,d^i_\bullet)$ 
for $i=1,2$ and $\bullet\in\{0,1,\infty\}$. Set 
$\Hbb_{\bullet,\star}=\Hbb_\bullet^1\otimes\Hbb_\star^2$ for 
$\bullet,\star\in\{0,1,\infty\}$. 
Abusing the notation, the map induced on homology by $f^i_\bullet$  
will also be denoted $\pphi^i_\bullet$ and  the map induced on homology 
by $\fbar_\bullet^i$ will be denoted $\pphibar_\bullet^i$.\\

Following the discussion of Subsection~\ref{subsec:simplifications}, we may choose 
 appropriate decompositions  
$C^i_\bullet=A^i_\bullet\oplus \Hbb^i_\bullet\oplus A^i_\bullet$ such that 
the differential 
$d^i_\bullet$
takes the form $d^i_\bullet=\left(\begin{array}{ccc}
0&0&I\\ 0&0&0\\ 0&0&0\end{array}\right)$. Correspondingly, we find the matrices
$G(f_\bullet^i)=(f_\bullet^i)_*$ and $G(\fbar_\bullet^i)=(\fbar_\bullet^i)_*$ 
which will be denoted by $\pphi_\bullet^i$ and $\pphibar_\bullet^i$, as well as 
the matrices
$$M(f_\bullet^i),M(\fbar_\bullet^i),P(f_\bullet^i),P(\fbar_\bullet^i),
Q(f_\bullet^i),Q(\fbar_\bullet^i),N(f_\bullet^i)\ \ \text{and}\ \  N(\fbar_\bullet^i).$$

 The maps $Q(\pphi_0^i)P(\pphi_\infty^i)$ and 
 $Q(\pphibar_0^i)P(\pphibar_\infty^i)$ 
 from $\Hbb^i_0$ to $\Hbb^i_\infty$ extend the homomorphisms 
$$\theta^i \co \Ker(\pphi_\infty^i)\ra \Coker(\pphi_0^i)\ \ \ \text{and}\ \ \ 
\thetabar^i \co \Ker(\pphibar_\infty^i)\ra \Coker(\pphibar_0^i),$$ 
associated with the knot $K_i\subset Y_i$. These extensions are still denoted by 
$\theta^i$ and $\thetabar^i$, respectively.\\

Lemma~\ref{lem:G-complex} implies that the homology of $(\M,d_\M)$ is isomorphic to 
the homology of the chain complex $(\Hbb,d_\Hbb)$ associated with the 
graph of chain complexes determined by the cube of Figure~\ref{fig:cube}.

\begin{prop}\label{prop:cube}
Let $(\Hbb,d_\Hbb)$ denote the complex obtained from the cube 
$(\M,d_\M)$ by applying Lemma~\ref{lem:G-complex}. Then $(\Hbb,d_\Hbb)$
is identified with the complex shown in Figure~\ref{fig:cube}, provided that the maps 
$\theta^i=\theta(K_i)$ and $\thetabar^i=\thetabar(K_i)$ are given as above. 
\end{prop}
\begin{proof}
If Lemma~\ref{lem:G-complex} is applied, we obtain the same oriented graph
(i.e. the same new edges) and the same complexes on the vertices. The directed edge 
from $\Hzerozero$ to $\Hinfinf$ is labelled by the map
\begin{displaymath}
\begin{split}
\pphibar_0^1 \pphi_\infty^1\otimes\pphibar_0^2\pphi_\infty^2+
Q(f_0^1\otimes I)N(f_\infty^1\otimes \fbar_0^2)P(I\otimes \fbar_\infty^2)+
Q(I\otimes f_0^2)N(\fbar_0^1\otimes f_\infty^2)P(\fbar_\infty^1\otimes I)
\end{split}
\end{displaymath}
which is, by Lemma~\ref{lem:theta-2}, equal to 
\begin{displaymath}
\pphibar_0^1 \pphi_\infty^1\otimes\pphibar_0^2\pphi_\infty^2
+\theta^1\otimes \thetabar^2
+\thetabar^1\otimes \theta^2.
\end{displaymath}
The map corresponding to the dashed edge from $\Hzerozero$ to $\Honeinf$ is,
by Lemma~\ref{lem:theta-1} 
\begin{displaymath}
Q(f_\infty^1\otimes \fbar_0^2)P(I\otimes \fbar_\infty^2)
=\pphi_\infty^1\otimes \left(Q(\fbar_0^2)P(\fbar_\infty^2)\right)
=\pphi_\infty^1\otimes \thetabar^2.
\end{displaymath}
The maps corresponding to the rest of dashed directed edges may be computed in 
a completely similar way. 
 This completes the proof of Proposition~\ref{prop:cube}.
\end{proof}

\begin{remark} We would like to make two comments:\\
(1) Note that $Y(K_1,K_2)=Y(-K_1,-K_2)$. One may assume that 
$\pphi_\bullet(-K)=\pphibar_\bullet(K)$ and
$\pphibar_\bullet(-K)=\pphi_\bullet(K)$, implying that 
$\ov{\HFT}(Y(-K_1,-K_2))$ is isomorphic to
  the homology of the complex determined by the oriented 
  graph in Figure~\ref{fig:cube}, where all barred maps change to the corresponding 
  unbarred maps, and all unbarred maps change to the corresponding barred maps.\\
(2) Proposition~\ref{prop:cube} is still weaker than Theorem~\ref{thm:cube}, since the 
extensions of $\theta^i$ and $\thetabar^i$ to maps from $\Hbb_0^i$ to $\Hbb_\infty^i$
are not arbitrary yet. In fact, without freedom in choosing these two extensions (which will
be proved by the end of the current section) Theorem~\ref{thm:cube} stays bound to the 
information from the corresponding nice Heegaard diagram, and has much less significance.
\end{remark}

\subsection{Simplifications of the splicing formula}\label{subsec:simplifications}
We now apply Lemma~\ref{lem:simplification} to the splicing formula of  
Proposition~\ref{prop:cube} and make some cancellations. The first cancellation 
comes from setting $C=\Hbb$, $A=\Honeone$ and 
$$B=\left(\Hinfinf\oplus \Honeinf\oplus \Hinfone\right)\oplus
\left( \Hzeroone\oplus \Honezero\oplus \Hzerozero\right)=\E_1\oplus \E_2.$$
We thus have 
$\ov\HFT (Y)=H_*(B,d_B)$, where
\begin{displaymath}
d_B=\colvec{
0&\pphi_0^1\otimes I&I\otimes \pphi_0^2&\theta^1\otimes \pphibar_0^2&
\pphibar_0^1\otimes \theta^2&\begin{array}{c}
\pphibar_0^1\circ \pphi_\infty^1)\otimes
\pphibar_0^2\circ \pphi_\infty^2+\\
\theta^1\otimes \thetabar^2+\thetabar^1\otimes \theta^1
\end{array}\\
&&&\\
0&0&0&\Phi&I\otimes (\pphi_0^2\circ \pphibar_\infty^2)&
\pphi_\infty^1\otimes\thetabar^2\\
0&0&0&(\pphi_0^1\circ \pphibar_\infty^1)\otimes I&\Psi &
\thetabar^1\otimes \pphi_\infty^2\\
0&0&0&0&0&I\otimes \pphibar_\infty^2\\
0&0&0&0&0&\pphibar_\infty^1\otimes I\\
0&0&0&0&0&0\\
},\end{displaymath}
with $\Phi=\pphibar_\infty^1\otimes \pphi_0^2+\pphi_\infty^1\otimes \pphibar_0^2$ and 
$\Psi=\pphi_0^1\otimes \pphibar_\infty^2+\pphibar_0^1\otimes \pphi_\infty^2$.\\

The dimension of the $\Fbb$-vector space $H_*(B,d_B)$ only depends on 
the rank of the kernel and the cokernel of the matrix $d_B$. Define a pair of matrices 
$M_1$ and $M_2$ equivalent if $\Ker(M_1)\simeq \Ker(M_2)$ and 
$\Coker(M_1)\simeq \Coker(M_2)$. For a matrix $M$ let 
$$\imath(M):=\Ker(M)\oplus \Coker(M)\ \ \ \text{and}\ \ \ i(M):=\rank(\imath(M)).$$
If $M_1$ and $M_2$ are equivalent matrices then $\imath(M_1)\simeq \imath(M_2)$
and $i(M_1)=i(M_2)$.\\
 
We make a change of basis for $\E_{2}$ which is given by the matrix 
\begin{displaymath}
\colvec{\tau_0(K_1)\otimes \tau_1(K_2)&0&0\\
0&\tau_1(K_1)\otimes \tau_0(K_2) &0\\
0&0& \tau_0(K_1)\otimes \tau_0(K_2)
}.\end{displaymath}
The matrix $d_B$ is thus equivalent to the matrix
\begin{displaymath}
d_B'=\colvec{
0&\pphi_0^1\otimes I&I\otimes \pphi_0^2&
\theta^1\tau_0^1\otimes \tau^2_\infty\pphi_0^2
&\tau^1_\infty\pphi_0^1\otimes\theta^2\tau_0^2
&\begin{array}{c}
\tau_\infty^1\pphi_0^1\pphibar_\infty^1\otimes
\tau_\infty^2\pphi_0^2\pphibar_\infty^2\\
+\theta^1\tau_0^1\otimes \thetabar^2\tau_0^2\\
+\thetabar^1\tau_0^1\otimes \theta^2\tau_0^2
\end{array}
\\
&&&\\
0&0&0&\begin{array}{c}
\tau_1^1\pphi_\infty^1\otimes \pphi_0^2\tau_1^2+\\
\pphi_\infty^1\tau_0^1\otimes \tau_\infty^2\pphi_0^2
\end{array}
&\tau_1^1\otimes \pphi_0^2\tau_1^2\pphi_\infty^2 &
\pphi_\infty^1\tau_0^1\otimes\thetabar^2\tau_0^2\\
&&&\\
0&0&0&\pphi_0^1\tau_1^1\pphi_\infty^1\otimes \tau_1^1 &
\begin{array}{c}
\pphi_0^1\tau_1^1\otimes \tau_1^2\pphi_\infty^2+\\
\tau_\infty^1\pphi_0^1\otimes \pphi_\infty^2\tau_0^2
\end{array}
&\thetabar^1\tau_0^1\otimes \pphi_\infty^2\tau_0^2\\
&&&\\
0&0&0&0&0&I\otimes \pphi_\infty^2\\
0&0&0&0&0&\pphi_\infty^1\otimes I\\
0&0&0&0&0&0\\
}.
\end{displaymath}
Let us use the decompositions of (\ref{eq:decompositions}) for $K_1$ and $K_2$
to obtain  a $24\times 24$ block decomposition of $d_B'$. 
Moreover, following the discussion at the end of Subsection~\ref{subsec:properties}
we may assume that in the corresponding decompositions 
$$\theta^i=\left(\begin{array}{cc}
0&I\\ 0&0
\end{array}\right) \ \ \ \text{and}\ \ \ 
\tau_\infty^i\thetabar^i\tau_0^i=\left(\begin{array}{cc}
M^i&I\\ P^iM^i&P^i
\end{array}\right).$$
Each entry in the above $6\times 6$ decomposition for $d_\B'$ 
corresponds to a $4\times 4$ sub-matrix of the aforementioned $24\times 24$
decomposition. For instance the $(1,4)$ entry 
$\theta^1\tau_0^1\otimes \tau_\infty^2\pphi_0^2$ 
corresponds to 
\begin{displaymath}
\colvec{0&I\\ 0&0}\colvec{A_0^1&B_0^1\\ C_0^1& D_0^1}\otimes 
\colvec{A_\infty^2& B_\infty^2\\ C_\infty^2& D_\infty^2}\colvec{0&0\\ I&0}=
\colvec{C_0^1& D_0^1\\0&0}\otimes \colvec{B_\infty^2&0\\ D_\infty^2&0}=
\colvec{
C_0^1\otimes B_\infty^2& 0& D_0^1\otimes B_\infty^2&0\\
C_0^1\otimes D_\infty^2& 0& D_0^1\otimes D_\infty^2 &0\\
0&0&0&0\\
0&0&0&0
}.
\end{displaymath}
For another instance, note that the $(3,5)$ entry  corresponds to 
\begin{displaymath}
\colvec{0&0&0&0\\
0&0&0&0\\
A_1^1\otimes B_1^2&0 & B_1^1\otimes B_1^2&0\\
A_1^1\otimes D_1^2&0 & B_1^1\otimes D_1^2&0 }+
\colvec{0 &0 &0&0\\
B_\infty^1\otimes A_0^2&B_\infty^1\otimes B_0^2&0&0\\
0 &0 &0&0\\
D_\infty^1\otimes A_0^2&D_\infty^1\otimes B_0^2&0&0
}=
\colvec{0&0&0&0\\
B_\infty^1\otimes A_0^2&B_\infty^1\otimes B_0^2&0&0\\
A_1^1\otimes B_1^2&0 & B_1^1\otimes B_1^2&0\\
&&&\\
\begin{array}{c}A_1^1\otimes D_1^2\\+D_\infty^1\otimes A_0^2\end{array}&
D_\infty^1\otimes B_0^2 & B_1^1\otimes D_1^2&0 
}.
\end{displaymath}
The aforementioned $24\times 24$ decomposition
 includes identity matrices as the entries determined by the following block coordinates:
\begin{displaymath}
(2,9), \ \ (3,5),\  \ (4,6), \ \  (14,21),\ \  (16,23)\ \ \text{and}\ \ (20,22).
\end{displaymath}
We use the above $6$ identity matrices for cancellation to obtain an equivalent 
matrix $d=\left(\begin{array}{cc}
0& D\\ 0&0\end{array}\right)$ 
over $\B_1\oplus \B_2$, where 
 $\A_{\bullet \star}=\A_\bullet(K_1)\otimes \A_\star(K_2)$ and
\begin{displaymath}
\begin{split}
&\B_1=\A_{11}\oplus \A_{\infty 1}\oplus \A_{\infty 0}\oplus \A_{1 \infty}\oplus
\A_{0 \infty}\oplus \A_{0 0}\ \ \text{and}\\
&\B_2=\A_{\infty 0}\oplus \A_{1 0}\oplus \A_{\infty \infty}\oplus \A_{0 \infty}\oplus
\A_{0 1}\oplus \A_{1 1}.
\end{split}
\end{displaymath}
Re-arrange the rows and the columns of the matrix $D$, so that $D$ corresponds to the 
rows $11,7,8,10,12,1$ and the columns $19,13,15,17,18,24$ in the above 
$24\times 24$ decomposition to obtain the  following matrix
\begin{displaymath}
\colvec{
B_1^1\otimes B_1^2&B_1^1\otimes A_1^2&0&A_1^1\otimes B_1^2&0&0\\
0&A_0^1\otimes B_\infty^2 &B_0^1\otimes B_\infty^2 &0&0&B_0^1\otimes 
(A_\infty^2+ B_\infty^2 P^2)\\
&&&&&\\
D_1^1\otimes B_1^2&
\begin{array}{c}
D_1^1\otimes A_1^2\\
+A_0^1\otimes D_\infty^2
\end{array}
& B_0^1\otimes D_\infty^2& C_1^1\otimes B_1^2&0&
B_0^1\otimes (C_\infty^2+ D_\infty^2 P^2)\\
&&&&&\\
0&0&0& B_\infty^1\otimes A_0^2   &   B_\infty^1\otimes B_0^2  
 &(A_\infty^1+B_\infty^1 P^1)\otimes B_0^2 \\
&&&&&\\
B_1^1\otimes D_1^2&B_1^1\otimes C_1^2&0  &
\begin{array}{c}
D_\infty^1\otimes A_0^2\\
+A_1^1\otimes D_1^2
\end{array}& D_\infty^1\otimes B_0^2
&(C_\infty^1+ D_\infty^1 P^1) \otimes B_0^2\\
&&&&&\\
0&C_0^1\otimes B_\infty^2& D_0^1\otimes B_\infty^2&
B_\infty^1\otimes C_0^2&B_\infty^1\otimes D_0^2& 
\begin{array}{c}
B_\infty^1 B_1^1 B_0^1\otimes B_\infty^2 B_1^2 B_0^2+\\
(A_\infty^1+ B_\infty^1 P^1)\otimes D_0^2+\\
D_0^1\otimes (A_\infty^2+ B_\infty^2 P^2)
\end{array}}.
\end{displaymath}
This matrix is in turn equivalent to the matrix $\Dd=\Dd(K_1,K_2)$ below,
which is obtained by adding $I\otimes P^2$ times the third column and $P^1\otimes I$
times the fifth column to the last column of the above matrix:
\begin{displaymath}
\Dd=\colvec{
B_1^1\otimes B_1^2&C_1^1\otimes A_1^2&0&A_1^1\otimes B_1^2&0&0\\
0&A_0^1\otimes B_\infty^2 &B_0^1\otimes B_\infty^2 &0&0&B_0^1\otimes A_\infty^2\\
&&&&&\\
D_1^1\otimes B_1^2&
\begin{array}{c}
D_1^1\otimes A_1^2\\
+A_0^1\otimes D_\infty^2
\end{array}
& B_0^1\otimes D_\infty^2& C_1^1\otimes B_1^2&0&B_0^1\otimes C_\infty^2\\
&&&&&\\
0&0&0& B_\infty^1\otimes A_0^2   &   B_\infty^1\otimes B_0^2  
 &A_\infty^1\otimes B_0^2 \\
&&&&&\\
B_1^1\otimes D_1^2&B_1^1\otimes C_1^2&0  &
\begin{array}{c}
D_\infty^1\otimes A_0^2\\
+A_1^1\otimes D_1^2
\end{array}& D_\infty^1\otimes B_0^2
&C_\infty^1 \otimes B_0^2\\
&&&&&\\
0&C_0^1\otimes B_\infty^2& D_0^1\otimes B_\infty^2&
B_\infty^1\otimes C_0^2&B_\infty^1\otimes D_0^2& 
\begin{array}{c}
A_\infty^1\otimes D_0^2\\+D_0^1\otimes A_\infty^2\\
+X^1\otimes X^2
\end{array}},
\end{displaymath}
where $X^i=X(K_i)=B_\infty^i B_1^i B_0^i$ for $i=1,2$. \\

Combining Proposition~\ref{prop:cube} with the above observations we find:
\begin{prop}\label{prop:splicing-2}
Let $K_i\subset Y_i$ denote null-homologous knots for $i=1,2$ and $Y(K_1,K_2)$ 
denote the three-manifold obtained by splicing the complement of $K_1$
with the complement of $K_2$. With the above definition of $\Dd(K_1,K_2)$
\begin{displaymath}
\begin{split}
\ov{\HFT}\left(Y(K_1,K_2),\Fbb\right)
&\simeq \imath\left(\Dd(K_1,K_2)\right)\\
\end{split}
\end{displaymath}
\end{prop}

\begin{cor}\label{cor:theta-independence}
The splicing formula of Proposition~\ref{prop:cube} is independent 
of the choice of extensions $\theta^i$ and $\thetabar^i$.
\end{cor}
\begin{proof}
The fact that the matrices $P^i$ and $M^i$ do not appear in the matrix 
$\Dd(K_1,K_2)$ implies that the choice of the extensions 
$\theta^i,\thetabar^i \co \Hbb_0^i\ra \Hbb_\infty^i$ does not change the rank of the 
homology group in the splicing formula of 
Proposition~\ref{prop:cube} or
Theorem~\ref{thm:cube}.
\end{proof}
With the above corollary in place, the proof of Theorem~\ref{thm:cube} is 
now complete.

\begin{defn}
 For a pair of knots $K_i\subset Y_i$, $i=1,2$ define 
\begin{displaymath}
\begin{split}
\chi(K_1,K_2):&=\Big(h_1(K_1)-h_\infty(K_1)\Big)\Big(h_1(K_2)-h_\infty(K_2)\Big)\\
&\ \ \ \ \ \ \ \ \ \ \ \ 
-\Big(h_1(K_1)-h_0(K_1)\Big)\Big(h_1(K_2)-h_0(K_2)\Big)\\
\end{split}
\end{displaymath}
\end{defn}
Note that $\chi(K_1,K_2)$ is in fact the difference between the ranks of $\B_1=\B_1(K_1,K_2)$
and $\B_2=\B_2(K_1,K_2)$. In the corresponding $\Z/2\Z$-grading on $\B_1\oplus \B_2$,
$\chi(K_1,K_2)$ is thus the Euler characteristic of the chain 
complex $(\B_1\oplus \B_2,d)$. 
\begin{cor}\label{cor:bound-on-rank}
With the above notation fixed, 
\begin{displaymath}
\rank\left(\ov\HFT(Y(K_1,K_2))\right)\geq 
\left|\chi(K_1,K_2)\right|.
\end{displaymath}
\end{cor}
\begin{proof}
It is enough to note that
\begin{displaymath}
\chi(K_1,K_2)
=\rank\big(\Ker(\Dd(K_1,K_2))\big)-\rank\big(\Coker(\Dd(K_1,K_2))\big).
\end{displaymath}
\end{proof}

Consider the matrices
\begin{displaymath}
P_L=\colvec{
I\otimes A_1^2&0&0&0&I\otimes B_1^2&0\\
0&I\otimes A_\infty^2&I\otimes B_\infty^2&0&0&0\\
0&I\otimes C_\infty^2&I\otimes D_\infty^2&0&0&0\\
0&0&0&I\otimes A_0^2&0&I\otimes B_0^2\\
I\otimes C_1^2&0&0&0&I\otimes D_1^2&0\\
0&0&0&I\otimes C_0^2&0&I\otimes D_0^2\\
}
\end{displaymath}
and
\begin{displaymath}
P_R=\colvec{
D_1^1\otimes I&0&0&C_1^1\otimes I&0&0\\
0&A_0^1\otimes I&B_0^1\otimes I&0&0&0\\
0&C_0^1\otimes I&D_0^1\otimes I&0&0&0\\
B_1^1\otimes I&0&0&A_1^1\otimes I&0&0\\
0&0&0&0&D_\infty^1\otimes I& C_\infty^1\otimes I\\
0&0&0&0&B_\infty^1\otimes I& A_\infty^1\otimes I\\
}.
\end{displaymath}
Since $P_R^2=P_L^2=Id$, both $P_R$ and $P_L$ are invertible and
 $\Dd(K_1,K_2)$ is  equivalent to  $\Dd'(K_1,K_2)=P_L\Dd(K_1,K_2) P_R$.
The matrix $\Dd'(K_1,K_2)$ has the following block presentation.
\begin{displaymath}
\colvec[.7]{
D_\infty^1B_1^1\otimes B_1^2A_0^2&B_1^1A_0^1\otimes I&
B_1^1B_0^1\otimes I&D_\infty^1 A_1^1\otimes B_1^2A_0^2&I\otimes B_1^2B_0^2&0\\
&&&\\
I\otimes B_\infty^2 B_1^2 &D_1^1A_0^1\otimes B_\infty^2 A_1^2 &
D_1^1B_0^1\otimes B_\infty^2 A_1^2& 0&B_0^1B_\infty^1\otimes I&
B_0^1 A_\infty^1\otimes I\\
&&&\\
I\otimes D_\infty^2 B_1^2& \begin{array}{c}I\otimes I+\\
D_1^1 A_0^1\otimes D_\infty^2 A_1^2\end{array}& D_1^1 B_0^1\otimes D_\infty^2 A_1^2 
&0&0&0\\
&&&\\
B_\infty^1 B_1^1\otimes I&0&I\otimes B_0^2B_\infty^2& B_\infty^1 A_1^1\otimes I
&\begin{array}{c} D_0^1 B_\infty^1\otimes B_0^2A_\infty^2\\ 
+X^1B_\infty^1\otimes B_0^2X^2\end{array}
&\begin{array}{c} D_0^1 A_\infty^1\otimes B_0^2A_\infty^2\\ 
+X^1A_\infty^1\otimes B_0^2X^2\end{array}\\
&&&\\
D_\infty^1 B_1^1\otimes D_1^2A_0^2&0&0&\begin{array}{c}I\otimes I+\\
D_\infty^1 A_1^1\otimes D_1^2A_0^2\end{array}&I\otimes D_1^2B_0^2&0\\
&&&\\
0&0&I\otimes D_0^2B_\infty^2 &0&\begin{array}{c} 
D_0^1 B_\infty^1\otimes D_0^2A_\infty^2\\
+X^1B_\infty^1\otimes D_0^2X^2\end{array}&\begin{array}{c}
I\otimes I+\\ D_0^1A_\infty^1\otimes D_0^2 A_\infty^2\\
+X^1 A_\infty^1\otimes D_0^2X^2\end{array}
},
\end{displaymath}
and is easier to use in actual computations. Note that 
$$\imath(\Dd'(K_1,K_2))\simeq \imath (\Dd(K_1,K_2))\simeq \ov\HFT(Y(K_1,K_2),\Fbb).$$
\newpage

\section{Splicing with the trefoil}\label{sec:example}
\subsection{The maps $\pphi_\bullet$ and $\pphibar_\bullet$  for the trefoils}
Let us now consider the case of the right handed trefoil, 
which will be denoted by $R$. Thus 
$h_\infty(R)=h_1(R)=3$ and $h_0(R)=4$.
Moreover, $y_\infty(R)=y_1(R)=1$, while $y_0(R)=2$ (c.f. Section 5 of \cite{Ef-Csurgery}). 
Since $\Hbb_\bullet(R,i)=\Fbb$ for $\bullet=1,\infty$ and $i=0,\pm 1$, the maps 
$\tau_1(R)$ and $\tau_\infty(R)$ are forced and we only need to determined $\tau_0(R)$.
\\

The decompositions of $\Hbb_\infty(R)=\Hbb_1(R)=\Fbb^3$ according to 
relative $\SpinC$ classes gives 
$$\Hbb_1(R)=\langle a,b,c\rangle_\Fbb \ \ \ \text{and}\ \ \ 
\Hbb_\infty(R)=\langle a',b',c'\rangle_\Fbb, $$
where $a,a'$ are generators in relative $\SpinC$ class $-1$, $b,b'$ are generators in 
relative $\SpinC$ class $0$ and $c,c'$ are generators in relative $\SpinC$ class $+1$.
The homomorphisms
$\pphi_0(R)$ and $\pphibar_0(R)$ have the 
following block forms in the corresponding basis:
\begin{equation}\label{eq:phi(R)}
\pphi_0(R)=\left(\begin{array}{ccc}\alpha&0&0\\ 0&\beta&0\\ 0&0&\gamma
\end{array}\right)\ \ \text{and}\ \ 
\pphibar_0(R)=\left(\begin{array}{ccc}\gamma&0&0\\ 0&\beta&0\\ 0&0&\alpha
\end{array}\right).
\end{equation} 

From (\ref{eq:ranks}) we know that the ranks of $\pphi_0(R)$ and $\pphibar_0(R)$
are equal to $1$, i.e. precisely one of $\alpha,\beta$ and $\gamma$ is equal to $1$ and 
the other two are zero. Moreover, the rank of $\pphi_0(R)+\pphibar_0(R)$ is $2$, i.e.
precisely two of $\alpha+\gamma,\alpha+\gamma,2\beta$ are non-zero. 
Since the coefficient ring is 
$\Fbb$, $2b$ is automatically zero. Thus $\alpha=1, \beta=\gamma=0$ or 
$\gamma=1, \alpha=\beta=0$. \\

The generator $a$ of $\Hbb_1(R)$
is not in the image of $\pphi_\infty(R)$, since $\pphi_\infty(R,-1)$ is trivial. 
Hence $a$ is not in the kernel of $\pphi_0(R,-1)$.
Thus, from the above two possibilities the former is the case, i.e.
in (\ref{eq:phi(R)}) we get $\alpha=1$ and $\beta=\gamma=0$.\\ 

The rank of $\pphibar_\infty(R)$ is equal to $2$ according to (\ref{eq:ranks}).
Moreover, $\langle a,b\rangle_\Fbb$ is already in the image of 
$\pphibar_0(R)$. Thus,
$\pphibar_0(R)$ is surjective  onto $\Hbb_1(R,-1)\oplus \Hbb_1(R,0)$. Let us use a basis 
$a'', b''$ for $\Hbb_0(R,-\half)$ which contains some pre-image $a''$ of $a$ under 
 $\pphibar_\infty$ and an element $b''$ in the kernel of  $\pphibar_\infty$.
Use the dual basis $\tau_0(b''),\tau_0(a'')$ for $\Hbb_0(R,\half)$. 
The basis $\{a'',b'',\tau_0(b''),\tau_0(a'')\}$
for $\Hbb_0(R)$ is thus invariant under $\tau_0=\tau_0(R)$.
Correspondingly, we get 
\begin{equation}\label{eq:psi(R)}
\pphibar_\infty(R)=\left(\begin{array}{cccc}1&0&0&0\\ 0&0&x&y\\ 0&0&0&0
\end{array}\right)\ \ \text{and}\ \ 
\pphi_\infty(R)=\left(\begin{array}{cccc}0&0&0&0\\ y&x&0&0\\ 0&0&0&1
\end{array}\right).
\end{equation} 

If $x=0$ then $y=1$, since the rank of $\pphi_\infty(R)$ is equal to $2$. 
The rank of $\pphi_\infty(R)+\pphibar_\infty(R)$ is then equal to $2$;
on the other hand, 
 (\ref{eq:ranks}) implies  that this rank is $3$, a contradiction. 
The contradiction implies that $x=1$. Replacing $a''$ with $a''-yb''$
 we obtain the presentation of $\pphi_\infty(R)$ and $\pphibar_\infty(R)$ in a 
new basis for $\Hbb_0(R)$ (which is still invariant under the involution 
$\tau_0(R)$) corresponding to the  values $x=1$ and $y=0$ in 
(\ref{eq:psi(R)}). From here, and by taking into account the fact that 
 the map $\theta(R)$ increases the $\SpinC$ grading by $\half$, while 
 $\thetabar(R)$ decreases the $\SpinC$ grading by $\half$,
\begin{equation}\label{eq:thetabar(R)}
\theta(R)=\left(\begin{array}{cccc}
0&0&0&0\\
1&0&0&0\\
0&0&1&0
\end{array}\right)\ \ \ \text{and}\ \ \ 
\thetabar(R)=\left(\begin{array}{cccc}
0&1&0&0\\
0&0&0&1\\
0&0&0&0
\end{array}\right). 
\end{equation}
The above computations imply that $a_0(R)=1$ while $a_1(R)=a_\infty(R)=2$.
Moreover, we may take
\begin{equation}\label{eq:tau(R)}
\begin{split}
&A_0(R)=D_0(R)=\left(\begin{array}{cc}0&0\\ 0&0\end{array}\right),\ \ \ \ 
A_1(R)=D_\infty(R)=\left(0\right),\\
&C_1(R)=B_\infty(R)=B_1^T(R)=C_\infty^T(R)=\left(\begin{array}{c}1\\ 0\end{array}\right),\\
&B_0(R)=C_0(R)=\left(\begin{array}{cc}0&1\\ 1&0\end{array}\right)\ \ \ \text{and}\ \ \ 
D_1(R)=A_\infty(R)=\left(\begin{array}{cc} 0&0\\ 0&1\end{array}\right).
\end{split}
\end{equation}

For the left handed trefoil a similar argument may be used for the computation, which 
is sketched below. The rank of $\pphi_0(L)$ is $2$ and the rank of $\pphi_\infty(L)$
is $3$. The latter implies that the rank of $\pphi_\infty(L,1)$ is one, the rank of 
$\pphi_\infty(L,0)$ is $2$ and the rank of $\pphi_\infty(L,-1)$ is zero. 
Correspondingly, the ranks of 
$\pphi_0(L,1), \pphi_0(L,0)$ and $\pphi_0(L,-1)$ are equal to $0,1$ and $1$ 
respectively. If the images of $\pphi_\infty(L,0)$ and $\pphibar_\infty(L,0)$ 
are identical the maps $\pphi_0(L,0)$ and $\pphibar_0(L,0)$ are forced to be identical, 
since $\Hbb_\infty(L,0)$ is 
$1$-dimensional. In particular, $\pphi_0(L,0)+\pphibar_0(L,0)$ is trivial. Hence 
the rank of $\pphi_0(L)+\pphibar_0(L)$ is at most $2$, which is in contradiction with 
$y_0(L)=2$. The $2$-dimensional subspaces $\Image(\pphi_\infty(L,0))$ and 
$\Image(\pphibar_\infty(L,0))$ of $\Hbb_1(L,0)$ are thus different. 
From here, their intersection
is $1$-dimensional and is generated by some $\tau_1(L)$-invariant element
$\pphi_\infty(b)$ with $b\in \Hbb_0(L,-\half)$. \\

Let $a\in\Hbb_0(L,\half)$ denote the 
unique non-trivial vector in the kernel of $\pphi_\infty(L)$. 
Let us first assume that $b=\tau_0(a)$.
Complete $a$ to a basis $(a,c)$ for $\Hbb_0(L,\half)$. Then
\begin{displaymath}
\big\{a,c,\tau_0(a),\tau_0(c)\big\}
\end{displaymath} 
is an ordered basis for $\Hbb_0(L)$. Correspondingly, we obtain the basis
\begin{displaymath}
\big\{\pphi_\infty(c),\pphibar_\infty(c),\pphi_\infty(\tau_0(a)), \pphi_\infty(\tau_0(c)),
\tau_1\pphi_\infty(c))\big\}
\end{displaymath}
for $\Hbb_1(L)$, and the matrices $\pphi_\infty(L)$ and $\pphibar_\infty(L)$ 
take the following forms, respectively:
\begin{displaymath}
\pphi_\infty(L)=\colvec{0&1&0&0\\ 0&0&0&0\\ 0&0&1&0\\ 0&0&0&1\\ 0&0&0&0}\ \ \ 
\text{and}\ \ \ 
\pphibar_\infty(L)=\colvec{0&0&0&0\\ 0&1&0&0\\ 1&0&0&0\\ 0&0&0&0\\ 0&0&0&1}.
\end{displaymath}
In particular, the matrix
$$\pphi_\infty(L)+\pphibar_\infty(L)=
\colvec{0&1&0&0\\ 0&1&0&0\\ 1&0&1&0\\ 0&0&0&1\\ 0&0&0&1}$$
is a matrix of rank $3$, while we should have 
$$\rank(\pphi_\infty(L)+\pphibar_\infty(L))=\frac{h_0(L)+h_1(L)-y_\infty(L)}{2}=4.$$ 
This contradiction implies that $b$ is different from $\tau_0(a)$, and we may take
$(\tau_0(a),b)$ as a basis for $\Hbb_0(L,\half)$. Correspondingly we obtain the 
basis 
\begin{displaymath}
\big\{\tau_0(a),b,\tau_0(b),a\big\}
\end{displaymath}
for $\Hbb_0(L)$. As a basis for $\Hbb_1(L,0)$ we obtain the three vectors
$\pphi_\infty(a),\pphibar_\infty(a)$ and $\pphibar_\infty(b)$. 
This basis is completed to the following (ordered) basis for $\Hbb_1(L)$:
\begin{displaymath}
\big\{\pphi_\infty(a),\pphibar_\infty(b),\pphi_\infty(b),\pphibar_\infty(a),
\tau_1(\pphi_\infty(b))\big\}.
\end{displaymath}
Finally, we choose the following basis for $\Hbb_\infty(L)$:
\begin{displaymath}
\big\{\pphibar_0(\pphi_\infty(b)), \pphi_0(\pphibar_\infty(a)),
\tau_\infty(\pphibar_0(\pphi_\infty(b)))\big\}.
\end{displaymath}
In these bases we may compute
\begin{displaymath}
\pphi_\infty(L)=\colvec{1&0&0&0\\ 0&1&0&0\\ 0&0&1&0\\ 0&0&0&0}\ \ \ \text{and}
\ \ \ \pphi_0(L)=\colvec{0&0&0&0&0\\ 0&0&0&1&0\\ 0&0&0&0&1}.
\end{displaymath}
Moreover, after re-ordering the elements of the above bases, we find the following 
presentations:
\begin{displaymath}
\begin{split}
&D_0(L)=A_\infty(L)=0,\ \ A_1(L)=\colvec{0&0\\ 0&0}, \ \ 
 B_1(L)=C_1^T(L)=\colvec{1&0&0\\ 0&0&1},\\   
&A_0(L)=\colvec{0&0&0\\ 0&0&1\\ 0&1&0},\ \ \  \ \ D_1(L)=\colvec{0&0&0\\ 0&1&0\\
0&0&0},\ \  \ \ \
D_\infty(L)=\colvec{1&0\\ 0&0},\\
&B_0(L)=C_0^T(L)=\colvec{1\\ 0\\ 0}\ \ \ 
\ \ \ \text{and} \ \ \ \ \ \ \ B_\infty(L)=C_\infty^T(L)=\colvec{0&1}.
\end{split}
\end{displaymath}

\subsection{Splicing a knot complement with the complement of 
a trefoil}\label{subsec:splicing-with-T}

For a knot $K\subset Y$, let $Y(R,K)$  denote the three-manifold
obtained by splicing the complement of ${K}\subset {Y}$ 
with the complements of the 
right handed  trefoil. We  study 
the rank $r_r(K)$  of $\ov\HFT(Y(R,K))$  in this subsection.
With the notation of Subsection~\ref{subsec:simplifications}, $r_r(K)=i(\Dd'(R,K))$.
 Replacing the block forms of (\ref{eq:tau(R)}) 
in $\Dd'(R,K)$, we find 
\begin{displaymath}
\Dd'(R,K)=\colvec{
0&0&0&0&0&I&0&B_1B_0&0&0\\
B_\infty B_1 &0&0&0&0&0&0&0&0&I \\
0&B_\infty B_1 &0&0&B_\infty A_1 &0&0&I&0&0 \\
D_\infty B_1 &0&I&0&0&0&0&0&0&0\\
0&D_\infty B_1& 0& I& D_\infty A_1&0&0&0&0&0\\
I&0&0&0&B_0B_\infty &0&0&0&0& B_0X\\
0&0&0&0&0&B_0B_\infty &0&0&0& 0\\
0&0&0&0&0&0& I& D_1B_0&0&0\\
0&0&0&0& D_0B_\infty &0&0&0&I&D_0X\\
0&0&0&0&0&D_0B_\infty &0&0& 0& I\\
},
\end{displaymath}
where $A_\bullet=A_\bullet(K), B_\bullet=B_\bullet(K), C_\bullet=C_\bullet(K),
D_\bullet=D_\bullet(K)$ and $X=X(K)$ for $\bullet\in\{0,1,\infty\}$.
Doing a series of cancellations that correspond to the identity matrices which appear 
as the 
$$(1,6), (3,8), (4,3), (5,4), (6,1), (8,7), (9,9) \ \text{and}\ (10,10)$$
entries in the above block presentation
we obtain the equivalent matrix
\begin{equation}\label{eq:rank-d2}
\begin{split}
R_r(K):&=
\left(\begin{array}{ccc}
0&B_0XB_\infty\\
XB_\infty B_1&X B_\infty A_1+D_0XB_\infty\\
\end{array}\right)
\end{split}
\end{equation}

\begin{cor} \label{cor:splicing-with-R}
For a knot $K\subset Y$ let $Y(R,K)$ denote the three-manifold obtained 
by splicing the complement of $K$ and the complement of the trefoil. Then
\begin{equation}\label{eq:rank-h1}
\begin{split}
\ov{\HFT}(Y(R,K)) &=\imath\big(R_r(K)\big).
\end{split}
\end{equation} 
\end{cor}
\begin{proof}
The claim follows immediately from the above discussion.
\end{proof}
For the  trefoils, our computations imply that 
\begin{displaymath}
\begin{split}
&X(R)B_\infty(R)=X(L)B_\infty(L)=0\\ 
\Rightarrow\ \ \ &R_r(R)=R_r(L)=0 \\
\Rightarrow\ \ \ &\left|\ov{\HFT}(Y(R,R))\right|=7\ \ \ \text{and}\ \ \ 
\left|\ov{\HFT}(Y(R,L))\right|=9.
\end{split}
\end{displaymath}
The above computations agree with the computations of Hedden and Levine 
\cite{Matt-Adam}.
\begin{cor}\label{cor:splicing-with-trefoil}
For every knot $K$ in a homology sphere $Y$ we have
\begin{displaymath}
\begin{split}
\left|\ov{\HFT}(Y(R,K))\right|&\geq
 \big(a_0(K)+a_1(K)+2a_\infty(K)\big)-4\min\big\{a_0(K),a_1(K),a_\infty(K)\big\}\\
 &=4\max\big\{h_0(K),h_1(K),h_\infty(K)\big\}-\big(h_0(K)+h_1(K)+2h_\infty(K)\big).
\end{split}
\end{displaymath}
Moreover, if $Y(R,K)$ is a homology sphere $L$-space $K$ is trivial and $Y$ is a 
homology sphere $L$-space.
\end{cor}
\begin{proof}
Let $M=M(K)=X(K)B_\infty(K)$ and note that 
\begin{displaymath}
\begin{split}
\rank(R_r(K))&=\rank\left(\begin{array}{cc}
0& B_0(K)M\\ MB_1(K)& MA_1(K)+D_0(K)M
\end{array}\right) \\
&\leq \rank\left(\begin{array}{cc}
 MB_1(K)& MA_1(K)
\end{array}\right) +\rank\left(\begin{array}{c}
B_0(K)M\\ D_0(K)M
\end{array}\right)\\
&=2\rank(M)\leq 2\rank(X(K)) 
\end{split}
\end{displaymath}
For every knot $K\subset Y$ as above note that the rank of $X=X(K)$ is at most 
equal to the minimum of the sizes of the matrices $B_0(K), B_1(K)$ and $B_\infty(K)$,
which is 
$$\min\big\{a_0(K),a_1(K),a_\infty(K)\big\}.$$
Since $R_r(K)$ is of size $h_0(K)\times h_1(K)=(a_1(K)+a_\infty(K))\times (a_0(K)+a_\infty(K)$ 
this proves the first part of the corollary.\\

Let us assume that $\rank\left(\ov{\HFT}(Y(R,K))\right)=1$. From here we find
\begin{displaymath}
\begin{split}
&\big(a_0(K)+a_1(K)+2a_\infty(K)\big)-4\min\big\{a_0(K),a_1(K),a_\infty(K)\big\}\\
&\ \ \ \ \ \ \ \ \ \ \ \ \ \ \ \ \ \ \ \ \ \ \ \
=\big(a_0(K)+a_1(K)+2a_\infty(K)\big)-4\rank(M)=1.
\end{split}
\end{displaymath}
Since $a_1(K)$ and $a_\infty(K)$ have the same parity while the parity of $a_0(K)$ is different 
from the parity of both $a_1(K)$ and $a_\infty(K)$, one can easily conclude that
 $a_0(K)-1=a_1(K)=a_\infty(K)$. Let $a$ denote the common value $a_1(K)=a_\infty(K)$. 
 Then the rank of $M$ is $a$
and both $B_0(K)$ and $X(K)$ are invertible. We may thus assume that $A_0(K)=D_0(K)=0$.
Since
\begin{displaymath}
\rank(\pphi_\infty(K)+\pphibar_\infty(K))=\rank \colvec{
B_1(K)A_0(K)& B_1(K)B_0(K)\\ I+D_1(K)A_0(K)& D_1(K)B_0(K)}=2a
\end{displaymath}
the three-manifold $Y$ is an $L$-space. Since splicing $K$ with the trefoil
is also a homology sphere $L$-space we conclude that $K$ is trivial, by 
Theorem 1 from  \cite{Matt-Adam}.
\end{proof}
\newpage
\section{Appendix; Bordered Floer homology for knot complements}
The first draft of this paper appeared while the theory of bordered Floer homology 
was being developed. With bordered Floer homology conventions widely known to the 
Heegaard Floer community, the referee recommended the inclusion of an appendix 
which addresses the contribution of this paper within the realm of bordered Floer homology.\\

Let $K\subset Y$ denote a null-homologous knot inside the three-manifold $Y$, and 
let $H=(\Sig,\alphas,\ov\betas \cup\{\lambda,\mu\};z)$ denote a special Heegaard 
diagram for $K$, as constructed in Lemma~\ref{lem:nice-to-nicer}. 
In particular, $H$ is a nice Heegaard diagram  for the bordered 
three-manifold $Y_K$ determined by $K\subset Y$ in the sense of \cite{LOT}. 
The Bordered Floer 
complex $\ov\CFDT(Y_K)$ may then be constructed from the chain complexes 
$M=M(K)$ and $L=L(K)$
(which are described in Proposition~\ref{prop:ML-mapping-cone} as the mapping 
cones of $\pphibar_\infty(K)\co C_0(K)\ra C_1(K)$ and 
$\pphi_0(K)\co C_1(K)\ra C_\infty(K)$, respectively) and the chain maps 
$\Phi=\Phi(K)\co L\ra M$ and $\Psi_{i}=\Psi_i(K)\co M\ra L, i=1,2,3$.\\

More precisely and following the notation of Subsection~4.2 from 
\cite{LOT-E}, the idempotents
$\imath_0$ and $\imath_1$, and the chords $\rho_1,\rho_2,\rho_3,\rho_{12}=\rho_1\rho_2,
\rho_{23}=\rho_2\rho_3$ and $\rho_{123}=\rho_1\rho_2\rho_3$ 
form a $\Fbb$-basis for the differential graded 
algebra associated with the torus boundary;
\begin{displaymath}
\begin{split}&\\
&\mathcal{A}(T^2,0)=
\Big\langle\begin{diagram}[w=3em]
\imath_0 \ \bullet\ & \upperarrow{\rho_1}
\lift{-2}
{\strarrow{\rho_2}}
\lowerarrow{\rho_3} &\  \bullet\ \imath_1
\end{diagram}\Big\rangle/\left(\rho_2\rho_1=\rho_3\rho_2=0\right).\\
&
\end{split}
\end{displaymath}  
The module $\ov\CFDT(Y_K)$ is generated (over $\mathcal{A}(T^2,0)$) by 
the generators of $M$ and $L$. For  a generator $\x$ of $L$ we have
\begin{equation}\label{eq:d_L}
\begin{split}
& I(\x)=\imath_0\ \ \ \ \text{and}\ \ \ \
\partial (\x)=d_L(\x)+\rho_1\Psi_1(\x)+\rho_3\Psi_2(\x)+\rho_{123}\Psi_3(\x),
\end{split}
\end{equation}
while for a generator $\y$ of $M$ we have 
\begin{equation}\label{eq:d_M}
\begin{split}
& I(\y)=\imath_1\ \ \ \ \ \text{and}\ \ \ \ \ 
\partial (\y)=d_M(\y)+\rho_2\Phi(\y).
\end{split}
\end{equation}

The splicing formula of (\ref{eq:complex-cube}) is then just the gluing formula for bordered Floer homology, i.e. Theorem 1.3 from \cite{LOT}. A related discussion is 
carried over in Section 8 of \cite{LOT}.\\

\begin{defn}\label{def:admissible-data}
The chain complexes $(C_\bullet(K),d_\bullet), \bullet\in\{0,1,\infty\}$ and 
the chain maps $f_\bullet(K),\fbar_\bullet(K), \bullet\in\{0,\infty\}$   
are called {\emph{admissible data}} associated with the knot $K$ if they satisfy the following 
conditions:
\begin{itemize}
\item The homology of the complex $(C_\bullet(K),d_\bullet)$ is $\Hbb_\bullet(K)$.
\item The maps induced by $f_\bullet(K)$ and $\fbar_\bullet(K)$ in homology (under the 
identification of  the homology of  $(C_\bullet(K),d_\bullet)$ with $\Hbb_\bullet(K)$)
are $\pphi_\bullet(K)$ and $\pphibar_\bullet(K)$, respectively.
\item We have $f_0(K)\circ f_\infty(K)=\fbar_0(K)\circ \fbar_\infty(K)=0$.
\item The corresponding maps 
$$\theta(K)\co \Ker(\pphi_\infty(K))\ra \Coker(\pphi_0(K))\ \ \text{and}\ \ 
\thetabar(K)\co \Ker(\pphibar_\infty(K))\ra \Coker(\pphibar_0(K))$$
are isomorphisms, and are the inverses of the maps induces by $\pphi_1(K)$ and 
$\pphibar_1(K)$, respectively.
\end{itemize}
\end{defn}

 The proof of Theorem~\ref{thm:cube}
implies that $(C_\bullet^i,d_\bullet^i)$ and 
the chain maps $f_\bullet^i,\fbar_\bullet^i$ for  $\bullet\in\{0,\infty\},i=1,2$
in (\ref{eq:complex-cube}) may be replaced by other admissible data corresponding 
to the knots $K_1$ and $K_2$. Correspondingly, the bordered Floer complex associated 
with any knot $K\subset Y$ may be constructed from admissible data associated with $K$.
More precisely, we have the following proposition.

\begin{prop}\label{prop:BFH}
Suppose that the chain complexes $(C_\bullet(K),d_\bullet), \bullet\in\{0,1,\infty\}$ and  
the chain maps $f_\bullet=f_\bullet(K),\fbar_\bullet=\fbar_\bullet(K), \bullet\in\{0,\infty\}$ 
are admissible data associated with the knot $K\subset Y$ and set
 \begin{displaymath}
 \begin{split}
 &M(K)=C_0(K)\oplus C_1(K),\ \ L(K)=C_1(K)\oplus C_\infty(K)\\
 \end{split}
 \end{displaymath}
 The bordered Floer complex $\ov\CFDT(Y_K)$ may then be constructed as the 
 left module over the differential graded algebra $\mathcal{A}(T^2,0)$ which is generated 
 by $\imath_0.L(K)$ and $\imath_1. M(K)$, and equipped with the differential 
 $\partial\co \ov\CFDT(Y_K)\ra \ov\CFDT(Y_K)$ defined by 
 \begin{equation}\label{eq:differential}
\partial \colvec{
\x\\ \y }=
\begin{cases}
\colvec{
d_0(\x)\\ \fbar_\infty(\x)+d_1(\y)}
+\rho_2.\colvec{
0\\ \x}\ \ \ \ &\text{if }\colvec{\x\\ \y}\in M(K)\\
\colvec{
d_1(\x)\\ f_0(\x)+d_\infty(\y)}
+\colvec{
\rho_1f_\infty(\x)\\ \rho_3\fbar_0(\y)+\rho_1\rho_2\rho_3\fbar_0(f_\infty(\x))}
\ \ \ \ &\text{if }\colvec{\x\\ \y }\in L(K)\\
\end{cases} 
 \end{equation}
 \end{prop}

In particular, let the $\Fbb$-modules $\A_\bullet=\A_\bullet(K),\bullet\in\{0,1,\infty\}$ 
and the matrices 
$A_\bullet=A_\bullet(K),B_\bullet=B_\bullet(K),C_\bullet=C_\bullet(K)$ and 
$D_\bullet=D_\bullet(K)$ be defined as in 
Subsection~\ref{subsec:properties}. Set 
\begin{displaymath}
\begin{split}
&(C_0(K),d_0)=\left(\A_\infty\oplus \A_1,0\right),\ \ \ \ 
(C_\infty(K),d_\infty)=\left(\A_1\oplus \A_0,0\right)\\
&C_1(K)=\A_1\oplus \A_0\oplus \A_\infty\oplus \A_1\ \ \ \text{and}\ \ \ 
d_1=\left(\begin{array}{cccc}
0&0&0&0\\ 0&0&0&0\\ 0&0&0&0\\ I_{\A_1}&0&0&0\\ 
\end{array}\right).
\end{split}
\end{displaymath}
Correspondingly, define
\begin{displaymath}
\begin{split}
&f_\infty(K)=\left(\begin{array}{cc}
0&0\\ 0&0\\ I&0\\ 0&I
\end{array}\right),\ 
f_0(K)=\left(\begin{array}{cccc}
I&0&0&0\\ 0&I&0&0
\end{array}\right)\ \text{and}\ 
\tau_1(K)=\left(\begin{array}{cccc}
0&0&0&0\\ 0&A_1&B_1&0\\ 0&C_1&D_1&0\\ 0&0&0&0
\end{array}\right),
\end{split}
\end{displaymath}
and set $\fbar_\infty(K)=\tau_1(K) f_\infty(K)\tau_0(K)$ and 
$\fbar_0(K)=\tau_\infty(K) f_0(K)\tau_1(K)$. The data associated with $K$ 
consisting of $(C_\bullet(K),d_\bullet)$ and 
 $f_\bullet(K),\fbar_\bullet(K), \bullet\in\{0,\infty\}$ is then admissible. \\
 
 Corresponding to the above admissible data and associated with $K\subset Y$ 
  we may construct
 the bordered Floer complex for $K$ via
 \begin{displaymath}
 \begin{split}
 &M(K)=C_0(K)\oplus C_1(K)=
 \A_\infty\oplus \A_1\oplus \A_1\oplus \A_0\oplus \A_\infty\oplus \A_1\\
&L(K)=C_1(K)\oplus C_\infty(K)=
 \A_1\oplus \A_0\oplus \A_\infty\oplus \A_1\oplus \A_1\oplus \A_0\\ 
 &d_M=\left(\begin{array}{cccccc}
0&0&0&0&0&0\\
 0&0&0&0&0&0\\
 0&0&0&0&0&0\\
 B_1A_0& B_1B_0&0&0&0&0\\
D_1A_0& D_1B_0&I&0&0&0\\
 0&0&0&0&0&0 
 \end{array}\right),\ \ \ \ \ \ \ 
 d_L=\left(\begin{array}{cccccc}
 0&0&0&0&0&0\\
 0&0&0&0&0&0\\
 0&0&0&0&0&0\\
 I&0&0&0&0& 0\\
 0&0&I&0&0& 0\\
 0&0&0&I&0&0
 \end{array}\right)\\ 
 \end{split}
 \end{displaymath}
  \begin{displaymath}
 \begin{split}
 &\Phi(K)=\left(\begin{array}{cccccc}
 0&0&0&0&0&0\\
 0&0&0&0&0&0\\
 I&0&0&0&0&0\\
 0&I&0&0&0& 0\\
 0&0&I&0&0& 0\\
 0&0&0&I&0&0
 \end{array}\right),\ \ \ \ \ \ \ \ \ \ \ \Psi_1(K)=\left(\begin{array}{cccccc}
 0&0&0&0&0&0\\
 0&0&0&0&0&0\\
 I&0&0&0&0&0\\
 0&I&0&0&0& 0\\
 0&0&0&0&0& 0\\
 0&0&0&0&0&0
 \end{array}\right)\\
 &\Psi_2(K)=\left(\begin{array}{cccccc}
 0&0&0&0&0&0\\
 0&0&0&0&0&0\\
 0&0&0&0&0&0\\
 0&0&0&0&0& 0\\
 0&0&0&B_\infty A_1&B_\infty B_1& 0\\
 0&0&0&D_\infty A_1&D_\infty B_1&0
 \end{array}\right)\ \ \text{and}\ \ \Psi_3(K)=\Psi_2(K)\Phi(K)\Psi_1(K).
 \end{split}
 \end{displaymath}
 as the left module over the differential graded algebra $\mathcal{A}(T^2,0)$
 generated by $\imath_0.L$ and $\imath_1.M$
 and equipped with the differential 
 $\partial\co \ov\CFDT(Y_K)\ra \ov\CFDT(Y_K)$ defined by 
 the equations~(\ref{eq:d_L}) and (\ref{eq:d_M}).\\
 
 \begin{remark}
 Simultaneous computation of the matrices $\tau_\bullet(K)=\left(\begin{array}{cc}
A_\bullet&B_\bullet\\ C_\bullet& D_\bullet  \end{array}\right)$ is {\emph{a priori}} 
quite difficult,
as we observed in the case of trefoils in Section~\ref{sec:example}. This 
makes the above description of the bordered Floer homology hard to use even for knots 
$K\subset Y$ where we have complete understanding of the Heegaard Floer complex 
associated with $K$. However, it is possible to construct admissible data associated 
with $K\subset Y$  completely
in terms of the filtered chain complex 
$\CFT^\infty(Y,K;\Fbb)$, as will be discussed in the revision of 
\cite{Ef-essential}.  
 \end{remark}
\newpage


\end{document}

%% file: command1.tex
\newcommand{\colvec}[2][.8]{%
  \scalebox{#1}{%
    \renewcommand{\arraystretch}{.8}%
    $\begin{pmatrix}#2\end{pmatrix}$%
  }
}
\newcommand{\lift}[2]{%
\setlength{\unitlength}{1pt}
\begin{picture}(0,0)(0,0)
\put(0,{#1}){\makebox(0,0)[b]{${#2}$}}
\end{picture}
}
\newcommand{\lowerarrow}[1]{%
\setlength{\unitlength}{0.03\DiagramCellWidth}
\begin{picture}(0,0)(0,0)
\qbezier(-28,-4)(0,-18)(28,-4)
\put(0,-14){\makebox(0,0)[t]{$\scriptstyle {#1}$}}
\put(28.6,-3.7){\vector(2,1){0}}
\end{picture}
}
\newcommand{\upperarrow}[1]{%
\setlength{\unitlength}{0.03\DiagramCellWidth}
\begin{picture}(0,0)(0,0)
\qbezier(-28,12)(0,25)(28,12)
\put(0,20){\makebox(0,0)[b]{$\scriptstyle {#1}$}}
\put(28.6,11.7){\vector(2,-1){0}}
\end{picture}
}
\newcommand{\strarrow}[1]{%
\setlength{\unitlength}{0.03\DiagramCellWidth}
\begin{picture}(0,0)(0,0)
\qbezier(-22,5)(0,5)(22,5)
\put(0,6){\makebox(0,0)[b]{$\scriptstyle {#1}$}}
\put(-22,4.8){\vector(-2,0){0}}
\end{picture}
}

\newcommand{\tteettaa}{{\pphi_0^1\pphibar_\infty^1\otimes\pphi_0^2\pphibar_\infty^2
+\theta^1\otimes \thetabar^2+\thetabar^1\otimes\theta^2}}

\newcommand{\co}{\colon\thinspace}
\newcommand{\fbar}{\overline{f}}
\newcommand{\f}{{f}}
\newcommand{\pphi}{\mathfrak{f}}
\newcommand{\pphibar}{\overline{\pphi}}
\newcommand{\CFDT}{\mathrm{CFD}}
\newcommand{\ovl}{\overline}
\newcommand{\ubar}{{\ovl{u}}}
\newcommand{\vbar}{{\ovl{v}}}
\newcommand{\wbar}{{\ovl{w}}}
\newcommand{\M}{{\text{\mancube}}}
\newcommand\relspincs{\underline{\mathfrak s}}

\newcommand{\Thet}{P}

\newcommand{\half}{\frac{1}{2}}

\newcommand{\Dd}{\mathfrak{D}}
\newcommand{\sig}{\sigma}

\newcommand{\thetabar}{\overline{\theta}}
\newcommand{\bb}{\mathfrak{b}}

\newcommand{\nd}{\mathrm{nd}}
\newcommand{\bfrak}{\mathfrak{b}}
\newcommand{\cfrak}{\mathfrak{c}}

\newcommand{\B}{\mathbb{B}}

\newcommand{\HFT}{\mathrm{HF}}
\newcommand{\CFT}{\mathrm{CF}}
\newcommand{\HFKT}{\mathrm{HFK}}

\newcommand{\rank}{\mathrm{rnk}}

\newcommand{\Dcal}{\mathcal{D}}

\newcommand{\ov}{\widehat}

\newcommand{\E}{\mathbb{E}}
\newcommand{\Sig}{\Sigma}
\newcommand{\ra}{\rightarrow}
\newcommand{\xra}{\xrightarrow}
\newcommand{\A}{\mathbb{A}}

\newcommand{\Hbbc}{\text{\mancube}}

\newcommand{\phibar}{{\overline{\phi}}}

\newcommand\Dual{\mathcal D}
\newcommand\Duality\Dual

\newcommand\RelSpinC{\underline{\SpinC}}
\newcommand\relspinc{s}

\newcommand\x{\mathbf x}
\newcommand{\cc}{\mathfrak{c} }

\newcommand\y{\mathbf y}

\newcommand\ModSphere{\ModFlow\left({\mathbb S}\longrightarrow
\Sym^{g-1}(\Sigma_{1})\times \Sym^2(\Sigma_{2})\right)}
\newcommand\ModSpheres\ModSphere

\newcommand\UnparModSp{\widehat \ModSp}
\newcommand\UnparModFlow\UnparModSp

\newcommand\ModMaps{\mathcal M}
\newcommand\ModSp\ModMaps

\newcommand\Ta{{\mathbb T}_{\alpha}}

\newcommand\alphas{\mbox{\boldmath$\alpha$}}

\newcommand\betas{\mbox{\boldmath$\beta$}}
\newcommand\gammas{\mbox{\boldmath$\gamma$}}

\newcommand\spincrel\relspinc

\newtheorem{thm}{Theorem}[section]
\newtheorem{prop}[thm]{Proposition}
\newtheorem{cor}[thm]{Corollary}

\newtheorem{lem}[thm]{Lemma}

\newtheorem{defn}[thm]{Definition}

\newtheorem{remark}[thm]{Remark}

\def\endproof{\relax\ifmmode\expandafter\endproofmath\else
  \unskip\nobreak\hfil\penalty50\hskip.75em\hbox{}\nobreak\hfil\bull
  {\parfillskip=0pt \finalhyphendemerits=0 \bigbreak}\fi}
\def\endproofmath$${\eqno\bull$$\bigbreak}
\def\bull{\vbox{\hrule\hbox{\vrule\kern3pt\vbox{\kern6pt}\kern3pt\vrule}\hrule}}

\newcommand{\Q}{\mathbb{Q}}

\newcommand{\C}{\mathbb{C}}

\newcommand{\Z}{\mathbb{Z}}

\newcommand{\Ker}{\mathrm{Ker}}

\newcommand{\Coker}{\mathrm{Coker}}

\newcommand{\Image}{\mathrm{Im}}

\newcommand{\ModSWfour}{\mathcal{M}}
\newcommand{\ModFlow}{\ModSWfour}

\newcommand{\SpinC}{{\mathrm{Spin}}^c}

\newcommand\abuts\Rightarrow
\newcommand\Sym{\mathrm{Sym}}

\newcommand{\Cinfinf}{{C^1_\infty\otimes C^2_\infty}}
\newcommand{\Cinfone}{{C^1_\infty\otimes C^2_1}}
\newcommand{\Coneinf}{{C^1_1\otimes C^2_\infty}}
\newcommand{\Coneone}{{C^1_1\otimes C^2_1}}
\newcommand{\Conezero}{{C^1_1\otimes C^2_0}}
\newcommand{\Czeroone}{{C^1_0\otimes C^2_1}}
\newcommand{\Czerozero}{{C^1_0\otimes C^2_0}}

\newcommand{\Hbb}{{\mathbb{H}}}
\newcommand{\Hinfinf}{{\Hbb_{\infty,\infty}}}
\newcommand{\Hinfone} {{\Hbb_{\infty,1}}}
\newcommand{\Honeinf} {{\Hbb_{1,\infty}}}
\newcommand{\Honeone} {{\Hbb_{1,1}}}
\newcommand{\Honezero} {{\Hbb_{1,0}}}
\newcommand{\Hzeroone} {{\Hbb_{0,1}}}
\newcommand{\Hzerozero} {{\Hbb_{0,0}}}

\newcommand{\Fbb}{\mathbb{F}}


%% file: standardcommand.tex
\newcommand{\lra}{\longrightarrow}



%% file: punctures-2.pdf_tex
\begingroup%
  \makeatletter%
  \providecommand\color[2][]{%
    \errmessage{(Inkscape) Color is used for the text in Inkscape, but the package 'color.sty' is not loaded}%
    \renewcommand\color[2][]{}%
  }%
  \providecommand\transparent[1]{%
    \errmessage{(Inkscape) Transparency is used (non-zero) for the text in Inkscape, but the package 'transparent.sty' is not loaded}%
    \renewcommand\transparent[1]{}%
  }%
  \providecommand\rotatebox[2]{#2}%
  \ifx\svgwidth\undefined%
    \setlength{\unitlength}{569.075bp}%
    \ifx\svgscale\undefined%
      \relax%
    \else%
      \setlength{\unitlength}{\unitlength * \real{\svgscale}}%
    \fi%
  \else%
    \setlength{\unitlength}{\svgwidth}%
  \fi%
  \global\let\svgwidth\undefined%
  \global\let\svgscale\undefined%
  \makeatother%
  \begin{picture}(1,0.42415323)%
    \put(0,0){\includegraphics[width=\unitlength]{punctures-2.pdf}}%
    \put(0.06396345,0.17644468){\color[rgb]{0,0,0}\makebox(0,0)[lb]{\smash{$\lambda_1$}}}%
    \put(0.16939769,0.06197326){\color[rgb]{0,0,0}\makebox(0,0)[lb]{\smash{$\lambda_0$}}}%
    \put(0.34210907,0.23870108){\color[rgb]{0,0,0}\makebox(0,0)[lb]{\smash{$\lambda_\infty$}}}%
    \put(0.14228603,0.29192033){\color[rgb]{0,0,0}\makebox(0,0)[lb]{\smash{$v$}}}%
    \put(0.1412819,0.14129992){\color[rgb]{0,0,0}\makebox(0,0)[lb]{\smash{$w$}}}%
    \put(0.25675751,0.28991203){\color[rgb]{0,0,0}\makebox(0,0)[lb]{\smash{$u$}}}%
    \put(0.63330844,0.23870113){\color[rgb]{0,0,0}\makebox(0,0)[lb]{\smash{$\lambda_1$}}}%
    \put(0.63130019,0.18246948){\color[rgb]{0,0,0}\makebox(0,0)[lb]{\smash{$\lambda_\infty$}}}%
    \put(0.72769724,0.12422958){\color[rgb]{0,0,0}\makebox(0,0)[lb]{\smash{$\ovl{u}$}}}%
    \put(0.8311232,0.26681693){\color[rgb]{0,0,0}\makebox(0,0)[lb]{\smash{$\ovl{w}$}}}%
    \put(0.79196186,0.06498562){\color[rgb]{0,0,0}\makebox(0,0)[lb]{\smash{$\lambda_0$}}}%
    \put(0.83614386,0.12523373){\color[rgb]{0,0,0}\makebox(0,0)[lb]{\smash{$\ovl{v}$}}}%
    \put(0.31198499,0.39835869){\color[rgb]{0,0,0}\makebox(0,0)[lb]{\smash{$\alpha$s}}}%
    \put(0.87731341,0.39635039){\color[rgb]{0,0,0}\makebox(0,0)[lb]{\smash{$\alpha$s}}}%
    \put(0.0991082,0.20456048){\color[rgb]{0,0,0}\makebox(0,0)[lb]{\smash{$\ubar$}}}%
    \put(0.24671616,0.20456048){\color[rgb]{0,0,0}\makebox(0,0)[lb]{\smash{$\vbar$}}}%
    \put(0.2179746,0.33531487){\color[rgb]{0,0,0}\makebox(0,0)[lb]{\smash{$\wbar$}}}%
    \put(0.7387427,0.20456048){\color[rgb]{0,0,0}\makebox(0,0)[lb]{\smash{$v$}}}%
    \put(0.75982955,0.09912622){\color[rgb]{0,0,0}\makebox(0,0)[lb]{\smash{$w$}}}%
    \put(0.89337961,0.20456048){\color[rgb]{0,0,0}\makebox(0,0)[lb]{\smash{$u$}}}%
  \end{picture}%
\endgroup%

%% file: Cube-1.pdf_tex
\begingroup%
  \makeatletter%
  \providecommand\color[2][]{%
    \errmessage{(Inkscape) Color is used for the text in Inkscape, but the package 'color.sty' is not loaded}%
    \renewcommand\color[2][]{}%
  }%
  \providecommand\transparent[1]{%
    \errmessage{(Inkscape) Transparency is used (non-zero) for the text in Inkscape, but the package 'transparent.sty' is not loaded}%
    \renewcommand\transparent[1]{}%
  }%
  \providecommand\rotatebox[2]{#2}%
  \ifx\svgwidth\undefined%
    \setlength{\unitlength}{397.23976416bp}%
    \ifx\svgscale\undefined%
      \relax%
    \else%
      \setlength{\unitlength}{\unitlength * \real{\svgscale}}%
    \fi%
  \else%
    \setlength{\unitlength}{\svgwidth}%
  \fi%
  \global\let\svgwidth\undefined%
  \global\let\svgscale\undefined%
  \makeatother%
  \begin{picture}(1,0.90899994)%
    \put(0,0){\includegraphics[width=\unitlength]{Cube-1.pdf}}%
    \put(0.00627769,0.86241677){\color[rgb]{0,0,0}\makebox(0,0)[lb]{\smash{$\Hinfinf$}}}%
    \put(0.62655799,0.86241677){\color[rgb]{0,0,0}\makebox(0,0)[lb]{\smash{$\Hinfone$}}}%
    \put(0.84808667,0.64088809){\color[rgb]{0,0,0}\makebox(0,0)[lb]{\smash{$\Honeone$}}}%
    \put(0.85010057,0.09713588){\color[rgb]{0,0,0}\makebox(0,0)[lb]{\smash{$\Honeone$}}}%
    \put(0.62857189,0.32067845){\color[rgb]{0,0,0}\makebox(0,0)[lb]{\smash{$\Honezero$}}}%
    \put(0.22579248,0.09914977){\color[rgb]{0,0,0}\makebox(0,0)[lb]{\smash{$\Hzeroone$}}}%
    \put(0.01433328,0.32067845){\color[rgb]{0,0,0}\makebox(0,0)[lb]{\smash{$\Hzerozero$}}}%
    \put(0.22780637,0.64088809){\color[rgb]{0,0,0}\makebox(0,0)[lb]{\smash{$\Honeinf$}}}%
    \put(0.48155741,0.12533043){\color[rgb]{0,0,0}\makebox(0,0)[lb]{\smash{$\pphibar_\infty^1\otimes I$}}}%
    \put(0.38086255,0.2864422){\color[rgb]{0,0,0}\makebox(0,0)[lb]{\smash{$\pphibar_\infty^1\otimes I$}}}%
    \put(0.07877801,0.24616428){\color[rgb]{0,0,0}\rotatebox{-45}{\makebox(0,0)[lb]{\smash{$I\otimes \pphibar_\infty^2$}}}}%
    \put(0.74336402,0.26630322){\color[rgb]{0,0,0}\rotatebox{-45}{\makebox(0,0)[lb]{\smash{$I\otimes \pphibar_\infty^2$}}}}%
    \put(0.01836108,0.40727603){\color[rgb]{0,0,0}\rotatebox{90}{\makebox(0,0)[lb]{\smash{$\tteettaa$}}}}%
    \put(0.53190483,0.65901316){\color[rgb]{0,0,0}\makebox(0,0)[lb]{\smash{$I\otimes \pphi_0^2$}}}%
    \put(0.26002873,0.89061132){\color[rgb]{0,0,0}\makebox(0,0)[lb]{\smash{$I\otimes \pphi_0^2$}}}%
    \put(0.1391949,0.78991646){\color[rgb]{0,0,0}\rotatebox{-45}{\makebox(0,0)[lb]{\smash{$\pphi_0^1\otimes I$}}}}%
    \put(0.72322501,0.81005548){\color[rgb]{0,0,0}\rotatebox{-45}{\makebox(0,0)[lb]{\smash{$\pphi_0^1\otimes I$}}}}%
    \put(0.29023718,0.36699808){\color[rgb]{0,0,0}\rotatebox{90}{\makebox(0,0)[lb]{\smash{$\pphi_\infty^1\otimes\pphibar_0^2$}}}}%
    \put(0.68294711,0.44755397){\color[rgb]{0,0,0}\rotatebox{90}{\makebox(0,0)[lb]{\smash{$\pphibar_0^1\otimes \pphi_\infty^2$}}}}%
    \put(0.87279625,0.36078726){\color[rgb]{0,0,0}\makebox(0,0)[lb]{\smash{$I$}}}%
    \put(0.07410287,0.38713703){\color[rgb]{0,0,0}\rotatebox{57.00000031}{\makebox(0,0)[lb]{\smash{$\pphi_\infty^1\otimes\thetabar^2$}}}}%
    \put(0.32044567,0.52810981){\color[rgb]{0,0,0}\rotatebox{42.99999925}{\makebox(0,0)[lb]{\smash{$\theta^1\otimes\pphi_\infty^2$}}}}%
    \put(0.30030662,0.79998599){\color[rgb]{0,0,0}\rotatebox{-26.00000032}{\makebox(0,0)[lb]{\smash{$\pphibar_0^1\otimes\theta^2$}}}}%
    \put(0.10898645,0.68922159){\color[rgb]{0,0,0}\rotatebox{-74.99999964}{\makebox(0,0)[lb]{\smash{$\theta^1\otimes\pphibar_0^2$}}}}%
  \end{picture}%
\endgroup%

%% file: Modification4.pdf_tex
\begingroup%
  \makeatletter%
  \providecommand\color[2][]{%
    \errmessage{(Inkscape) Color is used for the text in Inkscape, but the package 'color.sty' is not loaded}%
    \renewcommand\color[2][]{}%
  }%
  \providecommand\transparent[1]{%
    \errmessage{(Inkscape) Transparency is used (non-zero) for the text in Inkscape, but the package 'transparent.sty' is not loaded}%
    \renewcommand\transparent[1]{}%
  }%
  \providecommand\rotatebox[2]{#2}%
  \ifx\svgwidth\undefined%
    \setlength{\unitlength}{553.875bp}%
    \ifx\svgscale\undefined%
      \relax%
    \else%
      \setlength{\unitlength}{\unitlength * \real{\svgscale}}%
    \fi%
  \else%
    \setlength{\unitlength}{\svgwidth}%
  \fi%
  \global\let\svgwidth\undefined%
  \global\let\svgscale\undefined%
  \makeatother%
  \begin{picture}(1,1.1872715)%
    \put(0,0){\includegraphics[width=\unitlength]{Modification4.pdf}}%
    \put(0.2229767,0.8826175){\color[rgb]{0,0,0}\makebox(0,0)[lb]{\smash{$\mu$}}}%
    \put(0.40636425,1.08473979){\color[rgb]{0,0,0}\makebox(0,0)[lb]{\smash{$\lambda$}}}%
    \put(0.53351608,0.79017787){\color[rgb]{0,0,0}\makebox(0,0)[lb]{\smash{W}}}%
    \put(0.39485664,0.79017787){\color[rgb]{0,0,0}\makebox(0,0)[lb]{\smash{A}}}%
    \put(0.43096587,0.79017787){\color[rgb]{0,0,0}\makebox(0,0)[lb]{\smash{X}}}%
    \put(0.4858519,0.79017787){\color[rgb]{0,0,0}\makebox(0,0)[lb]{\smash{D}}}%
    \put(0.4670751,0.8681738){\color[rgb]{0,0,0}\makebox(0,0)[lb]{\smash{E}}}%
    \put(0.42518839,0.8681738){\color[rgb]{0,0,0}\makebox(0,0)[lb]{\smash{Y}}}%
    \put(0.38041295,0.8681738){\color[rgb]{0,0,0}\makebox(0,0)[lb]{\smash{B}}}%
    \put(0.37896858,0.93605916){\color[rgb]{0,0,0}\makebox(0,0)[lb]{\smash{C}}}%
    \put(0.4262201,0.7185784){\color[rgb]{0,0,0}\makebox(0,0)[lb]{\smash{O}}}%
    \put(0.27352962,0.93605916){\color[rgb]{0,0,0}\makebox(0,0)[lb]{\smash{Z}}}%
    \put(0.17474932,0.21938711){\color[rgb]{0,0,0}\makebox(0,0)[lb]{\smash{$\mu$}}}%
    \put(0.11697456,0.10383758){\color[rgb]{0,0,0}\makebox(0,0)[lb]{\smash{$\lambda$}}}%
    \put(0.68801625,0.21089642){\color[rgb]{0,0,0}\makebox(0,0)[lb]{\smash{$\mu'=\mu-\lambda$}}}%
    \put(0.7746784,0.10979057){\color[rgb]{0,0,0}\makebox(0,0)[lb]{\smash{$\lambda$}}}%
    \put(0.21756455,0.48016813){\color[rgb]{0,0,0}\makebox(0,0)[lb]{\smash{Double De-stablization}}}%
  \end{picture}%
\endgroup%

%% file: Modification3.pdf_tex
\begingroup%
  \makeatletter%
  \providecommand\color[2][]{%
    \errmessage{(Inkscape) Color is used for the text in Inkscape, but the package 'color.sty' is not loaded}%
    \renewcommand\color[2][]{}%
  }%
  \providecommand\transparent[1]{%
    \errmessage{(Inkscape) Transparency is used (non-zero) for the text in Inkscape, but the package 'transparent.sty' is not loaded}%
    \renewcommand\transparent[1]{}%
  }%
  \providecommand\rotatebox[2]{#2}%
  \ifx\svgwidth\undefined%
    \setlength{\unitlength}{460.43902924bp}%
    \ifx\svgscale\undefined%
      \relax%
    \else%
      \setlength{\unitlength}{\unitlength * \real{\svgscale}}%
    \fi%
  \else%
    \setlength{\unitlength}{\svgwidth}%
  \fi%
  \global\let\svgwidth\undefined%
  \global\let\svgscale\undefined%
  \makeatother%
  \begin{picture}(1,0.75274905)%
    \put(0,0){\includegraphics[width=\unitlength]{Modification3.pdf}}%
    \put(0.62679002,0.73403788){\color[rgb]{0,0,0}\makebox(0,0)[lb]{\smash{$\mu'$}}}%
    \put(0.55729114,0.73577532){\color[rgb]{0,0,0}\makebox(0,0)[lb]{\smash{$\lambda$}}}%
    \put(0.9073744,0.63007045){\color[rgb]{0,0,0}\makebox(0,0)[lb]{\smash{$\mu'$}}}%
    \put(0.88999968,0.73431878){\color[rgb]{0,0,0}\makebox(0,0)[lb]{\smash{$\lambda$}}}%
    \put(0.20892797,0.73536418){\color[rgb]{0,0,0}\makebox(0,0)[lb]{\smash{$\lambda$}}}%
    \put(0.26973949,0.68324002){\color[rgb]{0,0,0}\makebox(0,0)[lb]{\smash{D}}}%
    \put(0.20892797,0.68324002){\color[rgb]{0,0,0}\makebox(0,0)[lb]{\smash{X}}}%
    \put(0.14892129,0.6819538){\color[rgb]{0,0,0}\makebox(0,0)[lb]{\smash{A}}}%
    \put(0.12205436,0.57030433){\color[rgb]{0,0,0}\makebox(0,0)[lb]{\smash{W}}}%
    \put(0.20559369,0.59450138){\color[rgb]{0,0,0}\makebox(0,0)[lb]{\smash{O}}}%
    \put(0.11107814,0.51576137){\color[rgb]{0,0,0}\makebox(0,0)[lb]{\smash{$\alpha_{g}$}}}%
    \put(0.00043131,0.61374113){\color[rgb]{0,0,0}\makebox(0,0)[lb]{\smash{$\gamma$}}}%
    \put(0.26973949,0.73536418){\color[rgb]{0,0,0}\makebox(0,0)[lb]{\smash{$\mu'$}}}%
    \put(0.16793245,0.15070478){\color[rgb]{0,0,0}\makebox(0,0)[lb]{\smash{$\lambda$}}}%
    \put(0.49908582,0.12898638){\color[rgb]{0,0,0}\makebox(0,0)[lb]{\smash{$\lambda$}}}%
    \put(0.51646054,0.2766715){\color[rgb]{0,0,0}\makebox(0,0)[lb]{\smash{$\mu'$}}}%
    \put(0.81183081,0.18979789){\color[rgb]{0,0,0}\makebox(0,0)[lb]{\smash{$\lambda$}}}%
    \put(0.7249572,0.15504846){\color[rgb]{0,0,0}\makebox(0,0)[lb]{\smash{$\mu'$}}}%
    \put(0.18634083,0.31142094){\color[rgb]{0,0,0}\makebox(0,0)[lb]{\smash{$\mu'$}}}%
  \end{picture}%
\endgroup%

%% file: H-0-1-infty.pdf_tex
\begingroup%
  \makeatletter%
  \providecommand\color[2][]{%
    \errmessage{(Inkscape) Color is used for the text in Inkscape, but the package 'color.sty' is not loaded}%
    \renewcommand\color[2][]{}%
  }%
  \providecommand\transparent[1]{%
    \errmessage{(Inkscape) Transparency is used (non-zero) for the text in Inkscape, but the package 'transparent.sty' is not loaded}%
    \renewcommand\transparent[1]{}%
  }%
  \providecommand\rotatebox[2]{#2}%
  \ifx\svgwidth\undefined%
    \setlength{\unitlength}{444bp}%
    \ifx\svgscale\undefined%
      \relax%
    \else%
      \setlength{\unitlength}{\unitlength * \real{\svgscale}}%
    \fi%
  \else%
    \setlength{\unitlength}{\svgwidth}%
  \fi%
  \global\let\svgwidth\undefined%
  \global\let\svgscale\undefined%
  \makeatother%
  \begin{picture}(1,0.84864865)%
    \put(0,0){\includegraphics[width=\unitlength]{H-0-1-infty.pdf}}%
    \put(0.4963964,0.2081312){\color[rgb]{0,0,0}\makebox(0,0)[lb]{\smash{$r_1$}}}%
    \put(0.45135135,0.32164472){\color[rgb]{0,0,0}\makebox(0,0)[lb]{\smash{$r_2$}}}%
    \put(0.26403517,0.41558607){\color[rgb]{0,0,0}\makebox(0,0)[lb]{\smash{$r_3$}}}%
    \put(0.32653581,0.32163935){\color[rgb]{0,0,0}\makebox(0,0)[lb]{\smash{$p_1$}}}%
    \put(0.32653581,0.41528323){\color[rgb]{0,0,0}\makebox(0,0)[lb]{\smash{$p_2$}}}%
    \put(0.32844703,0.59811179){\color[rgb]{0,0,0}\makebox(0,0)[lb]{\smash{$p_3$}}}%
    \put(0.32844692,0.67646685){\color[rgb]{0,0,0}\makebox(0,0)[lb]{\smash{$p_4$}}}%
    \put(0.33035793,0.81916226){\color[rgb]{0,0,0}\makebox(0,0)[lb]{\smash{$p_n$}}}%
    \put(0.45495502,0.82120609){\color[rgb]{0,0,0}\makebox(0,0)[lb]{\smash{$q_n$}}}%
    \put(0.45495496,0.67659967){\color[rgb]{0,0,0}\makebox(0,0)[lb]{\smash{$q_4$}}}%
    \put(0.45368092,0.59888157){\color[rgb]{0,0,0}\makebox(0,0)[lb]{\smash{$q_3$}}}%
    \put(0.52348287,0.82213469){\color[rgb]{0,0,0}\makebox(0,0)[lb]{\smash{$s_n$}}}%
    \put(0.52539393,0.67752816){\color[rgb]{0,0,0}\makebox(0,0)[lb]{\smash{$s_4$}}}%
    \put(0.52603106,0.60044725){\color[rgb]{0,0,0}\makebox(0,0)[lb]{\smash{$s_3$}}}%
    \put(0.53963964,0.31833793){\color[rgb]{0,0,0}\makebox(0,0)[lb]{\smash{$s_2$}}}%
    \put(0.53963964,0.2081312){\color[rgb]{0,0,0}\makebox(0,0)[lb]{\smash{$s_1$}}}%
    \put(0.26114189,0.05801749){\color[rgb]{0,0,0}\makebox(0,0)[lb]{\smash{$\beta_{g-1}$}}}%
    \put(0.34324324,0.10002309){\color[rgb]{0,0,0}\makebox(0,0)[lb]{\smash{$\lambda_\infty$}}}%
    \put(0.17207207,0.34326634){\color[rgb]{0,0,0}\makebox(0,0)[lb]{\smash{$\lambda_1$}}}%
    \put(0.18468468,0.63335643){\color[rgb]{0,0,0}\makebox(0,0)[lb]{\smash{$\lambda_0$}}}%
    \put(0.2152771,0.54178412){\color[rgb]{0,0,0}\makebox(0,0)[lb]{\smash{$u$}}}%
    \put(0.25495495,0.56128436){\color[rgb]{0,0,0}\makebox(0,0)[lb]{\smash{$\vbar$}}}%
    \put(0.27297297,0.50002309){\color[rgb]{0,0,0}\makebox(0,0)[lb]{\smash{$w$}}}%
    \put(0.12162162,0.50002309){\color[rgb]{0,0,0}\makebox(0,0)[lb]{\smash{$v$}}}%
    \put(0.24060218,0.45028582){\color[rgb]{0,0,0}\makebox(0,0)[lb]{\smash{$\ubar$}}}%
    \put(0.27297297,0.59732039){\color[rgb]{0,0,0}\makebox(0,0)[lb]{\smash{$t_1$}}}%
    \put(0.71261261,0.2081312){\color[rgb]{0,0,0}\makebox(0,0)[lb]{\smash{$t_0$}}}%
    \put(0.73963964,0.31804111){\color[rgb]{0,0,0}\makebox(0,0)[lb]{\smash{$t_2$}}}%
    \put(0.78828829,0.59732039){\color[rgb]{0,0,0}\makebox(0,0)[lb]{\smash{$t_3$}}}%
    \put(0.79009009,0.67659967){\color[rgb]{0,0,0}\makebox(0,0)[lb]{\smash{$t_4$}}}%
    \put(0.79009009,0.81894201){\color[rgb]{0,0,0}\makebox(0,0)[lb]{\smash{$t_n$}}}%
    \put(0.66576577,0.2081312){\color[rgb]{0,0,0}\makebox(0,0)[lb]{\smash{$u_1$}}}%
    \put(0.66576577,0.31804111){\color[rgb]{0,0,0}\makebox(0,0)[lb]{\smash{$u_2$}}}%
    \put(0.66306306,0.59822129){\color[rgb]{0,0,0}\makebox(0,0)[lb]{\smash{$u_3$}}}%
    \put(0.66756757,0.67659967){\color[rgb]{0,0,0}\makebox(0,0)[lb]{\smash{$u_4$}}}%
    \put(0.66936937,0.81804111){\color[rgb]{0,0,0}\makebox(0,0)[lb]{\smash{$u_n$}}}%
    \put(0.86846847,0.81894201){\color[rgb]{0,0,0}\makebox(0,0)[lb]{\smash{$v_n$}}}%
    \put(0.86756757,0.67659967){\color[rgb]{0,0,0}\makebox(0,0)[lb]{\smash{$v_4$}}}%
    \put(0.86756757,0.59732039){\color[rgb]{0,0,0}\makebox(0,0)[lb]{\smash{$v_3$}}}%
    \put(0.86307695,0.31885602){\color[rgb]{0,0,0}\makebox(0,0)[lb]{\smash{$v_2$}}}%
    \put(0.86668056,0.09092809){\color[rgb]{0,0,0}\makebox(0,0)[lb]{\smash{$v_1$}}}%
    \put(0.17207207,0.55227535){\color[rgb]{0,0,0}\makebox(0,0)[lb]{\smash{$\wbar$}}}%
    \put(0.61520606,0.8212062){\color[rgb]{0,0,0}\makebox(0,0)[lb]{\smash{$S_n$}}}%
    \put(0.61711712,0.67659967){\color[rgb]{0,0,0}\makebox(0,0)[lb]{\smash{$S_4$}}}%
    \put(0.61775425,0.59951876){\color[rgb]{0,0,0}\makebox(0,0)[lb]{\smash{$S_3$}}}%
    \put(0.61531532,0.31833793){\color[rgb]{0,0,0}\makebox(0,0)[lb]{\smash{$S_2$}}}%
    \put(0.61531532,0.2081312){\color[rgb]{0,0,0}\makebox(0,0)[lb]{\smash{$S_1$}}}%
    \put(0.38235518,0.81964491){\color[rgb]{0,0,0}\makebox(0,0)[lb]{\smash{$Q_n$}}}%
    \put(0.38235512,0.67503849){\color[rgb]{0,0,0}\makebox(0,0)[lb]{\smash{$Q_4$}}}%
    \put(0.38108108,0.59732039){\color[rgb]{0,0,0}\makebox(0,0)[lb]{\smash{$Q_3$}}}%
    \put(0.38108108,0.409933){\color[rgb]{0,0,0}\makebox(0,0)[lb]{\smash{$R_3$}}}%
    \put(0.37927928,0.32164472){\color[rgb]{0,0,0}\makebox(0,0)[lb]{\smash{$R_2$}}}%
    \put(0.44234234,0.2081312){\color[rgb]{0,0,0}\makebox(0,0)[lb]{\smash{$R_1$}}}%
    \put(0.58468468,0.09641949){\color[rgb]{0,0,0}\makebox(0,0)[lb]{\smash{$\lambda$}}}%
  \end{picture}%
\endgroup%

%% file: pentagon.pdf_tex
\begingroup%
  \makeatletter%
  \providecommand\color[2][]{%
    \errmessage{(Inkscape) Color is used for the text in Inkscape, but the package 'color.sty' is not loaded}%
    \renewcommand\color[2][]{}%
  }%
  \providecommand\transparent[1]{%
    \errmessage{(Inkscape) Transparency is used (non-zero) for the text in Inkscape, but the package 'transparent.sty' is not loaded}%
    \renewcommand\transparent[1]{}%
  }%
  \providecommand\rotatebox[2]{#2}%
  \ifx\svgwidth\undefined%
    \setlength{\unitlength}{438.65003174bp}%
    \ifx\svgscale\undefined%
      \relax%
    \else%
      \setlength{\unitlength}{\unitlength * \real{\svgscale}}%
    \fi%
  \else%
    \setlength{\unitlength}{\svgwidth}%
  \fi%
  \global\let\svgwidth\undefined%
  \global\let\svgscale\undefined%
  \makeatother%
  \begin{picture}(1,0.3901993)%
    \put(0,0){\includegraphics[width=\unitlength]{pentagon.pdf}}%
    \put(0.40969856,0.36210068){\color[rgb]{0,0,0}\makebox(0,0)[lb]{\smash{$L_1$}}}%
    \put(0.36475546,0.00646415){\color[rgb]{0,0,0}\makebox(0,0)[lb]{\smash{$\gamma_1$}}}%
    \put(0.45203627,0.00516143){\color[rgb]{0,0,0}\makebox(0,0)[lb]{\smash{$\beta_1$}}}%
    \put(0.43249578,0.27481995){\color[rgb]{0,0,0}\makebox(0,0)[lb]{\smash{$Q_3$}}}%
    \put(0.47548481,0.32757915){\color[rgb]{0,0,0}\makebox(0,0)[lb]{\smash{$Y_1$}}}%
    \put(0.4754114,0.18967029){\color[rgb]{0,0,0}\makebox(0,0)[lb]{\smash{$X_1=Q_4$}}}%
    \put(0.47678753,0.06964494){\color[rgb]{0,0,0}\makebox(0,0)[lb]{\smash{$Q_8$}}}%
    \put(0.40579042,0.12305566){\color[rgb]{0,0,0}\makebox(0,0)[lb]{\smash{$R_1$}}}%
    \put(0.50284148,0.23573903){\color[rgb]{0,0,0}\makebox(0,0)[lb]{\smash{$T_1$}}}%
    \put(0.50349283,0.280682){\color[rgb]{0,0,0}\makebox(0,0)[lb]{\smash{$S_1$}}}%
    \put(0.5895555,0.27259611){\color[rgb]{0,0,0}\makebox(0,0)[lb]{\smash{$Q_7$}}}%
    \put(0.60797049,0.24094966){\color[rgb]{0,0,0}\makebox(0,0)[lb]{\smash{$\alpha_1$}}}%
    \put(0.60085198,0.0768097){\color[rgb]{0,0,0}\makebox(0,0)[lb]{\smash{$\alpha_2$}}}%
    \put(0.00112875,0.28099813){\color[rgb]{0,0,0}\makebox(0,0)[lb]{\smash{$\gamma_g$}}}%
    \put(0.00297106,0.05536162){\color[rgb]{0,0,0}\makebox(0,0)[lb]{\smash{$\beta_g$}}}%
    \put(0.07099702,0.16343913){\color[rgb]{0,0,0}\makebox(0,0)[lb]{\smash{$P=Q_1$}}}%
    \put(0.21820191,0.27742532){\color[rgb]{0,0,0}\makebox(0,0)[lb]{\smash{$Q_2$}}}%
    \put(0.22276138,0.04163687){\color[rgb]{0,0,0}\makebox(0,0)[lb]{\smash{$Q_6$}}}%
    \put(0.34130691,0.07355308){\color[rgb]{0,0,0}\makebox(0,0)[lb]{\smash{$Q_5$}}}%
    \put(0.26314502,0.15757708){\color[rgb]{0,0,0}\makebox(0,0)[lb]{\smash{$D_1$}}}%
    \put(0.89476053,0.24427986){\color[rgb]{0,0,0}\makebox(0,0)[lb]{\smash{$Y_i$}}}%
    \put(0.89356958,0.11899735){\color[rgb]{0,0,0}\makebox(0,0)[lb]{\smash{$X_i$}}}%
    \put(0.91050458,0.18608618){\color[rgb]{0,0,0}\makebox(0,0)[lb]{\smash{$T_i$}}}%
    \put(0.82517787,0.26620231){\color[rgb]{0,0,0}\makebox(0,0)[lb]{\smash{$L_i$}}}%
    \put(0.82387527,0.10987836){\color[rgb]{0,0,0}\makebox(0,0)[lb]{\smash{$R_i$}}}%
    \put(0.7958672,0.19064562){\color[rgb]{0,0,0}\makebox(0,0)[lb]{\smash{$D_i$}}}%
    \put(0.89213596,0.35869384){\color[rgb]{0,0,0}\makebox(0,0)[lb]{\smash{$\beta_i$}}}%
    \put(0.79786763,0.35972684){\color[rgb]{0,0,0}\makebox(0,0)[lb]{\smash{$\gamma_i$}}}%
    \put(0.33088531,0.01102359){\color[rgb]{0,0,0}\makebox(0,0)[lb]{\smash{$\ubar$}}}%
    \put(0.16283725,0.09960714){\color[rgb]{0,0,0}\makebox(0,0)[lb]{\smash{$\ubar$}}}%
    \put(0.55364671,0.13738529){\color[rgb]{0,0,0}\makebox(0,0)[lb]{\smash{$\ubar$}}}%
    \put(0.16674536,0.22206059){\color[rgb]{0,0,0}\makebox(0,0)[lb]{\smash{$\wbar$}}}%
    \put(0.01693508,0.16213652){\color[rgb]{0,0,0}\makebox(0,0)[lb]{\smash{$\vbar$}}}%
    \put(0.35172846,0.28459031){\color[rgb]{0,0,0}\makebox(0,0)[lb]{\smash{$\ubar$}}}%
  \end{picture}%
\endgroup%

%% file: pentagon-3.pdf_tex
\begingroup%
  \makeatletter%
  \providecommand\color[2][]{%
    \errmessage{(Inkscape) Color is used for the text in Inkscape, but the package 'color.sty' is not loaded}%
    \renewcommand\color[2][]{}%
  }%
  \providecommand\transparent[1]{%
    \errmessage{(Inkscape) Transparency is used (non-zero) for the text in Inkscape, but the package 'transparent.sty' is not loaded}%
    \renewcommand\transparent[1]{}%
  }%
  \providecommand\rotatebox[2]{#2}%
  \ifx\svgwidth\undefined%
    \setlength{\unitlength}{291.84233398bp}%
    \ifx\svgscale\undefined%
      \relax%
    \else%
      \setlength{\unitlength}{\unitlength * \real{\svgscale}}%
    \fi%
  \else%
    \setlength{\unitlength}{\svgwidth}%
  \fi%
  \global\let\svgwidth\undefined%
  \global\let\svgscale\undefined%
  \makeatother%
  \begin{picture}(1,1.14953912)%
    \put(0,0){\includegraphics[width=\unitlength]{pentagon-3.pdf}}%
    \put(0.73359078,0.71246096){\color[rgb]{0,0,0}\makebox(0,0)[lb]{\smash{$X_1$}}}%
    \put(0.18954554,0.67483697){\color[rgb]{0,0,0}\makebox(0,0)[lb]{\smash{$Q$}}}%
    \put(0.40903869,0.6645551){\color[rgb]{0,0,0}\makebox(0,0)[lb]{\smash{$D_1$}}}%
    \put(0.62790084,0.61099317){\color[rgb]{0,0,0}\makebox(0,0)[lb]{\smash{$R_1$}}}%
    \put(0.67694758,0.95199919){\color[rgb]{0,0,0}\makebox(0,0)[lb]{\smash{$p_1$}}}%
    \put(0.68110122,0.83293056){\color[rgb]{0,0,0}\makebox(0,0)[lb]{\smash{$p_2$}}}%
    \put(0.90123962,0.95338368){\color[rgb]{0,0,0}\makebox(0,0)[lb]{\smash{$t_1$}}}%
    \put(0.9040086,0.82877701){\color[rgb]{0,0,0}\makebox(0,0)[lb]{\smash{$t_2$}}}%
    \put(0.43188787,0.95338368){\color[rgb]{0,0,0}\makebox(0,0)[lb]{\smash{$r_2$}}}%
    \put(0.357124,0.83293065){\color[rgb]{0,0,0}\makebox(0,0)[lb]{\smash{$r_3$}}}%
    \put(0.74063545,1.03922381){\color[rgb]{0,0,0}\makebox(0,0)[lb]{\smash{$Y_1$}}}%
    \put(0.76417228,0.77754988){\color[rgb]{0,0,0}\makebox(0,0)[lb]{\smash{$T_1$}}}%
    \put(0.76971041,0.97276696){\color[rgb]{0,0,0}\makebox(0,0)[lb]{\smash{$S_1$}}}%
    \put(0.61464437,1.08214391){\color[rgb]{0,0,0}\makebox(0,0)[lb]{\smash{$L_1$}}}%
    \put(0.07191323,0.80524018){\color[rgb]{0,0,0}\makebox(0,0)[lb]{\smash{$\lambda_\infty$}}}%
    \put(0.06499062,0.54771989){\color[rgb]{0,0,0}\makebox(0,0)[lb]{\smash{ $\lambda_1$}}}%
    \put(0.7364819,0.28189233){\color[rgb]{0,0,0}\makebox(0,0)[lb]{\smash{$p_5$}}}%
    \put(0.73648198,0.37188584){\color[rgb]{0,0,0}\makebox(0,0)[lb]{\smash{$p_4$}}}%
    \put(0.73509741,0.48126288){\color[rgb]{0,0,0}\makebox(0,0)[lb]{\smash{$p_3$}}}%
    \put(0.73371292,0.10190481){\color[rgb]{0,0,0}\makebox(0,0)[lb]{\smash{$p_n$}}}%
    \put(0.59110754,0.10052032){\color[rgb]{0,0,0}\makebox(0,0)[lb]{\smash{$t_n$}}}%
    \put(0.5952611,0.47710924){\color[rgb]{0,0,0}\makebox(0,0)[lb]{\smash{$t_3$}}}%
    \put(0.437426,0.37327033){\color[rgb]{0,0,0}\makebox(0,0)[lb]{\smash{$q_4$}}}%
    \put(0.43604155,0.28604563){\color[rgb]{0,0,0}\makebox(0,0)[lb]{\smash{$q_5$}}}%
    \put(0.59249212,0.36911686){\color[rgb]{0,0,0}\makebox(0,0)[lb]{\smash{$t_4$}}}%
    \put(0.59387661,0.28050767){\color[rgb]{0,0,0}\makebox(0,0)[lb]{\smash{$t_5$}}}%
    \put(0.4401951,0.10605828){\color[rgb]{0,0,0}\makebox(0,0)[lb]{\smash{$q_n$}}}%
    \put(0.40558219,0.4784939){\color[rgb]{0,0,0}\makebox(0,0)[lb]{\smash{$q_3$}}}%
    \put(0.71432956,0.00775783){\color[rgb]{0,0,0}\makebox(0,0)[lb]{\smash{$\beta_{1}$}}}%
    \put(0.57033978,0.01052664){\color[rgb]{0,0,0}\makebox(0,0)[lb]{\smash{$\gamma_{1}$}}}%
    \put(0.25882317,0.7346297){\color[rgb]{0,0,0}\makebox(0,0)[lb]{\smash{$w$}}}%
    \put(0.11483329,0.64740516){\color[rgb]{0,0,0}\makebox(0,0)[lb]{\smash{$u$}}}%
    \put(0.25743868,0.5504887){\color[rgb]{0,0,0}\makebox(0,0)[lb]{\smash{$v$}}}%
    \put(0.25328511,0.84816012){\color[rgb]{0,0,0}\makebox(0,0)[lb]{\smash{$v$}}}%
    \put(0.53572682,0.99353456){\color[rgb]{0,0,0}\makebox(0,0)[lb]{\smash{$v$}}}%
    \put(0.80570781,0.59479338){\color[rgb]{0,0,0}\makebox(0,0)[lb]{\smash{$w$}}}%
    \put(0.94139066,0.85369816){\color[rgb]{0,0,0}\makebox(0,0)[lb]{\smash{$v$}}}%
  \end{picture}%
\endgroup%

%% file: Quadrangles.pdf_tex
\begingroup%
  \makeatletter%
  \providecommand\color[2][]{%
    \errmessage{(Inkscape) Color is used for the text in Inkscape, but the package 'color.sty' is not loaded}%
    \renewcommand\color[2][]{}%
  }%
  \providecommand\transparent[1]{%
    \errmessage{(Inkscape) Transparency is used (non-zero) for the text in Inkscape, but the package 'transparent.sty' is not loaded}%
    \renewcommand\transparent[1]{}%
  }%
  \providecommand\rotatebox[2]{#2}%
  \ifx\svgwidth\undefined%
    \setlength{\unitlength}{534.54594727bp}%
    \ifx\svgscale\undefined%
      \relax%
    \else%
      \setlength{\unitlength}{\unitlength * \real{\svgscale}}%
    \fi%
  \else%
    \setlength{\unitlength}{\svgwidth}%
  \fi%
  \global\let\svgwidth\undefined%
  \global\let\svgscale\undefined%
  \makeatother%
  \begin{picture}(1,0.30035649)%
    \put(0,0){\includegraphics[width=\unitlength]{Quadrangles.pdf}}%
    \put(0.56441584,0.00423549){\color[rgb]{0,0,0}\makebox(0,0)[lb]{\smash{$\lambda_0'$}}}%
    \put(0.60988814,0.06962045){\color[rgb]{0,0,0}\makebox(0,0)[lb]{\smash{$\lambda_1'$}}}%
    \put(0.50621187,0.06924247){\color[rgb]{0,0,0}\makebox(0,0)[lb]{\smash{$\lambda_1$}}}%
    \put(0.85014443,0.13387156){\color[rgb]{0,0,0}\makebox(0,0)[lb]{\smash{$\lambda_\infty$}}}%
    \put(0.37444029,0.13814751){\color[rgb]{0,0,0}\makebox(0,0)[lb]{\smash{$\lambda_\infty$}}}%
    \put(0.08443573,0.00637351){\color[rgb]{0,0,0}\makebox(0,0)[lb]{\smash{$\lambda_0'$}}}%
    \put(0.13311501,0.07175848){\color[rgb]{0,0,0}\makebox(0,0)[lb]{\smash{$\lambda_1'$}}}%
    \put(0.02943875,0.07138049){\color[rgb]{0,0,0}\makebox(0,0)[lb]{\smash{$\lambda_1$}}}%
    \put(0.69209603,0.16945685){\color[rgb]{0,0,0}\makebox(0,0)[lb]{\smash{$\Theta_{11}$}}}%
    \put(0.59802415,0.15342186){\color[rgb]{0,0,0}\makebox(0,0)[lb]{\smash{$\Theta_{1\infty}$}}}%
    \put(0.57878227,0.043315){\color[rgb]{0,0,0}\makebox(0,0)[lb]{\smash{$\Theta_{01}$}}}%
    \put(0.21746091,0.17052586){\color[rgb]{0,0,0}\makebox(0,0)[lb]{\smash{$\Theta_{11}$}}}%
    \put(0.12018206,0.15876692){\color[rgb]{0,0,0}\makebox(0,0)[lb]{\smash{$\Theta_{1\infty}$}}}%
    \put(0.04535222,0.15876683){\color[rgb]{0,0,0}\makebox(0,0)[lb]{\smash{$\Theta_{0\infty}$}}}%
    \put(0.10414715,0.04224599){\color[rgb]{0,0,0}\makebox(0,0)[lb]{\smash{$\Theta_{01}$}}}%
    \put(0.52105636,0.15449087){\color[rgb]{0,0,0}\makebox(0,0)[lb]{\smash{$\Theta_{0\infty}$}}}%
  \end{picture}%
\endgroup%

%% file: local.pdf_tex
\begingroup%
  \makeatletter%
  \providecommand\color[2][]{%
    \errmessage{(Inkscape) Color is used for the text in Inkscape, but the package 'color.sty' is not loaded}%
    \renewcommand\color[2][]{}%
  }%
  \providecommand\transparent[1]{%
    \errmessage{(Inkscape) Transparency is used (non-zero) for the text in Inkscape, but the package 'transparent.sty' is not loaded}%
    \renewcommand\transparent[1]{}%
  }%
  \providecommand\rotatebox[2]{#2}%
  \ifx\svgwidth\undefined%
    \setlength{\unitlength}{405.67370944bp}%
    \ifx\svgscale\undefined%
      \relax%
    \else%
      \setlength{\unitlength}{\unitlength * \real{\svgscale}}%
    \fi%
  \else%
    \setlength{\unitlength}{\svgwidth}%
  \fi%
  \global\let\svgwidth\undefined%
  \global\let\svgscale\undefined%
  \makeatother%
  \begin{picture}(1,0.35749174)%
    \put(0,0){\includegraphics[width=\unitlength]{local.pdf}}%
    \put(0.74744558,0.06017214){\color[rgb]{0,0,0}\makebox(0,0)[lb]{\smash{$\lambda_1\#\mu_2$}}}%
    \put(0.74603697,0.12778452){\color[rgb]{0,0,0}\makebox(0,0)[lb]{\smash{$\mu_1\#\lambda_2$}}}%
    \put(0.74181115,0.27991246){\color[rgb]{0,0,0}\makebox(0,0)[lb]{\smash{$\mu_1\#\lambda_2$}}}%
    \put(0.74321982,0.1968055){\color[rgb]{0,0,0}\makebox(0,0)[lb]{\smash{$\lambda_1\#\mu_2$}}}%
    \put(0.28542753,0.29118109){\color[rgb]{0,0,0}\makebox(0,0)[lb]{\smash{$W_1$}}}%
    \put(0.29528767,0.21370863){\color[rgb]{0,0,0}\makebox(0,0)[lb]{\smash{$A_1$}}}%
    \put(0.24035259,0.21370863){\color[rgb]{0,0,0}\makebox(0,0)[lb]{\smash{$B_1$}}}%
    \put(0.1882347,0.21511724){\color[rgb]{0,0,0}\makebox(0,0)[lb]{\smash{$C_1$}}}%
    \put(0.29669625,0.13764466){\color[rgb]{0,0,0}\makebox(0,0)[lb]{\smash{$X_1$}}}%
    \put(0.24457839,0.13905327){\color[rgb]{0,0,0}\makebox(0,0)[lb]{\smash{$Y_1$}}}%
    \put(0.18964328,0.13905339){\color[rgb]{0,0,0}\makebox(0,0)[lb]{\smash{$Z_1$}}}%
    \put(0.29528767,0.07566671){\color[rgb]{0,0,0}\makebox(0,0)[lb]{\smash{$D_1$}}}%
    \put(0.23190106,0.07566671){\color[rgb]{0,0,0}\makebox(0,0)[lb]{\smash{$E_1$}}}%
    \put(0.43051246,0.07707508){\color[rgb]{0,0,0}\makebox(0,0)[lb]{\smash{$W_2$}}}%
    \put(0.42910385,0.13764466){\color[rgb]{0,0,0}\makebox(0,0)[lb]{\smash{$A_2$}}}%
    \put(0.49671623,0.13764466){\color[rgb]{0,0,0}\makebox(0,0)[lb]{\smash{$B_2$}}}%
    \put(0.57700603,0.13905327){\color[rgb]{0,0,0}\makebox(0,0)[lb]{\smash{$C_2$}}}%
    \put(0.43192113,0.21370857){\color[rgb]{0,0,0}\makebox(0,0)[lb]{\smash{$X_2$}}}%
    \put(0.49530768,0.21511724){\color[rgb]{0,0,0}\makebox(0,0)[lb]{\smash{$Y_2$}}}%
    \put(0.57700603,0.21370863){\color[rgb]{0,0,0}\makebox(0,0)[lb]{\smash{$Z_2$}}}%
    \put(0.43262531,0.29329406){\color[rgb]{0,0,0}\makebox(0,0)[lb]{\smash{$D_2$}}}%
    \put(0.5023506,0.2939983){\color[rgb]{0,0,0}\makebox(0,0)[lb]{\smash{$E_2$}}}%
  \end{picture}%
\endgroup%